\documentstyle[10pt, epic, eepic, draft, amscd, amssymb]{article}
\newtheorem{theorem}{Theorem}
\newtheorem{lemma}{Lemma}
\newtheorem{corollary}{Corollary}
\newtheorem{definition}{Definition}
\newcommand{\real}{{\bf R}}
\newcommand{\lb}{\lfloor}
\newcommand{\rb}{\rfloor}
\newcommand{\rbe}{\rfloor_\epsilon}
\newcommand{\blb}{\lfloor}
\newcommand{\brb}{\rfloor}
\newcommand{\qed}{$\Box$}
\begin{document}

\title{Proof of the Riemannian Penrose Conjecture
            Using the Positive Mass Theorem      }

\author{             Hubert L. Bray
\thanks{
Mathematics Department, 2-179,
Massachusetts Institute of Technology,
77 Massachusetts Avenue,
Cambridge, MA  02139,
bray@math.mit.edu. 
Research supported in part by NSF grant \#DMS-9706006.}}

\date{   Preliminary Version: November 22, 1999   }

\maketitle


\begin{abstract} 
We prove the Riemannian Penrose conjecture, an important case of a 
conjecture made by Roger Penrose in 1973, by defining a new flow of metrics. 
This flow of metrics stays inside the class of 
asymptotically flat Riemannian 3-manifolds with nonnegative scalar curvature
which contain minimal spheres.
In particular, if we consider a 
Riemannian 3-manifold as a totally geodesic submanifold of a space-time 
in the context of general relativity, then outermost minimal spheres 
with total area $A$
correspond to apparent horizons of black holes
contributing a mass $\sqrt{A/16\pi}$, scalar curvature
corresponds to local energy density at each point, and the rate at which
the metric becomes flat at infinity corresponds to total mass.  The 
Riemannian Penrose conjecture then states
that the total mass of an
asymptotically flat 3-manifold with nonnegative scalar curvature 
is greater than or equal to the mass contributed by the black holes.

The flow of metrics we define continuously evolves the original 3-metric 
to a Schwarzschild 3-metric, which represents a spherically
symmetric black hole in vacuum.  We define the flow such that the
area of the minimal spheres (which flow outward) and hence 
the mass contributed by the black holes
in each of the metrics in the flow is constant, and then use the positive
mass theorem to show that the total mass of the
metrics is nonincreasing.  Then since the total mass equals the mass of the
black holes in a Schwarzschild metric, the Riemannian Penrose conjecture 
follows. 

This result improves upon the beautiful work of Huisken and Ilmanen
\cite{HI}, who used inverse mean curvature flows of surfaces to show  
that the total mass is at least the mass contributed
by the largest black hole. 
\end{abstract}

In section \ref{sec:introduction} and 
\ref{sec:definitions}, we motivate the problem, discuss important quantities
like total mass and horizons of black holes, and state the positive mass 
theorem and the Penrose conjecture for Riemannian 3-manifolds.
In section \ref{sec:overview}, we 
give the
proof of the Riemannian Penrose conjecture, with the supporting
arguments given in sections \ref{sec:existence} through 
\ref{sec:generalization}.  
In section \ref{sec:qlm}, we apply the techniques used in this paper to 
define several new quasi-local mass functions which have good monotonicity
properties.  Finally, in section \ref{sec:acknowledgments} we end with a 
brief discussion of some of the interesting problems which still remain open,
and the author thanks the many people who have made important 
contributions to the ideas in this paper.

\section{Introduction}\label{sec:introduction}

General relativity is a theory of gravity
which asserts that matter causes the four dimensional space-time 
to be curved, and that our perception of 
gravity is a consequence
of this curvature.  Let $(N^4,\bar{g})$ 
be the space-time manifold with
metric $\bar{g}$ of signature $(-+++)$.  Then the central formula of 
General Relativity is Einstein's equation,
\begin{equation}\label{Einstein}
G = 8\pi T,
\end{equation}
where $T$ is the stress-energy tensor,  
$G=Ric(\bar{g}) - \frac12 R(\bar{g})\cdot \bar{g}$ is the Einstein curvature
tensor, $Ric(\bar{g})$ is the Ricci curvature tensor, and 
$R(\bar{g})$ is the scalar curvature of $\bar{g}$.  
The beauty of General Relativity is that this simple formula 
explains gravity more accurately than Newtonian physics and is 
entirely consistent with large scale observations.

However, the nature of the behavior of matter in General Relativity
is still not well understood.  It is not even well 
understood how to define how much energy and 
momentum exist in a given region, except in special cases.  There
does exist a well defined notion of local energy and momentum 
density which is simply 
given by the stress-energy tensor which,
by equation \ref{Einstein}, can be computed in terms of the 
curvature of $N^4$.  Also, if we assume that the matter
of the space-time manifold $N^4$ is concentrated in 
some central region of the universe, then
$N^4$ becomes flatter as we get farther away from this
central region.  If the curvature of $N^4$ decays quickly enough,
then $N^4$ is said to be asymptotically flat (definition 
\ref{def:asymptotically_flat} in section \ref{sec:generalization}), 
and with these
assumptions it is then possible to define the total mass of the 
space-time $N^4$.  Interestingly enough, though, the definition
of local energy-momentum density, which involves curvature
terms of $N^4$, bears no obvious  
resemblance to the definition of the total mass of $N^4$, which
is a parameter related to how fast the metric becomes flat at infinity.

The Penrose conjecture (\cite{P}, \cite{JW}, \cite{HI}) 
and the positive mass theorem 
(\cite{SY3}, \cite{SY4}, \cite{SY5}, \cite{SY6}, \cite{Wi}) can both be
thought of as basic attempts at understanding the 
relationship between the local energy density of a space-time $N^4$
and the total mass of $N^4$.
In physical terms, the 
positive mass theorem states 
that an isolated gravitational system with 
nonnegative local energy density must have nonnegative total energy.  The
idea is that nonnegative energy densities must ``add up'' 
to something nonnegative.  The
Penrose conjecture, on the other hand, states that if an isolated
gravitational system with nonnegative local energy density contains 
black holes contributing a mass $m$, then the total energy of 
the system must be at least $m$. 

Important cases of the positive mass theorem 
and the Penrose conjecture can be translated
into statements about complete, asymptotically flat Riemannian  
$3$-manifolds $(M^3,g)$ 
with nonnegative scalar curvature. 
If we consider $(M^3,g, h)$ as a space-like hypersurface of 
$(N^4,\bar{g})$ with metric $g_{ij}$ and 
second fundamental form $h_{ij}$ in $N^4$,
then equation \ref{Einstein} implies that
\begin{equation}\label{eqn:c1}
\mu = G^0_0 =  \frac{1}{16\pi} 
      [R - \sum_{i,j} h^{ij}h_{ij} + (\sum_i h_i^i)^2],
\end{equation}
\begin{equation}\label{eqn:c2}
J^i = G_0^i = \frac{1}{8\pi}\sum_j \nabla_j[h^{ij} - (\sum_k h_k^k)g^{ij}],
\end{equation}
where $R$ is the scalar curvature of the metric $g$, $\mu$ is the
local energy density, and $J^i$ is the local current density.
The assumption of nonnegative energy density everywhere
in $N^4$, called the dominant energy condition, implies that we must have
\begin{equation}\label{eqn:c3}
\mu \ge \left(\sum_i J^i J_i \right)^\frac12
\end{equation}
at all points on $M^3$.  Equations \ref{eqn:c1},
\ref{eqn:c2}, and \ref{eqn:c3}
are called the constraint equations for $(M^3,g,h)$ in $(N^4, \bar{g})$.
Thus we see that 
if we restrict our attention 
to $3$-manifolds which have zero second fundamental form $h$ in $N^3$, 
the constraint equations are equivalent to the condition 
that the Riemannian manifold $(M^3,g)$ has nonnegative scalar 
curvature everywhere. 

An asymptotically flat $3$-manifold is a Riemannian  
manifold $(M^3,g)$ which, outside a compact set, is the disjoint union 
of one or more regions (called ends) diffeomorphic to 
($\real^3 \backslash B_1(0), \delta)$, where the 
metric $g$ in each of these $\real^3$ coordinate charts 
approaches the standard metric $\delta$ on $\real^3$ at
infinity (with certain asymptotic decay conditions
 - see definition \ref{def:asymptotically_flat}, \cite{SY5}, \cite{Ba2}).  
The positive mass theorem and the Penrose conjecture are both statements
which refer to a particular {\bf chosen end} of $(M^3,g)$.  The total mass
of $(M^3,g)$, also called the ADM mass \cite{ADM}, is then a parameter
related to how fast this chosen end of $(M^3,g)$ becomes flat at infinity.
The usual definition of the total mass is
given in equation \ref{eqn:ADM_mass} in section \ref{sec:generalization}.
In addition, we give an alternate  
definition of total mass in the next section.

The positive mass theorem was first proved by Schoen and Yau \cite{SY3}
in 1979 using minimal surfaces.  The Riemannian positive mass theorem 
is a special case of the positive mass theorem which comes from considering 
the space-like hypersurfaces which have zero second fundamental form in 
the spacetime.

\vspace{.1in}\noindent
{\bf The Riemannian Positive Mass Theorem}
\newline {\it
Let $(M^3,g)$ be a complete, smooth, asymptotically flat 3-manifold with 
nonnegative scalar curvature and total mass $m$.
Then 
\begin{equation}
m \ge 0 , 
\end{equation}
with equality if and only if
$(M^3,g)$ is isometric to $\real^3$ with the standard flat metric.}
\vspace{.1in}

Apparent horizons of black holes in $N^4$ 
correspond to outermost minimal surfaces 
of $M^3$ if we assume $M^3$ has zero second fundamental form
in $N^4$.  A minimal surface is a surface which has zero mean curvature (and 
hence is a critical point for the area functional).
An 
outermost minimal surface is a minimal surface which is not contained entirely inside
another minimal surface.  Again, there is a chosen end of $M^3$, and 
``contained entirely inside'' is defined with respect to this end.
Interestingly, it follows from a stability argument
\cite{SY1} that outermost minimal surfaces are always spheres.
There could be more than one outermost sphere, with each minimal sphere
corresponding to a different black hole, and we note 
that outermost minimal spheres never intersect.

As an example, consider the Schwarzschild manifolds 
$(\real^3 \backslash  \{0\}, s)$ where $s_{ij} = (1 + m/2r)^4 \; \delta_{ij}$ and 
$m$ is a positive constant and equals
the total mass of the manifold.  This manifold has zero scalar curvature
everywhere, is spherically symmetric, and it can be checked that it 
has an outermost minimal sphere at $r=m/2$.  

We define the horizon of $(M^3,g)$ to be 
the union of all of the outermost minimal spheres in $M^3$, so that 
the horizon of a manifold can have multiple connected components.
We note that it is usually more common to call each outermost minimal sphere
a horizon, so that their union is referred to as ``horizons'', but it turns
out to be mathematically more convenient for our purposes 
to refer to the union of all of the 
outermost minimal spheres as one object, which we will call the horizon
of $(M^3,g)$.
    
There is a very convincing (but not rigorous)
physical motivation to define the mass that a collection of black holes 
contributes to be
$\sqrt{\frac{A}{16\pi}}$, where $A$ is the total surface area of the 
horizon of $(M^3,g)$. 
Then the physical statement that a system
with nonnegative energy density containing 
black holes contributing a mass $m$ must 
have total mass
at least $m$ can be translated into the following geometric statement 
(\cite{P}, \cite{JW}), 
the proof of which is the object of this paper.

\vspace{.1in}\noindent
{\bf The Riemannian Penrose Conjecture} \newline {\it
Let $(M^3,g)$ be a complete, smooth, asymptotically flat 3-manifold
with nonnegative scalar curvature and total mass $m$ whose  
outermost minimal spheres have total surface area $A$.
Then 
\begin{equation}\label{eqn:RPI}
m \ge \sqrt{\frac{A}{16\pi}}, 
\end{equation}
with equality if and only if 
$(M^3,g)$ is isometric to the Schwarzschild metric  
$(\real^3 \backslash \{0\}, s)$ of mass $m$ outside their respective horizons.
}\vspace{.1in}

The overview of the proof of this result is given in section 
\ref{sec:overview}. 
The basic idea of the approach to the problem is to flow the original
metric continuously to a Schwarzschild metric (outside the horizon).  The
particular flow we define has the important  
property that the area of the horizon stays constant while
the total mass of the manifold is non-increasing.  Then since the Schwarzschild
metric gives equality in the Penrose inequality, the inequality
follows for the original metric.

The first breakthrough on the Riemannian Penrose conjecture was made by
Huisken and Ilmanen who 
proved the above theorem in the case
that the horizon of $(M^3,g)$ has only one component \cite{HI}.  
Specifically, they proved that $m \ge \sqrt{\frac{A_{max}}{16\pi}}$, 
where $A_{max}$ is the area of the largest component of the horizon of 
$(M^3,g)$. 
Their proof is as interesting as the result itself.
In the seventies, Geroch \cite{G} observed that in a
manifold with nonnegative scalar curvature, the 
Hawking mass of a sphere (but not surfaces with multiple components) 
was monotone increasing under a $1/H$ flow, where
$H$ is the mean curvature of the sphere.  Jang and Wald \cite{JW} 
proposed using this to attack the Riemannian Penrose conjecture 
by flowing the horizon
of the manifold out to infinity.  However, it is not
hard to concoct situations in which the $1/H$ flow of a sphere develops 
singularities, preventing the idea from working much of the time.  
Huisken and Ilmanen's approach then was  
to define a generalized $1/H$ flow
which sometimes ``jumps'' in order to prevent singularities from developing.

Other contributions have also been made by 
Herzlich \cite{He} using the Dirac
operator which Witten \cite{Wi} used to prove the positive mass theorem, 
by Gibbons \cite{Gi} in the special case of collapsing shells, 
by Tod \cite{T}, by Bartnik \cite{Ba3} for quasi-spherical metrics, and by
the author \cite{Bray_thesis} using isoperimetric surfaces.  There is also 
some interesting work of Ludvigsen and Vickers \cite{LV} using 
spinors and Bergqvist \cite{Berg}, both concerning the Penrose inequality for 
null slices of a space-time.

For more on physical discussions related to the Penrose inequality, 
gravitational collapse, and cosmic censorship, see also 
\cite{HP}, \cite{HE}, \cite{H2}, \cite{H1}, \cite{G}, and \cite{C}.

\section{Definitions and Setup}\label{sec:definitions}

Without loss of generality, we will be able to assume 
that an asymptotically flat metric (see definition 
\ref{def:asymptotically_flat}) 
has an even nicer form
at infinity in each end because of the following lemma and definition.

\begin{lemma}\label{lem:sy}
{\bf (Schoen, Yau \cite{SY5})}
Let $(M^3,g)$ be any asymptotically flat metric with nonnegative scalar curvature.  Then given any
$\epsilon > 0$, there exists
a metric $g_0$ with nonnegative scalar curvature which is harmonically flat at infinity (defined in the next definition) such that
\begin{equation}\label{eqn:epsilon}
1 - \epsilon \le \frac{g_0(\vec{v},\vec{v})}{g(\vec{v},\vec{v})} \le 1 + \epsilon
\end{equation}
for all nonzero vectors $\vec{v}$ in the tangent space at every point in $M$ and 
\begin{equation}
|\bar{m}_k - m_k| \le \epsilon
\end{equation}  
where $\bar{m}_k$ and $m_k$ are respectively the total masses of $(M^3,g_0)$ and $(M^3,g)$ in 
the $k$th end. 
\end{lemma}

Notice that because of equation \ref{eqn:epsilon}, the percentage difference 
in areas as well as 
lengths between the 
two metrics is arbitrarily small.  
Hence, since the mass changes arbitrarily little also and since  
inequality \ref{eqn:RPI} is a closed
condition, it follows that the Riemannian Penrose inequality for 
asymptotically flat manifolds 
follows from proving the inequality
for manifolds which are harmonically flat at infinity. 

\begin{definition}\label{def:harmonically_flat}
A Riemannian manifold is defined to be {\bf harmonically flat at infinity} if, 
outside a compact set, it is the disjoint union
of regions (which we will again call ends) 
with zero scalar curvature which are conformal to 
$(\real^3 \backslash B_1(0), \delta)$ with the 
conformal factor approaching a positive constant at infinity in each region.   
\end{definition}

Now it is fairly easy to define the total mass of an end of 
a manifold $(M^3,g_0)$ which is harmonically flat at infinity.
Define $g_{flat}$ to be a smooth metric on $M^3$ conformal to $g_0$ such that
in each end of $M^3$ in the above definition $(M^3,g_{flat})$ is isometric to 
$(\real^3 \backslash  B_1(0),\delta)$.  Define ${\cal U}_0(x)$ such that 
\begin{equation}
g_0 = {\cal U}_0(x)^4 g_{flat}.
\end{equation}
Then since $(M^3,g_0)$ has zero scalar curvature in each end, 
$(\real^3 \backslash B_1(0),{\cal U}_0(x)^4 \delta)$ must have zero scalar 
curvature.
This implies that ${\cal U}_0(x)$ is harmonic in \\
$(\real^3 \backslash B_1(0), \delta)$
(see equation \ref{eqn:scalar_curv} in appendix \ref{sec:harmonic}).  
Since ${\cal U}_0(x)$ is a harmonic function going to a constant at infinity,
we may expand it in terms of spherical harmonics to get 
\begin{equation}\label{eqn:mass}
{\cal U}_0(x) = a + \frac{b}{|x|} + {\cal O} \left(\frac{1}{|x|^2}\right), 
\end{equation}
where $a$ and $b$ are constants.
\begin{definition}\label{def:totalmass}
The {\bf total mass} (of an end) of a Riemannian 3-manifold 
which is harmonically flat at infinity
is defined to be $2ab$ in the above equation.
\end{definition}
While the constants $a$ and $b$ scale depending on how we  
represent $(M^3,g_0)$ as the disjoint union of a compact set and ends
in definition \ref{def:harmonically_flat}, it can be checked 
that $2ab$ does not.
Furthermore, this definition 
agrees with the standard definition of the 
total mass of an asymptotically flat manifold 
(defined in equation \ref{eqn:ADM_mass}) 
in the case that the manifold is harmonically
flat at infinity.  We choose to work with this definition because it is more 
convenient for the calculations we will be doing in this paper.

Now we turn our attention to the definition and properties of horizons.  
For convenience, we modify the topology of $M^3$ by compactifying all of the
ends of $M^3$ except for the chosen end 
by adding the points $\{\infty_k\}$.
\begin{definition}\label{def:calS}
Define ${\bf\cal S}$ to be the collection of surfaces which are smooth
compact boundaries of open sets in $M^3$ containing the points $\{\infty_k\}$. 
\end{definition}
All of the surfaces that we will be dealing with in this paper will be in ${\cal S}$. 
Also, we see that all of the 
surfaces in ${\cal S}$ divide $M^3$ into two regions, 
an inside (the open set)
and an outside (the complement of the open set).  
Thus, the notion of one surface in ${\cal S}$ (entirely) enclosing 
another surface
in ${\cal S}$ is well defined.
\begin{definition}
A {\bf horizon} of $(M^3,g)$ is 
any zero mean curvature surface in ${\cal S}$.  
\end{definition}

A horizon may have multiple 
components.  Furthermore, by minimizing area over surfaces in ${\cal S}$, a horizon is
guaranteed to exist when $M^3$ has more than one end.

\begin{definition}
A horizon is defined  
to be {\bf outermost} if it is not enclosed by another horizon.
\end{definition}

We note that when at least one horizon exists, there is always a unique
outermost horizon, with respect to the chosen end.

\begin{definition}\label{def:outerminimizing}
A surface $\Sigma \in {\cal S}$ is defined to be 
{\bf (strictly) outer-minimizing} if every other 
surface $\tilde{\Sigma} \in {\cal S}$ which encloses it
has (strictly) greater area.
\end{definition}

An outer-minimizing surface must have nonnegative mean curvature since 
otherwise the first variation formula would
imply that an outward variation would yield a surface with less area.  
Also, in the case that $\Sigma$ is a 
horizon, it is an outer-minimizing horizon if and only if it is not 
enclosed by a horizon with less area.  
Interestingly, every component of an outer-minimizing horizon must be a
2-sphere.  When the horizon is strictly outer-minimizing (or outermost),
this fact follows from a second variation argument (\cite{SY1}, or see section 
\ref{sec:stability}).  (Without the strictness assumption, tori become 
possibilities, but are then ruled out by \cite{GG}.) 
We will not use outermost
horizons in this paper, but simply point out that outermost horizons 
are always strictly outer-minimizing.
Hence, the following theorem, which is the main result of this paper,
is a slight generalization of the Riemannian Penrose conjecture.

\begin{theorem}\label{thm:Penrose}
Let $(M^3,g)$ be a complete, smooth, asymptotically flat  
$3$-manifold with nonnegative scalar curvature, total 
mass $m$, and an outer-minimizing horizon (with one or
more components) of total area $A$.
Then 
\begin{equation}\label{eqn:Penrose}
  m \ge \sqrt{\frac{A}{16\pi}} 
\end{equation}
with equality if and only if $(M^3,g)$ is isometric to a 
Schwarzschild manifold outside their respective outermost horizons.
\end{theorem}

Besides the flat metric on $\real^3$, the Schwarzschild manifolds are the only
other complete spherically symmetric 
$3$-manifolds with zero scalar curvature, and as previously mentioned can be
described explicitly as $(\real^3 \backslash  \{0\}, s)$ where
\begin{equation} s_{ij} = \left(1 + \frac{m}{2r}\right)
^{4} \delta_{ij},\end{equation} 
$r$ is the distance from the origin in $\real^3$, and
$m$ is a positive constant and equals the total mass of the 
manifold.  Then since the Schwarzschild manifolds have a single
minimal sphere which is the coordinate sphere of radius 
$m/2$, we can verify they give equality in the above theorem.

\section{Overview of the Proof}\label{sec:overview}

In this section we give the overview of the proof of theorem \ref{thm:Penrose}
which is a slight generalization of the Riemannian Penrose conjecture.  The
remainder of the paper is then devoted to proving and finding applications
for the claims made in this section.

As discussed in the previous section, without loss of generality for proving 
the Riemannian Penrose inequality for asymptotically flat manifolds we 
may restrict our attention to harmonically flat manifolds.  

\vspace{.1in}\noindent
{\bf Assumption:}
From this point on, we will assume that $(M^3,g_0)$ is a complete, smooth,
harmonically flat 3-manifold with nonnegative scalar curvature and an
outer-minimizing horizon $\Sigma_0$ (with one or more components) of total
area $A_0$, unless otherwise stated. 
\vspace{.1in}

\noindent
We will generalize our results to the asymptotically flat case 
and handle the case of equality of theorem \ref{thm:Penrose} 
in section \ref{sec:generalization}.

We define a continuous 
family of conformal metrics $\{g_t\}$ on $M^3$, where
\begin{equation}\label{eqn:ODE1}
   g_t = u_t(x)^4 g_0 
\end{equation}
and $u_0(x) \equiv 1$.  Given the metric $g_t$, define
\begin{equation}\label{eqn:ODE2}
   \Sigma(t) = \mbox{the outermost minimal area enclosure of }
               \Sigma_0 \mbox{ in } (M^3,g_t) 
\end{equation}
where $\Sigma_0$ is the original outer-minimizing horizon in $(M^3,g_0)$
and we stay inside the collection of surfaces ${\cal S}$ defined in the
previous section.
In the cases in which we are interested, 
$\Sigma(t)$ will not touch $\Sigma_0$, from which it follows
that $\Sigma(t)$ is actually a strictly outer-minimizing horizon of $(M^3,g_t)$.
Then given the horizon $\Sigma(t)$, define $v_t(x)$ such that
\begin{equation}\label{eqn:ODE3}
\left\{
\begin{array}{r l l l}
\Delta_{g_0} v_t(x) & \equiv & 0 & \mbox{ outside } \Sigma(t) \\
v_t(x) & = & 0 & \mbox{ on } \Sigma(t) \\
\lim_{x \rightarrow \infty} v_t(x) & = & -e^{-t} & \\
\end{array}
\right.
\end{equation}
and $v_t(x) \equiv 0$ inside $\Sigma(t)$.  Finally, given $v_t(x)$, define
\begin{equation}\label{eqn:ODE4}
u_t(x) = 1 + \int_0^t v_s(x) ds
\end{equation}
so that $u_t(x)$ is continuous in $t$ and has $u_0(x) \equiv 1$.  

\begin{theorem}\label{thm:existence}
Taken together, equations \ref{eqn:ODE1}, \ref{eqn:ODE2}, \ref{eqn:ODE3},
and \ref{eqn:ODE4} define a first order o.d.e.~in $t$ for $u_t(x)$ which has
a solution which is Lipschitz in the $t$ variable, 
$C^1$ in the $x$ variable everywhere, 
and smooth in the $x$ variable outside $\Sigma(t)$.
Furthermore, $\Sigma(t)$ is a smooth, 
strictly outer-minimizing horizon in $(M^3,g_t)$ 
for all $t \ge 0$, and 
$\Sigma(t_2)$ encloses but does not touch $\Sigma(t_1)$ 
for all $t_2 > t_1 \ge 0$.
\end{theorem}

Since $v_t(x)$ is a superharmonic function in $(M^3,g_0)$, it follows that
$u_t(x)$ is superharmonic as well, and from equation \ref{eqn:ODE4} we see
that $\lim_{x \rightarrow \infty} u_t(x) = e^{-t}$ and consequently that 
$u_t(x) > 0$ for all $t$.  Then since
\begin{equation}
R(g_t) = u_t(x)^{-5}(-8 \Delta_g + R(g))u_t(x)
\end{equation}
it follows that $(M^3,g_t)$ is an asymptotically flat manifold with 
nonnegative scalar curvature.

Even so, it still may not seem like $g_t$ is particularly naturally
defined since the rate of change of $g_t$ 
appears to depend on $t$ and the original metric $g_0$ 
in equation \ref{eqn:ODE3}.  We would prefer a flow where the rate
of change of $g_t$ is only a function of $g_t$ (and $M^3$ and $\Sigma_0$
perhaps), and interestingly enough this 
actually does turn out to be the case.  
In appendix \ref{sec:harmonic} we prove this very important fact and 
provide a very natural motivation for defining this conformal flow of 
metrics.

\begin{definition}
The function $A(t)$ is defined 
to be the total area of the horizon $\Sigma(t)$ in
$(M^3,g_t)$.
\end{definition}

\begin{definition}
The function $m(t)$ is defined to be the total mass of $(M^3,g_t)$ in 
the chosen end. 
\end{definition}

The next theorem is the key property of the conformal
flow of metrics.

\begin{theorem}\label{thm:monotone}
The function 
$A(t)$ is constant in $t$ and $m(t)$ is non-increasing in $t$, 
for all $t \ge 0$.
\end{theorem}

The fact that $A'(t) = 0$ follows from the fact that to first order the
metric is not changing on $\Sigma(t)$ (since $v_t(x) = 0$ there) and from 
the fact that to first order the area of $\Sigma(t)$ does not change as it
moves outward since $\Sigma(t)$ has zero mean curvature in $(M^3,g_t)$. 
We make this rigorous in section \ref{sec:A(t)}.
Hence, the interesting part of theorem \ref{thm:monotone} is proving that
$m'(t) \le 0$. Curiously, this follows from a nice trick 
using the Riemannian positive mass theorem. 

Another important aspect of this conformal flow of the metric is that 
outside the horizon $\Sigma(t)$, the manifold $(M^3,g_t)$ becomes more
and more spherically symmetric and ``approaches'' a Schwarzschild manifold
$(\real^3 \backslash \{0\}, s)$
in the limit as $t$ goes to 
$\infty$.  More precisely,

\begin{theorem}\label{thm:limit}
For sufficiently large $t$, there exists a 
diffeomorphism $\phi_t$ between $(M^3,g_t)$ outside the horizon $\Sigma(t)$ 
and a fixed Schwarzschild manifold $(\real^3 \backslash \{0\}, s)$ outside
its horizon.  Furthermore, for all $\epsilon > 0$, there exists a $T$
such that for all $t>T$, the metrics $g_t$ and $\phi^{*}_t(s)$ 
(when determining the lengths of unit vectors
of $(M^3,g_t)$) are within $\epsilon$ of each other 
and the total masses of the two manifolds are within
$\epsilon$ of each other.  Hence,
\begin{equation}
\lim_{t \rightarrow \infty} \frac{m(t)}{\sqrt{A(t)}} = \sqrt{\frac{1}{16\pi}}.
\end{equation}
\end{theorem}
   
Inequality \ref{eqn:Penrose} of theorem \ref{thm:Penrose} then 
follows from theorems 
\ref{thm:existence}, \ref{thm:monotone} and \ref{thm:limit},
for harmonically flat manifolds.  
In addition, in section \ref{sec:generalization} we will see that
the case of equality in theorem \ref{thm:Penrose} follows from the fact that
$m'(0) = 0$ if and only if 
$(M^3,g_0)$ is isometric to a Schwarzschild manifold
outside their respective outermost horizons. 
We will also generalize our results to the asymptotically flat case 
in that section.

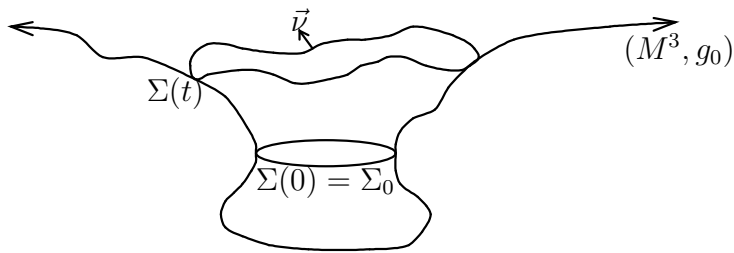
\vspace{.5in}
\begin{center}
\input{evolve_horizons_v2.eepic}
\end{center}
\vspace{.5in}

The diagrams above and below are meant to help illustrate some of the 
properties of the conformal flow of the metric.  The above picture is 
the original metric which has a strictly outer-minimizing horizon
$\Sigma_0$.  As $t$ increases, $\Sigma(t)$ moves outwards, but never 
inwards.  In the diagram below, we can observe one of the consequences 
of the fact that $A(t) = A_0$ is constant in $t$.  Since the metric
is not changing inside $\Sigma(t)$, all of the horizons $\Sigma(s)$,
$0 \le s \le t$ have area $A_0$ in $(M^3,g_t)$.  Hence, inside $\Sigma(t)$, 
the manifold $(M^3,g_t)$ becomes 
cylinder-like
in the sense that it is laminated by all of the previous horizons
which all have the same area $A_0$ with respect to the metric $g_t$.

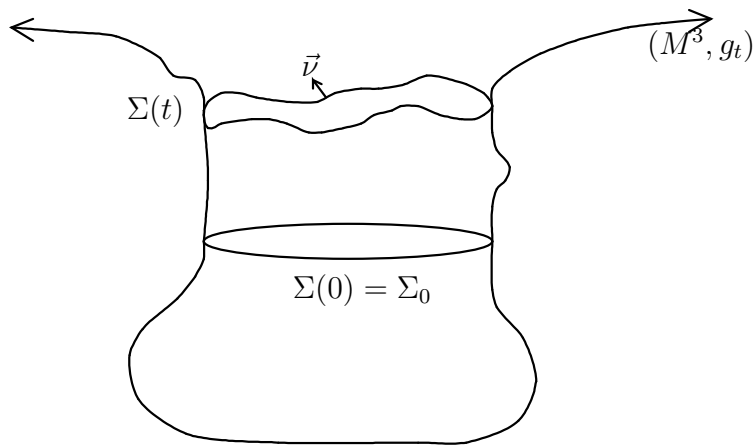
\vspace{.5in}
\begin{center}
\input{cylinder_v2.eepic}
\end{center}
\vspace{.5in}

Now let us suppose that the 
original horizon $\Sigma_0$ of $(M^3,g)$ had two components, for
example.
Then each of the components of the horizon will move outwards as $t$
increases, and at some point before they touch they 
will suddenly jump outwards
to form a horizon with a single component 
enclosing the previous horizons with two components.  Even horizons with
only one component will sometimes jump outwards, and it is interesting that
this phenomenon of surfaces jumping is also found in the Huisken-Ilmanen 
approach to the Penrose conjecture using their generalized $1/H$ flow.

\section{Existence and Regularity of the Flow of Metrics $\{g_t\}$}
\label{sec:existence}

In this section we will prove theorem \ref{thm:existence} which claims 
that there exists a solution $u_t(x)$ to the o.d.e.~in $t$ defined by equations
\ref{eqn:ODE1}, \ref{eqn:ODE2}, \ref{eqn:ODE3}, and \ref{eqn:ODE4} with 
certain regularity properties.  To do this, for each 
$\epsilon \in (0,\frac12)$ we will
define another family of conformal factors 
$u^\epsilon_t(x)$ which will be easy to 
prove exists, and then define
\begin{equation}
   u_t(x) = \lim_{\epsilon \rightarrow 0} u^\epsilon_t(x)
\end{equation}
which we will then show satisfies the original o.d.e.   

Define $\lb z \rb$ to be the greatest integer less than or equal to $z$.
Define
\begin{equation}
   \lb z \rbe = \epsilon \blb \frac{z}{\epsilon} \brb
\end{equation}
which we see is the greatest integer multiple of $\epsilon$ less than or 
equal to $z$.

Define
\begin{equation}\label{eqn:ODE1_ep}
   g^\epsilon_t = u^\epsilon_t(x)^4 g_0 
\end{equation}
and $u^\epsilon_0(x) \equiv 1$.  Given the metric $g^\epsilon_t$, 
define (for $t \ge 0$)
\begin{equation}\label{eqn:ODE2_ep}
\Sigma^\epsilon(t) = \left\{ \begin{array}{l l}
\Sigma_0 & \mbox{ if } t = 0 \\ & \\
\mbox{ the outermost minimal area enclosure } & \\
\mbox{ of } \Sigma^\epsilon(t - \epsilon) \mbox{ in }
(M^3, g^\epsilon_{t}) & \mbox { if } t = k \epsilon, \; 
k \in {\bf Z}^+ \\ & \\
\Sigma^\epsilon(\lb t \rbe) & \mbox{ otherwise, }
\end{array}
\right.
\end{equation}
where $\Sigma_0$ is the original outer-minimizing horizon in $(M^3,g_0)$
and we stay inside the collection of surfaces ${\cal S}$ defined in 
section \ref{sec:definitions}.
Given $\Sigma^\epsilon(t)$, define $v^\epsilon_t(x)$ such that
\begin{equation}\label{eqn:ODE3_ep}
\left\{
\begin{array}{r l l l}
\Delta_{g_0} v^\epsilon_t(x) & \equiv & 0 & \mbox{ outside } 
\Sigma^\epsilon(t) \\
v^\epsilon_t(x) & = & 0 & \mbox{ on } \Sigma^\epsilon(t) \\
\lim_{x \rightarrow \infty} v^\epsilon_t(x) & = & 
-(1-\epsilon)^{\blb \frac{t}{\epsilon} \brb} & \\
\end{array}
\right.
\end{equation}
and $v^\epsilon_t(x) \equiv 0$ inside $\Sigma^\epsilon(t)$.  
Finally, given $v^\epsilon_t(x)$, define
\begin{equation}\label{eqn:ODE4_ep}
u^\epsilon_t(x) = 1 + \int_0^t v^\epsilon_s(x) ds
\end{equation}
so that $u^\epsilon_t(x)$ is continuous in $t$ and 
has $u^\epsilon_0(x) \equiv 1$.

Notice that $\Sigma^\epsilon(t)$ and hence $v^\epsilon_t(x)$ 
are fixed for $t \in [k\epsilon,(k+1)\epsilon)$.  Furthermore,
for $t = k \epsilon$,  $k \in {\bf Z}^+$, 
$\Sigma^\epsilon(t)$ does not touch $\Sigma^\epsilon(t-\epsilon)$,
because it can be shown that $\Sigma^\epsilon(t-\epsilon)$
has negative mean curvature in $(M^3,g^\epsilon_t)$ and thus acts
as a barrier.  Hence, 
$\Sigma^\epsilon(t)$ is actually a strictly outer-minimizing horizon 
of $(M^3,g^\epsilon_t)$ and is smooth since $g^\epsilon_t$ is
smooth outside $\Sigma^\epsilon(t-\epsilon)$. 

From these considerations it follows that a solution 
$u^\epsilon_t(x)$ to  
equations \ref{eqn:ODE1_ep}, \ref{eqn:ODE2_ep}, \ref{eqn:ODE3_ep}, and 
\ref{eqn:ODE4_ep} always exists.  We can think of $u^\epsilon_t(x)$
as what results when we approximate the original o.d.e.~with a stepping procedure
where the step size equals $\epsilon$. 

Initially it might seem a little strange that we are requiring 
\begin{equation}\label{eqn:vvv}
\lim_{x \rightarrow \infty} v^\epsilon_t(x) =  
-(1-\epsilon)^{\blb \frac{t}{\epsilon} \brb}.
\end{equation}
However, this is done so that
\begin{equation}\label{eqn:v=-u}
\lim_{x \rightarrow \infty} v^\epsilon_t(x) =
- \lim_{x \rightarrow \infty} u^\epsilon_t(x)
\end{equation}
for $t$ values which are an integral multiple of $\epsilon$.  This is 
necessary to prove that the rate of change of $g^\epsilon_t$ is just a function
of $g^\epsilon_t$ and not of $g_0$ (when $t$ is an integral multiple of 
$\epsilon$), and the argument is the same as the one given for the original
o.d.e.~in appendix \ref{sec:harmonic}.

In corollary \ref{cor:regularity} of 
appendix \ref{sec:regularity}, we show that we have 
upper bounds on the $C^{k,\alpha}$ ``norms'' 
of the horizons $\{\Sigma^\epsilon(t)\}$, as defined in 
definition \ref{def:k,alpha,S}.
Furthermore, these upper bounds do not depend on $\epsilon$, and depend only
on $T$, $\Sigma_0$, $g_0$, and the choice of coordinate charts for $M^3$.
Hence, not only are these surfaces smooth, but any 
limits of these surfaces, which
we will be dealing with later, will also be smooth.

\begin{lemma}\label{lem:enclose_ep}
The horizon $\Sigma^\epsilon(t_2)$ encloses $\Sigma^\epsilon(t_1)$
for all $t_2 \ge t_1 \ge 0$. 
\end{lemma}
{\it Proof.}
Follows from the definition of $\Sigma^\epsilon(t)$ in equation
\ref{eqn:ODE2_ep}.  \qed

\begin{lemma}\label{lem:omae_ep}
The horizon $\Sigma^{\epsilon}(t)$ is the outermost
minimal area enclosure of $\Sigma_0$ in $(M^3,g_t^{\epsilon})$ when 
$t = k \epsilon$, $k \in {\bf Z}^+$.
\end{lemma}
{\it Proof.}
The proof is by induction on $k$.  The case when $k=1$ follows by definition
from equation \ref{eqn:ODE2_ep}.  Now assume that the lemma is true 
for $t = (k-1) \epsilon$.  Then since the metric is not changing inside
$\Sigma^{\epsilon}((k-1)\epsilon)$ for $(k-1)\epsilon \le t \le k \epsilon$,
the outermost minimal area enclosure of $\Sigma_0$ in $(M^3,g_{k\epsilon})$
must be outside $\Sigma^{\epsilon}((k-1)\epsilon)$.  Then the lemma follows
for $t = k \epsilon$
from equation \ref{eqn:ODE2_ep} again.  \qed

\begin{lemma}\label{lem:Lipschitz}
The functions $u_t^\epsilon(x)$ are positive, 
bounded, locally Lipschitz functions (in $x$ and $t$),
with uniform
Lipschitz constants independent of $\epsilon$.
\end{lemma}
{\it Proof.}  Since $u_t^\epsilon(x)$ is superharmonic and goes to a positive
constant at infinity (which is seen by integrating $v_t^\epsilon(x)$ 
in equation \ref{eqn:ODE4_ep}), $u_t^\epsilon(x) > 0$.  And since
$v_t^\epsilon(x) \le 0$, then from equation \ref{eqn:ODE4_ep} it follows that
$u_t^\epsilon(x) \le 1$.  The fact that $u_t^\epsilon(x)$ is Lipschitz in $t$ 
follows from equation \ref{eqn:ODE4_ep} and the fact that 
$v_t^\epsilon(x)$ is bounded, and the fact that $u_t^\epsilon(x)$ is Lipschitz
in $x$ follows from the fact that $v_t^\epsilon(x)$ is Lipschitz (with
Lipschitz constant depending on $t$) by 
corollary \ref{cor:regularity} of appendix \ref{sec:regularity}.  \qed

\begin{corollary}\label{cor:limit}
There exists a subsequence $\{\epsilon_i\}$ converging to zero such that 
\begin{equation}\label{eqn:uexists}
   u_t(x) = \lim_{\epsilon_i \rightarrow 0} u_t^{\epsilon_i}(x)
\end{equation}
exists, is locally Lipschitz (in $x$ and $t$), and the convergence is 
locally uniform.  
Hence, we may define the metric 
\begin{equation}
   g_t  = \lim_{\epsilon_i \rightarrow 0} g_t^{\epsilon_i} = u_t(x)^4 g_0
\end{equation}
for $t \ge 0$ as well.
\end{corollary}
{\it Proof.}  Follows from lemma \ref{lem:Lipschitz}.  We will stay inside
this subsequence for the remainder of this section.  \qed

\begin{definition}
Define $\{\tilde{\Sigma}_\gamma(t)\}$ 
to be the collections of limit surfaces of 
$\Sigma^{\epsilon_i}(t)$ in the limit as $\epsilon_i$ approaches $0$.
\end{definition} 
Initially we define the limits in the measure theoretic sense, which must 
exist since it is possible to bound the areas of the surfaces
$\Sigma^{\epsilon_i}(t)$ from above and below and to show that they are all
contained in a compact region.
However, we also have bounds on the $C^{k,\alpha}$ ``norms'' 
(see definition \ref{def:k,alpha,S}) of the horizons
$\Sigma^{\epsilon}(t)$
by corollary \ref{cor:regularity}
in appendix \ref{sec:regularity}.  Hence, since the bounds given in corollary 
\ref{cor:regularity} are independent
of $\epsilon$, the
above limits are also 
true in the Hausdorff distance sense, and the limit surfaces are all smooth.

\begin{theorem}\label{thm:enclose} 
The limit surface $\tilde{\Sigma}_{\gamma_2}(t_2)$ encloses 
$\tilde{\Sigma}_{\gamma_1}(t_1)$ for all $t_2 > t_1 \ge 0$ and for 
any $\gamma_1$ and $\gamma_2$.
\end{theorem}
{\it Proof.}
This theorem would be trivial and would follow directly from lemma
\ref{lem:enclose_ep} if not for the fact that there can be multiple limit
surfaces for $\Sigma^{\epsilon_i}(t)$ as $\epsilon_i$ goes to zero.  Instead,
we have some work to do.

Given any $\delta > 0$, choose $\bar{\epsilon}$ such that
\begin{equation}\label{eqn:e1}
   |u_t^{\epsilon_i}(x) - u_t(x)| < \delta  \;\;\mbox{ for all }
                                    \epsilon_i < \bar{\epsilon} ,
\end{equation}
which we can do since the convergence in equation \ref{eqn:uexists} is 
locally Lipschitz and since the $u^\epsilon_t(x)$ are all harmonic outside
a compact set.  Next,
choose $\epsilon_1, \epsilon_2 < \bar{\epsilon}$ such that
\begin{eqnarray}
\mbox{dis}(\Sigma^{\epsilon_1}(t_1), \tilde{\Sigma}_{\gamma_1}(t_1)) < \delta \\ 
\mbox{dis}(\Sigma^{\epsilon_2}(t_2), \tilde{\Sigma}_{\gamma_2}(t_2)) < \delta
\end{eqnarray} 
which is possible since we have convergence in the Hausdorff distance sense.  
Note
that we define 
\begin{equation}
   \mbox{dis}(S,T) = \max(\sup_{x \in S} \inf_{y \in T} d(x,y), \;
   \sup_{x \in T} \inf_{y \in S} d(x,y))
\end{equation}
where $d(x,y)$ is the usual distance function in $(M^3, g_0)$.
Then by equation \ref{eqn:e1} and the triangle inequality,
\begin{equation}\label{eqn:e2}
   |u_t^{\epsilon_1}(x) - u_t^{\epsilon_2}(x)| < 2 \delta 
\end{equation}
so that by equation \ref{eqn:ODE4_ep} we have that
\begin{equation}
   \left| \int_0^t (v_s^{\epsilon_1}(x) - v_s^{\epsilon_2}(x)) \; ds \right| 
   < 2 \delta 
\end{equation}
which by the triangle inequality once again implies that
\begin{equation}\label{eqn:star}
   \left| \int_{t_1}^{t_2} 
   (v_s^{\epsilon_1}(x) - v_s^{\epsilon_2}(x)) \; ds \right| < 4 \delta .
\end{equation}

Since $\Sigma^{\epsilon_2}(t_2)$ encloses 
$\Sigma^{\epsilon_2}(s)$ for $s < t_2$, it follows from the maximum 
principle that 
\begin{equation}\label{eqn:e3}
   v_s^{\epsilon_2}(x) \le v_{t_2}^{\epsilon_2}(x) \; (1-\epsilon_2)^{
   (\lb \frac{s}{\epsilon_2} \rb - \lb \frac{t_2}{\epsilon_2} \rb)}.
\end{equation}
Similarly, since $\Sigma^{\epsilon_1}(t_1)$ is enclosed by 
$\Sigma^{\epsilon_1}(s)$ for $s > t_1$, it also follows from the maximum 
principle that 
\begin{equation}\label{eqn:e4}
   v_s^{\epsilon_1}(x) \ge v_{t_1}^{\epsilon_1}(x) \; (1-\epsilon_1)^{
   (\lb \frac{s}{\epsilon_1} \rb - \lb \frac{t_1}{\epsilon_1} \rb)}.
\end{equation}

Now suppose that 
$\tilde{\Sigma}_{\gamma_2}(t_2)$ did not enclose 
$\tilde{\Sigma}_{\gamma_1}(t_1)$ for some $t_2 > t_1 \ge 0$.  Then choose
$x_0$ strictly inside $\tilde{\Sigma}_{\gamma_1}(t_1)$ and strictly outside 
$\tilde{\Sigma}_{\gamma_2}(t_2)$, and choose any positive $\delta < \frac12
\max(\mbox{dis}(x_0,\tilde{\Sigma}_{\gamma_1}(t_1) ) , \;
\mbox{dis}(x_0,\tilde{\Sigma}_{\gamma_2}(t_2) ) )$.  Then we must have
\begin{equation}\label{eqn:e5}
v_{t_2}^{\epsilon_2}(x_0) < 0  \;\;\mbox{ and }\;\; 
v_{t_1}^{\epsilon_1}(x_0) = 0
\end{equation}
Then combining equations \ref{eqn:star}, \ref{eqn:e3}, \ref{eqn:e4}
and \ref{eqn:e5} yields
\begin{equation}\label{eqn:e6}
    - v_{t_2}^{\epsilon_2}(x_0) \int_{t_1}^{t_2} (1-\epsilon_2)^{
   (\lb \frac{s}{\epsilon_2} \rb - \lb \frac{t_2}{\epsilon_2} \rb)} 
   \; ds < 4 \delta.
\end{equation}
Let $\lim_{\epsilon_2 \rightarrow 0} v_{t_2}^{\epsilon_2}(x_0) = - \alpha$,
which must be negative by the maximum principle since
$\Sigma^{\epsilon_2}(t_2)$ is approaching $\tilde{\Sigma}_{\gamma_2}(t_2)$
smoothly.  Then taking the limit of inequality \ref{eqn:e6} 
as $\delta$ and $\epsilon_2$ both go to zero yields
\begin{equation}
   \alpha \int_{t_1}^{t_2} e^{(t_2 - s)} \; ds \le 0,
\end{equation}
a contradiction since $t_2$ is strictly greater than $t_1$.  Hence,
$\tilde{\Sigma}_{\gamma_2}(t_2)$ must enclose 
$\tilde{\Sigma}_{\gamma_1}(t_1)$ for all $t_2 > t_1 \ge 0$, proving the
theorem.  \qed

\begin{theorem}\label{thm:Ap=0}
Let $A_0$ be the area of the original outer-minimizing horizon $\Sigma_0$
with respect to the original metric $g_0$.
Then
\begin{equation}
   |\tilde{\Sigma}_{\gamma}(t)|_{g_t} = A_0
\end{equation}
for all $t \ge 0$, where $|\cdot|_{g_t}$ denotes area with respect to the 
metric $g_t$.
\end{theorem}
{\it Proof.}  Given in the next section.  \qed

\begin{corollary}\label{cor:Ap=0}
In addition,
\begin{equation}
   |\tilde{\Sigma}_{\gamma}(t_1)|_{g_{t_2}} = A_0
\end{equation}
for all $t_2 \ge t_1 \ge 0$.
\end{corollary}
{\it Proof.}  This statement follows from the previous two theorems and the
fact that $v_t^{\epsilon_i}(x)$ is defined to be zero inside 
$\Sigma^{\epsilon_i}(t)$
so that $g_{t_1} = g_{t_2}$ inside $\Sigma(t_2)$ for $t_2 \ge t_1 \ge 0$.  \qed

\begin{definition}\label{def:sigma}
Define $\Sigma(t)$ to be the outermost minimal area enclosure of the 
original horizon $\Sigma_0$ in $(M^3,g_t)$.
\end{definition}
We note that the outermost minimal area enclosure of a smooth region is 
well-defined in that it always exists and is unique \cite{BT}.

\begin{lemma}\label{lem:SigmaA_0}
It is also true that 
\begin{equation}
   |\Sigma(t)|_{g_t} = A_0.
\end{equation}
for all $t \ge 0$.
\end{lemma}
{\it Proof.}
Since by lemma \ref{lem:omae_ep} $\Sigma^{\epsilon_i}(t)$ is the outermost
minimal area enclosure of $\Sigma_0$ in $(M^3,g_t^{\epsilon_i})$, the 
result follows from theorem \ref{thm:Ap=0} and the fact that 
$u^{\epsilon_i}_t(x)$ is locally uniformly Lipschitz and 
is converging to $u_t(x)$ locally uniformly.
This is also discussed in the next section.  \qed

\begin{lemma}\label{lem:enclose1}
The surface $\Sigma(t_2)$ encloses $\tilde{\Sigma}_\gamma(t_1)$ for all
$\gamma$ and $t_2 > t_1 \ge 0$.
\end{lemma}
{\it Proof.} 
By lemma \ref{lem:SigmaA_0}, the minimal area enclosures of $\Sigma_0$
in $(M^3,g_{t_2})$ have area $A_0$.  But by corollary \ref{cor:Ap=0},
$\tilde{\Sigma}_\gamma(t_1)$ has area $A_0$ in $(M^3,g_{t_2})$.  Hence,
since $\Sigma(t_2)$ is the outermost minimal area enclosure by definition, 
it follows that $\Sigma(t_2)$ encloses $\tilde{\Sigma}_\gamma(t_1)$.  \qed

\begin{lemma}\label{lem:enclose2}
The surface
$\tilde{\Sigma}_\gamma(t_2)$ encloses $\Sigma(t_1)$ for all $\gamma$ and 
$t_2 > t_1 \ge 0$.
\end{lemma}
{\it Proof.}  
Suppose $\tilde{\Sigma}_\gamma(t_2)$ did not (entirely) enclose $\Sigma(t_1)$
for some $t_2 > t_1 \ge 0$.  Then we choose a subsequence 
$\{\epsilon_i'\} \subset \{\epsilon_i\}$ converging to zero such that 
$\Sigma^{\epsilon_i'}(t_2)$ is converging to $\tilde{\Sigma}_\gamma(t_2)$
in the Hausdorff distance sense.

On the other hand, by lemma \ref{lem:SigmaA_0} we have 
\begin{equation}
   \lim_{\epsilon_i' \rightarrow 0} |\Sigma(t_1)|_{g^{\epsilon_i'}_{t_1}} =
   |\Sigma(t_1)|_{g_{t_1}} = A_0.
\end{equation}
However, for a given $\epsilon_i'$, the metric is shrinking outside of 
$\Sigma^{\epsilon_i'}(t_2)$ for $t_1 \le t \le t_2$ by a uniform amount 
which can be made independent of $\epsilon_i'$, and follows from 
lemma \ref{lem:enclose_ep} and 
equations \ref{eqn:ODE3_ep} and \ref{eqn:ODE4_ep}.  Hence,
\begin{equation}
   \lim_{\epsilon_i' \rightarrow 0} |\Sigma(t_1)|_{g^{\epsilon_i'}_{t_2}} =
   |\Sigma(t_1)|_{g_{t_2}} < A_0,
\end{equation}
which by definition \ref{def:sigma} violates lemma \ref{lem:SigmaA_0}. \qed

\begin{corollary}\label{cor:enclose}
The surface
$\Sigma(t_2)$ encloses $\Sigma(t_1)$ for all 
$t_2 > t_1 \ge 0$.
\end{corollary}
{\it Proof.} Follows directly from the two previous lemmas.  \qed

\begin{definition}\label{def:limits}
Define
\begin{equation}
\begin{array}{ll}
   \Sigma^+(t) = \lim_{s \rightarrow t^+} \Sigma(s), \hspace{.5in} &
   \tilde{\Sigma}^+(t) = \lim_{s \rightarrow t^+} \tilde{\Sigma}_\gamma(s), \\
& \\   
   \Sigma^-(t) = \lim_{s \rightarrow t^-} \Sigma(s), \hspace{.5in} &
   \tilde{\Sigma}^-(t) = \lim_{s \rightarrow t^-} \tilde{\Sigma}_\gamma(s),
\end{array}
\end{equation}
where we define $\Sigma^-(0) = \Sigma_0 = \tilde{\Sigma}^-(0)$.
\end{definition}
We note that by the inclusion properties of theorem \ref{thm:enclose}
and corollary \ref{cor:enclose} that these limits are always unique.
\begin{definition}
Define the jump times $J$ to be the set of all $t \ge 0$ with  
$\Sigma^+(t) \ne \Sigma^-(t)$. 
\end{definition}
\begin{theorem} \label{thm:equal}
We have the inclusion property that 
for all $t_2 > t_1 \ge 0$, the surfaces 
$\Sigma(t_2), \tilde{\Sigma}_{\gamma}(t_2),
\Sigma^+(t_2), \tilde{\Sigma}^+(t_2), 
\Sigma^-(t_2), \tilde{\Sigma}^-(t_2) $
all enclose the surfaces \\ 
$\Sigma(t_1), \tilde{\Sigma}_{\gamma}(t_1),
\Sigma^+(t_1), \tilde{\Sigma}^+(t_1), 
\Sigma^-(t_1), \tilde{\Sigma}^-(t_1)$.  
Also, for $t \ge 0$ 
\begin{equation}\label{eqn:equal}
\Sigma(t) = \Sigma^+(t) = \tilde{\Sigma}^+(t) \;\;\mbox{ and encloses }\;\;
\Sigma^-(t) = \tilde{\Sigma}^-(t),
\end{equation}
and all five surfaces are smooth  
with area $A_0$ in $(M^3,g_t)$. 
Furthermore,
except for $t \in J$, all five surfaces in equation \ref{eqn:equal} are equal,
$\tilde{\Sigma}_\gamma(t)$ is single valued, and 
$\Sigma(t) = \tilde{\Sigma}_\gamma(t)$.  In addition, 
the set $J$ is countable.
\end{theorem}
{\it Proof.}
The inclusion property follows directly from lemmas \ref{lem:enclose1}
and \ref{lem:enclose2}.  These two lemmas also prove that 
$\Sigma^+(t)$ encloses $\tilde{\Sigma}^+(t)$ and that 
$\tilde{\Sigma}^+(t)$ encloses $\Sigma^+(t)$, proving that they are equal.
Similarly it follows that $\Sigma^-(t) = \tilde{\Sigma}^-(t)$, and these
two lemmas also imply that $\Sigma^+(t) = \tilde{\Sigma}^+(t)$ encloses
$\Sigma^-(t) = \tilde{\Sigma}^-(t)$.  

By corollary \ref{cor:limit} and theorem \ref{thm:Ap=0}, 
all five surfaces have area $A_0$ in $(M^3,g_t)$.
By corollary \ref{cor:enclose}, 
$\Sigma(t)$ is enclosed by $\Sigma^+(t)$.  On the other hand,  
$\Sigma^+(t)$ has area $A_0$ and $\Sigma(t)$ encloses all other minimal area
enclosures of $\Sigma_0$, so $\Sigma(t)$ encloses $\Sigma^+(t)$.  Hence,  
$\Sigma(t) = \Sigma^+(t)$.

By the definition of $J$, all five surfaces are equal for  $t \notin J$.
Since each $\tilde{\Sigma}_\gamma(t)$ is sandwiched between the surfaces
$\tilde{\Sigma}^-(t)$ and $\tilde{\Sigma}^+(t)$ which are equal, 
it follows that $\tilde{\Sigma}_\gamma(t)$ is a single limit surface and 
equals $\tilde{\Sigma}^-(t) = \tilde{\Sigma}^+(t)$ for these $t \notin J$,
from which it follows that $\Sigma(t) = \tilde{\Sigma}_\gamma(t)$.

Define $\Delta V(t)$ to equal the volume enclosed by $\Sigma^+(t)$ but 
not by $\Sigma^-(t)$.  Then $\sum_{t \in J \cap [0,T]} \Delta V(t)$ 
is finite for all $T > 0$ since it is less than or equal to the volume
enclosed by $\Sigma(T+1)$ but not by $\Sigma_0$ 
which is finite.  Also, $\Delta V(t) > 0$
for $t \in J$ since by corollary \ref{cor:regularity} 
of appendix \ref{sec:regularity}  
these surfaces are all uniformly smooth.
Hence, $J$ is countable.  \qed

\begin{definition} \label{def:v} 
Given the horizon $\Sigma(t)$, define $v_t(x)$ such that
\begin{equation}\label{eqn:ODE3_copy}
\left\{
\begin{array}{r l l l}
\Delta_{g_0} v_t(x) & \equiv & 0 & \mbox{ outside } \Sigma(t) \\
v_t(x) & = & 0 & \mbox{ on } \Sigma(t) \\
\lim_{x \rightarrow \infty} v_t(x) & = & -e^{-t} & \\
\end{array}
\right.
\end{equation}
and $v_t(x) \equiv 0$ inside $\Sigma(t)$.  
\end{definition}

\begin{lemma}\label{lem:integral}
Using the definition of $v_t(x)$ given above and the definition of 
$u_t(x)$ given in corollary \ref{cor:limit}, we have that
\begin{equation}\label{eqn:ODE4_copy}
u_t(x) = 1 + \int_0^t v_s(x) ds .
\end{equation}
\end{lemma}
{\it Proof.}
By theorem \ref{thm:equal}, $\lim_{\epsilon_i \rightarrow 0}
v^{\epsilon_i}_s(x) = v_s(x)$ for almost every value of $s \ge 0$.
Then using 
the inclusion property of theorem \ref{thm:equal} it follows that we can
pass the limit into the integral in equation \ref{eqn:ODE4_ep}, proving the
lemma.  \qed

\begin{lemma}\label{lem:C1}
The surfaces $\{\Sigma(t)\}$ do not touch for different values of $t \ge 0$.  
Furthermore, the metric $g_t$ and the conformal factor $u_t(x)$ are 
$C^1$ in $x$, and $\Sigma(t_1)$, $\Sigma^+(t_1)$, and $\Sigma^-(t_1)$
are outer-minimizing horizons with area $A_0$ in $(M^3,g_{t_2})$
for $t_2 \ge t_1 \ge 0$.  
\end{lemma}
{\it Proof.}
First we note that by equation
\ref{eqn:ODE3_ep} and corollary \ref{cor:regularity} in 
appendix \ref{sec:regularity}
that there exists a 
locally uniform
constant $c$ such that $v^\epsilon_t(x) - c f(x)$ is convex, where $f(x)$ 
is any locally defined 
smooth function with all of its second derivatives greater than or equal
to one.  Hence, by equation \ref{eqn:ODE4_ep}, the same statement is true 
for $u^\epsilon_t(x)$, and then also for 
$u_t(x)$ after taking the limit.  Hence, for all $x \in M$ the directional
derivatives of $u_t(x)$ exist in all directions and 
\begin{equation}\label{eqn:uconvex}
   \nabla_{\vec{w}} u_t(x) \le -\nabla_{-\vec{w}} u_t(x).
\end{equation}

The mean curvature of a surface $\Sigma$ in $(M^3,g_t)$ is 
\begin{equation}\label{eqn:meancurvature2}
   H = u_t(x)^{-2} H_0 + 4u_t(x)^{-3}\nabla_{\vec{\nu}}u_t(x)
\end{equation}
where $H_0$ is the mean curvature and $\vec{\nu}$ is the
outward pointing unit normal vector of $\Sigma$ in $(M^3,g_0)$.
Hence, by equation \ref{eqn:uconvex}, 
the mean curvature of $\Sigma$ on the outside is always less than or equal
to the mean curvature of $\Sigma$ on the inside (but still using an outward
pointing unit normal vector). 
 Since $\Sigma_0$ has zero mean curvature in $(M^3,g_0)$ and
since $u_t(x) = 1$ on $\Sigma_0$ and $u_t(x) \le 1$ everywhere
else, it follows that the mean curvature of $\Sigma_0$ (on the outside) 
is nonpositive
in $(M^3,g_t)$ and thus acts as a barrier for $\Sigma(t)$.

Suppose $\Sigma(t)$ touched $\Sigma_0$ at $x_0$ for some $t>0$.  
Since $\Sigma(t)$ is outside $\Sigma_0$, 
it must have nonpositive mean curvature in $(M^3,g_0)$.  Then by equation
\ref{eqn:ODE4_copy} we have $\nabla_{\vec{\nu}}u_t(x_0) < 0$, where 
$\vec{\nu}$ is the outward pointing unit normal vector to $\Sigma(t)$ and
$\Sigma_0$ at $x_0$ in $(M^3,g_0)$.  
Hence, by equation \ref{eqn:meancurvature2}
it follows that $\Sigma(t)$ has negative mean curvature on the outside in
$(M^3,g_t)$.  This is a contradiction, since by the first variation formula 
and the pseudo-convexity of $u_t(x)$ we could then flow $\Sigma(t)$ out and
decrease its area, contradicting the fact that it is defined to be a minimal
area enclosure of $\Sigma_0$ in $(M^3,g_t)$.  Hence, $\Sigma(t)$ does 
not touch $\Sigma_0$.

Then since the mean curvature of $\Sigma(t)$ on the outside is less than or 
equal to the mean curvature on the inside, it follows that both mean 
curvatures for $\Sigma(t)$ in $(M^3,g_t)$ must be zero.  Otherwise, it would
be possible to do a variation of $\Sigma(t)$ which decreased its area in 
$(M^3,g_t)$.  Similarly,
by corollary \ref{cor:Ap=0} and theorem \ref{thm:equal}, 
$|\Sigma(t_1)|_{g_{t_2}} = A_0$
too, so that by the same first variation argument 
the mean curvatures of $\Sigma(t_1)$ on the outside and inside 
in $(M^3,g_2)$ must both also be zero, for $t_2 \ge t_1 \ge 0$.  Hence, 
$\Sigma(t_1)$ is a horizon in $(M^3,g_{t_2})$.  
 
Hence, $\Sigma(t_1)$ acts as a barrier for $\Sigma(t_2)$ for $t_2 > t_1 \ge 0$.
As discussed in appendix 
\ref{sec:harmonic}, the o.d.e.~for $u_t(x)$ is actually translation 
invariant in $t$ since the rate of change of $g_t$ is a function of $g_t$
and not of $t$.  Thus, the argument proving that $\Sigma(t_2)$ does not touch 
$\Sigma(t_1)$ is essentially the same as the argument two paragraphs above
which proved that $\Sigma(t)$ did not touch $\Sigma_0$.

To see that $u_t(x)$ is $C^1$ in $x$, 
we first observe that 
by equation \ref{eqn:ODE4_copy}
the directional derivatives in $x$ for $u_t(x)$ equal 
\begin{equation}
   \nabla_{\vec{w}} u_t(x) = \int_0^t \nabla_{\vec{w}} v_s(x) ds.
\end{equation} 
Furthermore, $v_t(x)$ is smooth everywhere except on $\Sigma(t)$
and has locally uniform bounds on all of its derivatives off of 
$\Sigma(t)$ because of the 
regularity of $\Sigma(t)$ coming from theorem
\ref{thm:equal} and corollary \ref{cor:regularity} 
of appendix \ref{sec:regularity}.  
Hence, $u_t(x)$ will be $C^1$ at $x = x_0$ if and only
if the set $I(x_0) =\{t \;\;|\;\; x_0 \in \Sigma(t) \}$ 
has measure zero in $\real$.  But since we just got through proving that
the horizons $\Sigma(t)$ do not touch for different values of $t$, $I(x_0)$
is at most one point, so $u_t(x)$ is $C^1$ in $x$.  \qed

Theorem \ref{thm:existence} then follows from  
corollary \ref{cor:limit},
definition \ref{def:sigma},
corollary \ref{cor:enclose},
theorem \ref{thm:equal},
definition \ref{def:v}, 
lemma \ref{lem:integral},  
lemma \ref{lem:C1}, and 
corollary \ref{cor:regularity} in appendix \ref{sec:regularity}.

\section{Proof That $A(t)$ Is a Constant}
\label{sec:A(t)}

The fact that $A(t)$, defined to be the area of the horizon $\Sigma(t)$
in $(M^3,g_t)$, is constant in $t$ 
was proven already in lemma \ref{lem:SigmaA_0} of the previous section.  
However, this lemma relied entirely on theorem \ref{thm:Ap=0}, the proof
of which we have postponed until now.

{\it Proof of theorem \ref{thm:Ap=0}:  }
We continue with the same notation as in the previous section.
\begin{definition}
We define $A^\epsilon(t) = |\Sigma^\epsilon(t)|_{g^\epsilon_t}$ and 
\begin{equation}\label{eqn:Adiff}
   \Delta A^\epsilon(t) = A^\epsilon(t + \epsilon) 
   - A^\epsilon(t),
\end{equation}
where $t$ is a nonnegative integer multiple of $\epsilon$.
\end{definition}
Then if we can prove that for all $T > 0$,
\begin{equation}\label{eqn:mainarea}
   \lim_{\epsilon \rightarrow 0} A^\epsilon(t) = A_0
\end{equation}
for all $t \in [0,T]$,
theorem \ref{thm:Ap=0} follows since each limit surface 
$\tilde{\Sigma}_\gamma(t)$ is the limit of a $\Sigma^{\epsilon_i}(t)$ for
some choice of $\{\epsilon_i\}$ converging to zero, all of the surfaces 
involved are uniformly smooth (corollary \ref{cor:regularity}), 
and $\{g^{\epsilon_i}_t\}$ are all uniformly Lipschitz
and are converging uniformly to $g_t$ (corollary \ref{cor:limit}). 

Our first observation is that since $v^\epsilon_t(x) \le 0$, the metric 
$g^\epsilon_t$ gets smaller pointwise as $t$ increase.  Hence,
$\Delta A^\epsilon(t)$ is always negative, where we are requiring that 
$t = k\epsilon$ for some nonnegative integer $k$.  Then 
since $A^\epsilon(0) = A_0$ by definition, it follows that 
\begin{equation}\label{eqn:area_ub}
   A^\epsilon(t) \le A_0.
\end{equation}
Hence, all that
we need to prove equation \ref{eqn:mainarea} is a lower bound on
$\Delta A^\epsilon(t) / \epsilon$ 
which goes to zero as $\epsilon$ goes to zero 
for ``most'' values of $t$.

For convenience, we first set $t = 0$ and 
estimate $\Delta A^\epsilon(0)$.  This estimate will then be generalizable for
all values of $t$ since the flow is independent of the base metric 
$g_0$ and $t$ as discussed in section \ref{sec:existence} and 
described in appendix $\ref{sec:harmonic}$.  
Since $u^\epsilon_\epsilon(x) = 1 + \epsilon v^\epsilon_0(x)$, it follows that
\begin{equation}
   A^\epsilon(\epsilon) = \int_{\Sigma^\epsilon(\epsilon)} 
   \left(1 + \epsilon v^\epsilon_0(x)\right)^4 \; dA_{g^\epsilon_0}
\end{equation}
where as usual $dA_{g^\epsilon_0}$ is the area form of 
$\Sigma^\epsilon(\epsilon)$
with respect to $g^\epsilon_0$.  
Then since $\Sigma^\epsilon(0)$ is outer-minimizing in
$(M^3,g^\epsilon_0)$, we also have
\begin{equation}
   A^\epsilon(0) \le \int_{\Sigma^\epsilon(\epsilon)} 1 \; dA_{g^\epsilon_0}.
\end{equation} 
Hence, 
\begin{eqnarray}
   \Delta A^\epsilon(0) & \ge & \int_{\Sigma^\epsilon(\epsilon)} 
   [\left(1 + \epsilon v^\epsilon_0(x)\right)^4 - 1] \; 
   dA_{g^\epsilon_0}
   \\ & \ge & 4\epsilon 
   \int_{\Sigma^\epsilon(\epsilon)} \left( v^\epsilon_0(x)  
   - \epsilon^2 \right) \; dA_{g^\epsilon_0}
\end{eqnarray}
which follows from expanding and the fact that $-1 \le v^\epsilon_0(x) \le 0$,
so that
\begin{equation}
   \Delta A^\epsilon(0) \ge 4\epsilon 
   |\Sigma^\epsilon(\epsilon)|_{g^\epsilon_0}
   \left(\min_{\Sigma^\epsilon(\epsilon)} v^\epsilon_0(x) - \epsilon^2 \right).
\end{equation}

More generally, from the discussion at the end of appendix \ref{sec:harmonic},
we see that if we define 
\begin{equation}
   \bar{v}^\epsilon_t(x) = v^\epsilon_t(x) / u^\epsilon_t(x),
\end{equation}
then 
\begin{equation}
   \Delta A^\epsilon(t) \ge 4\epsilon 
   |\Sigma^\epsilon(t + \epsilon)|_{g^\epsilon_t}
   \left(\min_{\Sigma^\epsilon(t + \epsilon)} \bar{v}^\epsilon_t(x) 
   - \epsilon^2 \right),
\end{equation}
where as before $t$ is a nonnegative integral multiple of $\epsilon$.  
Furthermore, 
\begin{equation}
|\Sigma^\epsilon(t + \epsilon)|_{g^\epsilon_t} \;\;\le\;\; 
|\Sigma^\epsilon(t + \epsilon)|_{g^\epsilon_{t+\epsilon}}(1-\epsilon)^{-4} 
\;\;\le\;\; A_0 (1-\epsilon)^{-4},
\end{equation}
since $g^\epsilon_t$ and $g^\epsilon_{t+\epsilon}$ are conformally related
within a factor of $(1-\epsilon)^{-4}$ and by inequality \ref{eqn:area_ub}.  
Thus,
\begin{equation}\label{eqn:deltaAlb}
   \Delta A^\epsilon(t) \ge 4\epsilon A_0 (1-\epsilon)^{-4} 
   \left(\min_{\Sigma^\epsilon(t + \epsilon)} \bar{v}^\epsilon_t(x) 
   - \epsilon^2 \right).
\end{equation}

\begin{definition}
We define $V^\epsilon(t)$ to be the volume of the region 
enclosed by $\Sigma^\epsilon(t)$ which is outside the original horizon
$\Sigma_0$ and define 
\begin{equation}
\Delta V^\epsilon(t) = V^\epsilon(t+\epsilon) - 
V^\epsilon(t),
\end{equation}
where $t$ is a nonnegative integer multiple of $\epsilon$.  
\end{definition}

We now require $t \in [0,T]$, for any fixed $T > 0$.  
Since $(M^3,g_0)$ is asymptotically flat and $u^\epsilon_t(x)$ has uniform
upper and lower bounds, there exists a uniform upper bound 
$V_0$ for $V^\epsilon(T)$
which is independent of $\epsilon$.  Hence, 
by lemma \ref{lem:enclose_ep} it follows that
\begin{equation}
   \Delta V^\epsilon(k\epsilon) \le V_0 \sqrt{\epsilon}
\end{equation}
for $k = 0$ to $\lb \frac{T}{\epsilon} \rb - 1$ except at most 
$\lb \frac{1}{\sqrt{\epsilon}} \rb$ values of $k$.

Furthermore, it is possible to define a 
function $f$ which is independent of $\epsilon$ and 
only depends on $T$, $\Sigma_0$, $g_0$, and 
the uniform regularity bounds on the surfaces 
$\Sigma^\epsilon(t)$ in corollary \ref{cor:regularity} in 
appendix \ref{sec:regularity} such that
\begin{equation}
   \max_{x \in \Sigma^\epsilon((k+1)\epsilon)} \mbox{dis}(x,
   \Sigma^\epsilon(k\epsilon)) \le f(\Delta V^\epsilon(k\epsilon)),
\end{equation}
where $f$ is a continuous, increasing function which equals zero at zero.  
Also, from corollary \ref{cor:regularity} we have
\begin{equation}
   |\nabla v^\epsilon_t(x)|_{g_0} \le K
\end{equation}
for $0 \le t \le T$, where $K$ is also independent of $\epsilon$.  Hence,
since ${v}^\epsilon_{k\epsilon}(x) = 0$ on $\Sigma^\epsilon(k\epsilon)$,
\begin{equation}
   \min_{\Sigma^\epsilon((k+1)\epsilon)} {v}^\epsilon_{k\epsilon}(x) \ge - K
   f(V_0 \sqrt{\epsilon}),
\end{equation}  
and since $u^\epsilon_t(x) \ge (1-\epsilon)^{\lb \frac{T}{\epsilon} \rb}$ 
for $0 \le t \le T$,
\begin{equation}
   \min_{\Sigma^\epsilon((k+1)\epsilon)} {\bar{v}}^\epsilon_{k\epsilon}(x) 
   \ge - K f(V_0 \sqrt{\epsilon})
   (1-\epsilon)^{-\lb \frac{T}{\epsilon} \rb},
\end{equation}  
for $k = 0$ to $\lb \frac{T}{\epsilon} \rb - 1$ except at most 
$\lb \frac{1}{\sqrt{\epsilon}} \rb$ values of $k$.  At these exceptional
values of $k$ we will just 
use the fact that ${\bar{v}}^\epsilon_{k\epsilon}(x) > -1$.

Hence, from equation \ref{eqn:deltaAlb} 
\begin{equation}\label{eqn:deltaAlb2}
   \Delta A^\epsilon(k\epsilon) \ge - 4\epsilon A_0 (1-\epsilon)^{-4} 
   \left( K f(V_0 \sqrt{\epsilon}) (1-\epsilon)^{-\lb \frac{T}{\epsilon} \rb}
   + \epsilon^2 \right)
\end{equation}
for $k = 0$ to $\lb \frac{T}{\epsilon} \rb - 1$ except at most 
$\lb \frac{1}{\sqrt{\epsilon}} \rb$ values of $k$ and 
\begin{equation}\label{eqn:deltaAlb3}
   \Delta A^\epsilon(k\epsilon) \ge - 4\epsilon A_0 (1-\epsilon)^{-4} 
   \left(1 + \epsilon^2 \right)
\end{equation}
for all values of $k$ including the exceptional ones.

Then from equation \ref{eqn:Adiff} we have that
\begin{eqnarray}
   A^\epsilon(n\epsilon) & = & A_0 + \sum_{k=0}^{n-1} 
   \Delta A^\epsilon(k\epsilon) \\ 
   & \ge & \begin{array}{l}
   A_0 - 4 T A_0 (1-\epsilon)^{-4} 
   \left( K f(V_0 \sqrt{\epsilon})(1-\epsilon)^{-\lb \frac{T}{\epsilon} \rb} 
   + \epsilon^2 \right) \\
   - 4\sqrt{\epsilon} A_0 (1-\epsilon)^{-4} 
   \left(1 + \epsilon^2 \right) 
\end{array}\label{eqn:area_lb}
\end{eqnarray}
for $1 \le n \le \lb \frac{T}{\epsilon} \rb$.  Then since $f$ goes to
zero at zero, equation \ref{eqn:mainarea} follows 
from inequalities \ref{eqn:area_ub} and \ref{eqn:area_lb} 
for $0 \le t \le T$.  Since $T > 0$ was arbitrary, this proves 
theorem \ref{thm:Ap=0}.  \qed

\section{Green's Functions at Infinity and the Riemannian
         Positive Mass Theorem}\label{sec:green}

In this section we will prove a few theorems about certain Green's functions
on asymptotically flat 3-manifolds with nonnegative scalar curvature which 
will be needed in the next section to prove that $m(t)$ is nonincreasing
in $t$.  However, the theorems in this section, which follow from and 
generalize the 
Riemannian positive mass theorem, are also of independent 
interest.

The results of this section are closely related to the beautiful 
ideas used by
Bunting and Masood-ul-Alam in \cite{BM} to prove the non-existence of 
multiple black holes in asymptotically flat, static, vacuum space-times.  
Hence, while theorems \ref{thm:green1} and \ref{thm:green2} 
do not appear in their
paper, these two theorems follow from a natural extension of their techniques.

\begin{definition}\label{def:energy1}
Given a complete, asymptotically flat manifold $(M^3,\bar{g})$ with multiple 
asymptotically flat ends (with one chosen end), define 
\begin{equation}
   {\cal E}(\bar{g}) = \inf_{\phi} \left\{ \frac{1}{2\pi}
          \int_{(M^3,\bar{g})} | \nabla \phi |^2 \; dV \right\}
\end{equation}
where the infimum is taken over all smooth 
$\phi(x)$ which go to one in the chosen 
end and zero in the other ends. 
\end{definition}
Without loss of generality, we may assume that $(M^3,\bar{g})$ is actually
harmonically flat at infinity (as defined in section \ref{sec:definitions}).
Then since such a modification can be done so as to change the metric 
uniformly pointwise as small as one likes (by lemma \ref{lem:sy}), 
it follows that ${\cal E}(\bar{g})$ changes as small as one likes as well.  We remind the
reader that the total mass of $(M^3,\bar{g})$ also changes arbitrarily little 
with such a deformation. 

From standard theory it follows that the infimum in the above definition
is achieved by the Green's function $\phi(x)$ which satisfies
\begin{equation}\label{eqn:phi}
\left\{\begin{array}{rll}
\lim_{x \rightarrow \infty_0} \phi(x) &=& 1 \\
\Delta \phi &=& 0 \\
\lim_{x \rightarrow \infty_k} \phi(x) &=& 0 \;\mbox{ for all } k \neq 0
\end{array}\right.
\end{equation}
where $\infty_k$ are the points at infinity of the various asymptotically
flat ends and $\infty_0$ is infinity in the chosen end.
Define the level sets of $\phi(x)$ to be
\begin{equation}
\Sigma_l = \{ x \;|\; \phi(x) = l  \}
\end{equation}
for $0 < l < 1$.  Then it follows from Sard's theorem and the smoothness 
of $\phi(x)$ that $\Sigma(l)$ is a smooth surface for almost every $l$. 
Then by the co-area formula it follows that
\begin{equation}
   {\cal E}(\bar{g}) \;=\; \frac{1}{2\pi} \int_0^1 dl \; \int_{\Sigma_l} |\nabla \phi| \;dA 
        \;=\; \frac{1}{2\pi} \int_0^1 dl \; \int_{\Sigma_l} 
            \frac{d\phi}{d\eta} \;dA 
\end{equation}
since $\nabla \phi$ is orthogonal to the unit normal vector $\eta$ of 
$\Sigma_l$.  But by the divergence theorem (and since $\phi(x)$ is harmonic),
$\int_\Sigma \frac{d\phi}{d\eta} \;dA $ is constant for all homologous
$\Sigma$.  Hence,
\begin{equation}\label{eqn:energy1}
   {\cal E}(\bar{g}) = \frac{1}{2\pi} \int_\Sigma \frac{d\phi}{d\eta} \;dA
\end{equation} 
where $\Sigma$ is any surface in $M^3$ in ${\cal S}$ (which is the set of 
smooth, compact boundaries of open regions
which contains the points at infinity $\{\infty_k\}$ in all of the ends except
the chosen one).  Then since $(M^3,\bar{g})$ is harmonically flat at infinity, 
we know that
in the chosen end $\phi(x) = 1 - c/|x| + {\cal O}(1/|x|^2)$, so that 
from equation \ref{eqn:energy1} it follows that 
\begin{equation}\label{eqn:phi_exp}
   \phi(x) = 1 - \frac{{\cal E}(\bar{g})}{2|x|} + {\cal O}\left(\frac{1}{|x|^2}\right)
\end{equation}
by letting $\Sigma$ in equation \ref{eqn:energy1} be a large sphere in the 
chosen end.

\begin{theorem}\label{thm:green1}
Let $(M^3,\bar{g})$ be a complete, smooth, 
asymptotically flat 3-manifold with nonnegative 
scalar curvature which has multiple asymptotically flat ends and total
mass $\bar{m}$ in the chosen end.  Then
\begin{equation}\label{eqn:m>e}
   \bar{m} \ge {\cal E}(\bar{g})
\end{equation}
with equality if and only if $(M^3,\bar{g})$ has zero scalar curvature and is 
conformal to $(\real^3,\delta)$ minus a finite number of points.
\end{theorem} 
{\it Proof.}  Again, without loss of generality we will assume that 
$(M^3,\bar{g})$ is harmonically flat at infinity, which by definition means that
all of the ends are harmonically flat at infinity.  In the picture below,
$(M^3,\bar{g})$ has three ends, the chosen end (at the top of the picture) and
two other ends. 
\vspace{.5in}
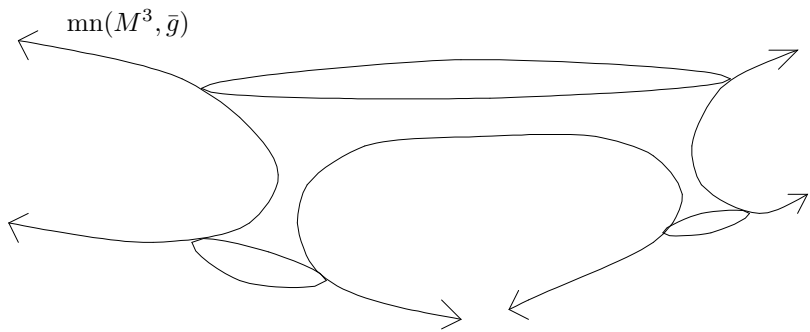
\begin{center}
\setlength{\unitlength}{0.00066667in}
\begingroup\makeatletter\ifx\SetFigFont\undefined%
\gdef\SetFigFont#1#2#3#4#5{%
  \reset@font\fontsize{#1}{#2pt}%
  \fontfamily{#3}\fontseries{#4}\fontshape{#5}%
  \selectfont}%
\fi\endgroup%
{\renewcommand{\dashlinestretch}{30}
\begin{picture}(6249,2496)(0,-10)
\path(87,2262)(237,2337)
\path(87,2262)(162,2112)
\path(5937,2187)(6162,2187)(6012,2037)
\path(12,837)(162,912)
\path(12,837)(87,687)
\path(6087,1062)(6237,1062)(6162,912)
\path(3987,312)(3912,162)
\path(3912,162)(4062,87)
\path(3387,237)(3537,87)
\path(3537,87)(3387,12)
\path(87,2262)	(155.431,2249.747)
	(221.434,2237.817)
	(285.068,2226.196)
	(346.390,2214.871)
	(405.460,2203.828)
	(462.335,2193.053)
	(517.073,2182.532)
	(569.734,2172.251)
	(620.376,2162.198)
	(669.056,2152.357)
	(760.767,2133.261)
	(845.334,2114.853)
	(923.224,2097.023)
	(994.904,2079.660)
	(1060.840,2062.656)
	(1121.500,2045.900)
	(1177.351,2029.283)
	(1228.860,2012.694)
	(1276.493,1996.024)
	(1362.000,1962.000)

\path(1362,1962)	(1411.287,1938.225)
	(1468.993,1906.460)
	(1531.903,1868.739)
	(1596.803,1827.097)
	(1660.477,1783.569)
	(1719.711,1740.187)
	(1812.000,1662.000)

\path(1812,1662)	(1896.787,1577.015)
	(1945.942,1524.168)
	(1994.692,1466.273)
	(2039.357,1404.718)
	(2076.251,1340.894)
	(2101.693,1276.191)
	(2112.000,1212.000)

\path(2112,1212)	(2107.297,1159.430)
	(2092.034,1105.401)
	(2068.144,1051.456)
	(2037.562,999.139)
	(1964.060,905.557)
	(1887.000,837.000)

\path(1887,837)	(1812.819,796.028)
	(1722.619,763.297)
	(1673.140,749.727)
	(1621.612,737.863)
	(1568.687,727.586)
	(1515.015,718.777)
	(1461.249,711.320)
	(1408.041,705.095)
	(1356.042,699.984)
	(1305.904,695.869)
	(1213.819,690.153)
	(1137.000,687.000)

\path(1137,687)	(1057.586,686.702)
	(967.275,690.967)
	(917.004,694.948)
	(862.772,700.235)
	(804.167,706.882)
	(740.779,714.945)
	(672.194,724.478)
	(598.001,735.537)
	(517.788,748.176)
	(431.143,762.451)
	(337.653,778.416)
	(288.213,787.049)
	(236.908,796.125)
	(183.685,805.652)
	(128.494,815.635)
	(71.283,826.082)
	(12.000,837.000)

\path(3537,87)	(3475.926,97.368)
	(3417.040,107.651)
	(3360.289,117.867)
	(3305.622,128.032)
	(3252.987,138.165)
	(3202.334,148.282)
	(3153.610,158.400)
	(3106.765,168.537)
	(3018.502,188.934)
	(2937.133,209.612)
	(2862.247,230.707)
	(2793.431,252.356)
	(2730.274,274.698)
	(2672.363,297.869)
	(2619.286,322.006)
	(2570.632,347.248)
	(2484.943,401.592)
	(2412.000,462.000)

\path(2412,462)	(2354.969,538.771)
	(2330.361,587.569)
	(2308.879,640.020)
	(2290.912,693.656)
	(2276.850,746.009)
	(2262.000,837.000)

\path(2262,837)	(2264.837,907.587)
	(2279.812,991.342)
	(2304.632,1072.927)
	(2337.000,1137.000)

\path(2337,1137)	(2424.505,1227.121)
	(2482.615,1271.851)
	(2545.972,1314.521)
	(2611.333,1353.751)
	(2675.453,1388.160)
	(2735.090,1416.370)
	(2787.000,1437.000)

\path(2787,1437)	(2877.266,1460.813)
	(2930.063,1470.272)
	(2986.884,1478.285)
	(3046.919,1484.990)
	(3109.355,1490.525)
	(3173.378,1495.028)
	(3238.178,1498.639)
	(3302.940,1501.495)
	(3366.853,1503.736)
	(3429.105,1505.499)
	(3488.882,1506.923)
	(3545.373,1508.147)
	(3597.764,1509.308)
	(3687.000,1512.000)

\path(3687,1512)	(3776.050,1515.841)
	(3828.402,1518.073)
	(3884.880,1520.364)
	(3944.669,1522.604)
	(4006.951,1524.681)
	(4070.907,1526.483)
	(4135.721,1527.900)
	(4200.575,1528.820)
	(4264.651,1529.132)
	(4327.133,1528.724)
	(4387.201,1527.486)
	(4444.040,1525.305)
	(4496.831,1522.072)
	(4587.000,1512.000)

\path(4587,1512)	(4633.925,1503.306)
	(4690.033,1490.893)
	(4752.066,1475.183)
	(4816.763,1456.601)
	(4880.864,1435.569)
	(4941.111,1412.509)
	(4994.243,1387.845)
	(5037.000,1362.000)

\path(5037,1362)	(5102.152,1306.865)
	(5171.797,1231.552)
	(5230.295,1146.464)
	(5262.000,1062.000)

\path(5262,1062)	(5259.191,1005.020)
	(5239.980,942.334)
	(5213.029,883.230)
	(5187.000,837.000)

\path(5187,837)	(5117.280,756.345)
	(5037.000,687.000)

\path(5037,687)	(4944.378,630.332)
	(4888.613,600.649)
	(4827.841,570.352)
	(4763.035,539.655)
	(4695.167,508.774)
	(4625.210,477.922)
	(4554.135,447.315)
	(4482.916,417.167)
	(4412.524,387.692)
	(4343.932,359.105)
	(4278.112,331.621)
	(4216.036,305.454)
	(4158.677,280.818)
	(4107.008,257.929)
	(4062.000,237.000)

\path(4062,237)	(4009.815,211.279)
	(3968.369,190.463)
	(3912.000,162.000)

\path(6162,2187)	(6097.116,2162.564)
	(6036.908,2139.399)
	(5981.158,2117.394)
	(5929.644,2096.441)
	(5838.450,2057.249)
	(5761.568,2020.943)
	(5697.239,1986.643)
	(5643.707,1953.473)
	(5562.000,1887.000)

\path(5562,1887)	(5490.437,1796.684)
	(5454.570,1739.134)
	(5420.734,1676.974)
	(5390.493,1612.942)
	(5365.413,1549.777)
	(5347.060,1490.217)
	(5337.000,1437.000)

\path(5337,1437)	(5337.421,1365.586)
	(5351.003,1281.960)
	(5376.333,1200.854)
	(5412.000,1137.000)

\path(5412,1137)	(5499.736,1056.607)
	(5556.269,1019.364)
	(5617.793,985.732)
	(5681.717,956.947)
	(5745.453,934.241)
	(5806.410,918.848)
	(5862.000,912.000)

\path(5862,912)	(5928.045,917.665)
	(6006.701,942.210)
	(6053.505,962.661)
	(6106.756,989.150)
	(6167.555,1022.116)
	(6237.000,1062.000)

\path(1437,687)	(1526.396,716.210)
	(1606.323,710.744)
	(1699.027,696.206)
	(1747.662,686.546)
	(1796.479,675.831)
	(1844.474,664.464)
	(1890.645,652.849)
	(1973.496,630.495)
	(2037.000,612.000)

\path(2037,612)	(2089.522,595.005)
	(2157.250,571.049)
	(2233.113,542.022)
	(2310.041,509.816)
	(2380.964,476.321)
	(2438.812,443.426)
	(2487.000,387.000)

\path(2487,387)	(2396.445,343.848)
	(2308.279,335.273)
	(2257.558,333.848)
	(2203.969,334.185)
	(2148.677,336.172)
	(2092.850,339.699)
	(2037.655,344.655)
	(1984.259,350.930)
	(1933.828,358.412)
	(1887.530,366.992)
	(1812.000,387.000)

\path(1812,387)	(1761.067,408.990)
	(1699.998,441.933)
	(1634.442,482.567)
	(1570.046,527.632)
	(1512.457,573.869)
	(1467.321,618.016)
	(1437.000,687.000)

\path(5112,762)	(5143.633,800.212)
	(5228.659,845.460)
	(5279.718,867.202)
	(5330.355,886.478)
	(5412.000,912.000)

\path(5412,912)	(5453.079,919.892)
	(5506.403,927.303)
	(5566.677,933.345)
	(5628.608,937.132)
	(5686.900,937.778)
	(5736.259,934.396)
	(5787.000,912.000)

\path(5787,912)	(5755.367,873.788)
	(5670.341,828.540)
	(5619.282,806.798)
	(5568.645,787.522)
	(5487.000,762.000)

\path(5487,762)	(5445.921,754.108)
	(5392.597,746.697)
	(5332.323,740.655)
	(5270.392,736.868)
	(5212.100,736.222)
	(5162.741,739.604)
	(5112.000,762.000)

\path(1512,1887)	(1588.557,1857.737)
	(1676.813,1846.008)
	(1731.547,1840.817)
	(1792.634,1836.055)
	(1859.533,1831.708)
	(1931.706,1827.761)
	(2008.613,1824.199)
	(2089.714,1821.006)
	(2174.470,1818.169)
	(2262.343,1815.672)
	(2352.791,1813.501)
	(2445.278,1811.640)
	(2492.116,1810.821)
	(2539.261,1810.075)
	(2586.646,1809.399)
	(2634.204,1808.791)
	(2681.866,1808.249)
	(2729.565,1807.773)
	(2777.235,1807.359)
	(2824.806,1807.006)
	(2872.213,1806.712)
	(2919.388,1806.476)
	(3012.771,1806.167)
	(3104.415,1806.065)
	(3193.782,1806.155)
	(3280.332,1806.422)
	(3363.526,1806.851)
	(3442.823,1807.427)
	(3517.686,1808.136)
	(3587.575,1808.963)
	(3651.949,1809.892)
	(3710.271,1810.909)
	(3762.000,1812.000)

\path(3762,1812)	(3853.931,1814.474)
	(3907.693,1816.152)
	(3966.051,1818.113)
	(4028.556,1820.353)
	(4094.757,1822.864)
	(4164.202,1825.639)
	(4236.443,1828.673)
	(4311.026,1831.959)
	(4387.504,1835.491)
	(4465.424,1839.261)
	(4544.336,1843.263)
	(4623.789,1847.490)
	(4703.333,1851.937)
	(4782.518,1856.597)
	(4860.892,1861.463)
	(4938.006,1866.528)
	(5013.408,1871.786)
	(5086.648,1877.230)
	(5157.275,1882.854)
	(5224.839,1888.652)
	(5288.890,1894.617)
	(5348.975,1900.742)
	(5404.646,1907.020)
	(5455.451,1913.446)
	(5500.940,1920.013)
	(5574.167,1933.542)
	(5637.000,1962.000)

\path(5637,1962)	(5566.524,1993.935)
	(5484.567,2009.122)
	(5433.668,2016.477)
	(5376.826,2023.658)
	(5314.545,2030.652)
	(5247.325,2037.447)
	(5175.671,2044.030)
	(5100.086,2050.391)
	(5021.071,2056.516)
	(4939.130,2062.393)
	(4854.766,2068.011)
	(4768.481,2073.356)
	(4680.779,2078.418)
	(4592.161,2083.183)
	(4503.131,2087.641)
	(4414.192,2091.777)
	(4325.845,2095.581)
	(4238.596,2099.041)
	(4152.945,2102.143)
	(4069.395,2104.876)
	(3988.451,2107.229)
	(3910.613,2109.187)
	(3836.386,2110.741)
	(3766.272,2111.876)
	(3700.773,2112.582)
	(3640.393,2112.846)
	(3585.635,2112.656)
	(3537.000,2112.000)

\path(3537,2112)	(3489.945,2110.771)
	(3436.985,2108.837)
	(3378.605,2106.231)
	(3315.292,2102.985)
	(3247.532,2099.134)
	(3175.809,2094.710)
	(3100.611,2089.746)
	(3022.422,2084.276)
	(2941.729,2078.332)
	(2859.018,2071.947)
	(2774.774,2065.155)
	(2689.483,2057.989)
	(2603.631,2050.481)
	(2517.704,2042.665)
	(2432.188,2034.574)
	(2347.568,2026.241)
	(2264.330,2017.699)
	(2182.960,2008.982)
	(2103.944,2000.121)
	(2027.769,1991.151)
	(1954.918,1982.105)
	(1885.879,1973.014)
	(1821.138,1963.914)
	(1761.179,1954.836)
	(1706.489,1945.814)
	(1657.554,1936.881)
	(1578.892,1919.413)
	(1512.000,1887.000)

\put(462,2337){\makebox(0,0)[lb]{\smash{{{\SetFigFont{10}{12.0}{\rmdefault}{\mddefault}{\updefault}$(M^3,\bar{g})$}}}}}
\end{picture}
}
\end{center}
Then we consider the metric
$(M^3, {{\tilde{g}}} )$, with $ {{\tilde{g}}}  = \phi(x)^4 \bar{g}$, drawn below.  
Since $\phi(x)$ goes to zero 
(and is bounded above by $C/|x|$)
in all of the harmonically flat ends other than the chosen one, 
the metric $ {{\tilde{g}}}  = \phi(x)^4 \bar{g}$ in each end is conformal to a punctured ball
with the conformal factor being a bounded harmonic function 
to the fourth power in the
punctured ball.  Hence, by the removable singularity theorem, this harmonic
function can be extended to the whole ball, which proves that the metric
$ {{\tilde{g}}} $ can be extended smoothly over all of the points at infinity
in the compactified ends.  
\vspace{.5in}
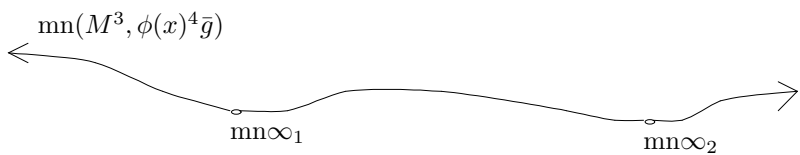
\begin{center}
\setlength{\unitlength}{0.00066667in}
\begingroup\makeatletter\ifx\SetFigFont\undefined%
\gdef\SetFigFont#1#2#3#4#5{%
  \reset@font\fontsize{#1}{#2pt}%
  \fontfamily{#3}\fontseries{#4}\fontshape{#5}%
  \selectfont}%
\fi\endgroup%
{\renewcommand{\dashlinestretch}{30}
\begin{picture}(6174,1086)(0,-10)
\put(5010,239){\ellipse{74}{36}}
\put(1785,314){\ellipse{74}{36}}
\path(162,852)(12,777)(162,702)
\path(6012,552)(6162,477)
\path(6162,477)(6012,327)
\path(12,777)	(82.770,772.223)
	(148.576,767.586)
	(209.664,763.062)
	(266.282,758.623)
	(318.677,754.241)
	(367.097,749.890)
	(452.996,741.169)
	(525.959,732.239)
	(587.962,722.880)
	(640.983,712.874)
	(687.000,702.000)

\path(687,702)	(741.571,684.355)
	(806.170,658.643)
	(877.362,627.398)
	(951.713,593.160)
	(1025.786,558.464)
	(1096.148,525.848)
	(1159.365,497.847)
	(1212.000,477.000)

\path(1212,477)	(1288.340,451.720)
	(1336.392,437.124)
	(1393.084,420.637)
	(1459.952,401.821)
	(1538.535,380.236)
	(1630.372,355.442)
	(1681.741,341.704)
	(1737.000,327.000)

\path(1812,327)	(1887.462,323.528)
	(1952.560,321.049)
	(2008.391,319.561)
	(2056.054,319.065)
	(2131.271,321.049)
	(2187.000,327.000)

\path(2187,327)	(2233.404,337.838)
	(2288.247,355.458)
	(2348.717,377.599)
	(2412.000,402.000)
	(2475.283,426.401)
	(2535.753,448.542)
	(2590.596,466.162)
	(2637.000,477.000)

\path(2637,477)	(2712.121,485.640)
	(2803.430,490.750)
	(2853.457,492.147)
	(2905.494,492.861)
	(2958.862,492.957)
	(3012.881,492.502)
	(3066.873,491.564)
	(3120.158,490.206)
	(3172.057,488.498)
	(3221.892,486.503)
	(3312.650,481.924)
	(3387.000,477.000)

\path(3387,477)	(3432.703,472.903)
	(3484.599,467.208)
	(3541.795,460.107)
	(3603.401,451.791)
	(3668.525,442.452)
	(3736.274,432.281)
	(3805.757,421.471)
	(3876.082,410.213)
	(3946.359,398.697)
	(4015.693,387.117)
	(4083.195,375.664)
	(4147.972,364.528)
	(4209.133,353.902)
	(4265.786,343.978)
	(4317.039,334.947)
	(4362.000,327.000)

\path(4362,327)	(4444.472,310.098)
	(4495.246,298.887)
	(4548.795,287.062)
	(4602.432,275.542)
	(4653.471,265.246)
	(4737.000,252.000)

\path(4737,252)	(4815.251,249.128)
	(4877.427,249.846)
	(4962.000,252.000)

\path(5037,252)	(5120.943,246.712)
	(5182.909,244.950)
	(5262.000,252.000)

\path(5262,252)	(5331.634,280.101)
	(5413.312,325.418)
	(5494.335,371.526)
	(5562.000,402.000)

\path(5562,402)	(5650.018,423.235)
	(5705.128,432.682)
	(5769.983,441.641)
	(5846.339,450.333)
	(5935.956,458.977)
	(5986.286,463.350)
	(6040.590,467.793)
	(6099.088,472.334)
	(6162.000,477.000)

\put(237,927){\makebox(0,0)[lb]{\smash{{{\SetFigFont{10}{12.0}{\rmdefault}{\mddefault}{\updefault}$(M^3,\phi(x)^4\bar{g})$}}}}}
\put(1737,102){\makebox(0,0)[lb]{\smash{{{\SetFigFont{10}{12.0}{\rmdefault}{\mddefault}{\updefault}$\infty_1$}}}}}
\put(4962,27){\makebox(0,0)[lb]{\smash{{{\SetFigFont{10}{12.0}{\rmdefault}{\mddefault}{\updefault}$\infty_2$}}}}}
\end{picture}
}
\end{center}
Furthermore, $(M^3, {{\tilde{g}}} )$ has nonnegative scalar curvature since 
$(M^3,\bar{g})$ has nonnegative scalar curvature and $\phi(x)$ is harmonic with 
respect to $\bar{g}$ (see equation \ref{eqn:scalar_curv} in appendix 
\ref{sec:harmonic}).  Moreover, 
$(M^3 \cup \{\infty_k\}, {{\tilde{g}}} )$ has nonnegative scalar curvature too
since in a neighborhood of each $\infty_k$ the manifold is conformal 
to a ball, with the conformal factor being a positive harmonic function to the
fourth power.  Hence, since $(M^3 \cup \{\infty_k\}, {{\tilde{g}}} )$ is a 
complete 3-manifold with nonnegative scalar curvature with a single 
harmonically flat end, we may apply the Riemannian positive mass theorem
to this manifold to conclude that the total mass of this manifold,
which we will call ${{\tilde{m}}}$, is nonnegative.

Now we will compute ${{\tilde{m}}}$ in terms of $\bar{m}$ and ${\cal E}(\bar{g})$.  Since 
$\bar{g}$ is harmonically flat at infinity, we know that by definition 
\ref{def:totalmass} we have $\bar{g} = \bar{{\cal U}}(x)^4 \bar{g}_{flat}$, 
where $(M^3,\bar{g}_{flat})$ is isometric to 
$(\real^3 \backslash B_r(0), \delta)$ 
in the harmonically flat end of $(M^3, \bar{g})$, 
\begin{equation}\label{eqn:u_exp}
   \bar{{\cal U}}(x) = 1 + \frac{\bar{m}}{2|x|} +{\cal O}\left(\frac{1}{|x|^2}\right),
\end{equation}
and the scale of the harmonically flat coordinate chart 
$\real^3 \backslash B_r(0)$ has been chosen so that $\bar{{\cal U}}(x)$ goes
to one at infinity in the chosen end.  Furthermore, since 
$ {{\tilde{g}}}  = \phi(x)^4 \bar{g}$ is also harmonically flat at infinity (which follows
from equation \ref{eqn:identity} in appendix \ref{sec:harmonic}), we have
$ {{\tilde{g}}}  = {\tilde{{\cal U}}}(x)^4 \bar{g}_{flat} = \phi(x)^4\bar{{\cal U}}(x)^4 \bar{g}_{flat}$
where
\begin{equation}\label{eqn:ubar_exp}
   {\tilde{{\cal U}}}(x) = 1 + \frac{{{\tilde{m}}}}{2|x|} 
                         +{\cal O}\left(\frac{1}{|x|^2}\right).
\end{equation}  
Then comparing equations \ref{eqn:phi_exp},
\ref{eqn:u_exp}, and \ref{eqn:ubar_exp} with 
${\tilde{{\cal U}}}(x) = \bar{{\cal U}}(x)\phi(x)$ yields
\begin{equation}\label{eqn:mtilde}
  {{\tilde{m}}} = \bar{m} - {\cal E}(\bar{g}) \ge 0
\end{equation}
by the Riemannian positive mass theorem, which proves 
inequality \ref{eqn:m>e} for harmonically flat manifolds.
Then since asymptotically flat manifolds can be arbitrarily well approximated
by harmonically flat manifolds by lemma \ref{lem:sy}, inequality 
\ref{eqn:m>e} follows for asymptotically flat manifolds as well.

To prove the case of equality, we require a generalization of the case of 
equality of the positive mass theorem given in \cite{BF} as theorem 5.3.  In
that paper, we say that a singular manifold has 
{\it generalized nonnegative scalar curvature} 
if it is the limit (in the sense given in \cite{BF}) of
smooth manifolds with nonnegative scalar curvature.    

Note that if $(M^3,\bar{g})$ is only asymptotically flat and not harmonically
flat, then we have not shown that $(M^3,\tilde{g})$ can be extended smoothly
over the missing points $\{\infty_k\}$.  However, as a possibly singular
manifold, it does have generalized nonnegative scalar curvature since 
it is the limit of smooth manifolds with nonnegative scalar curvature 
(since $(M^3,\bar{g})$ 
can be arbitrarily well approximated by harmonically flat manifolds).

Then theorem 5.3 in \cite{BF} states that if $(M^3,\tilde{g})$ has generalized
nonnegative scalar curvature, zero mass, and positive isoperimetric constant,
then $(M^3,\tilde{g})$ is flat (outside the singular set).  Hence, it follows
that if we have equality in inequality \ref{eqn:m>e}, then $\tilde{g}$ is flat.
Hence, $(M^3,\tilde{g})$ is isometric to $(\real^3,\delta)$ minus a finite 
number of points, and since the harmonic conformal factor $\phi(x)$ preserves
the sign of the scalar curvature by equation \ref{eqn:scalar_curv}, 
the case of equality of the theorem follows.  \qed

\begin{definition}
Given a complete, asymptotically flat manifold $(M^3,g)$ with horizon $\Sigma
\in {\cal S}$ (defined in section \ref{sec:definitions}), define 
\begin{equation}
   {\cal E}(\Sigma,g) = \inf_{\varphi} \left\{ \frac{1}{2\pi}
          \int_{M^3} | \nabla \varphi |^2 \; dV \right\}
\end{equation}
where the infimum is taken over all smooth 
$\varphi(x)$ which go to one at infinity and equal zero on the horizon $\Sigma$
(and are zero inside $\Sigma$).
(By definition of ${\cal S}$, all of the ends other than the chosen end are
contained inside $\Sigma$.)  
\end{definition}
The infimum in the above definition
is achieved by the Green's function $\varphi(x)$ which satisfies
\begin{equation}\label{eqn:varphi}
\left\{\begin{array}{rll}
\lim_{x \rightarrow \infty_0} \varphi(x) &=& 1 \\
\Delta \varphi &=& 0 \\
\varphi(x) &=& 0 \;\mbox{ on }\; \Sigma
\end{array}\right.
\end{equation}
and as before,
\begin{equation}\label{eqn:phi_expansion}
   \varphi(x) = 1 - \frac{{\cal E}(\Sigma,g)}{2|x|}+ {\cal O}\left(\frac{1}{|x|^2}\right)
\end{equation}
in the chosen end.

\begin{theorem}\label{thm:green2}
Let $(M^3,g)$ be a complete, smooth, 
asymptotically flat 3-manifold with nonnegative 
scalar curvature with a horizon $\Sigma \in {\cal S}$ and total
mass $m$ (in the chosen end).  Then
\begin{equation}\label{eqn:horizon_energy}
   m \ge \frac12 {\cal E}(\Sigma,g)
\end{equation}
with equality if and only if $(M^3,g)$ is a Schwarzschild manifold
outside the horizon $\Sigma$.
\end{theorem} 
{\it Proof.}  
Let $M_\Sigma^3$ be the closed region of $M^3$ which is outside (or on)
$\Sigma$.  Since $\Sigma \in {\cal S}$, $(M_\Sigma^3,g)$ has only one end,
and we recall that 
$\Sigma$ could have multiple components.  For example, in the picture below 
$\Sigma$ has two components.
\vspace{.5in}
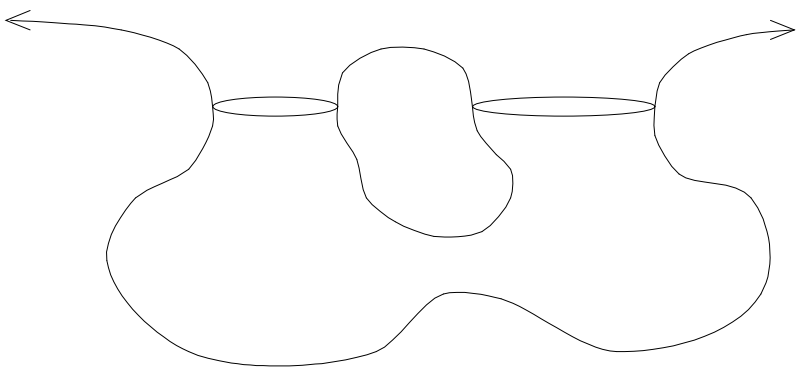
\begin{center}
\setlength{\unitlength}{0.00033333in}
\begingroup\makeatletter\ifx\SetFigFont\undefined%
\gdef\SetFigFont#1#2#3#4#5{%
  \reset@font\fontsize{#1}{#2pt}%
  \fontfamily{#3}\fontseries{#4}\fontshape{#5}%
  \selectfont}%
\fi\endgroup%
{\renewcommand{\dashlinestretch}{30}
\begin{picture}(12324,5591)(0,-10)
\put(4212,4064){\ellipse{1950}{300}}
\put(8712,4064){\ellipse{2850}{300}}
\path(11937,5414)(12312,5264)(11937,5114)
\path(387,5564)(12,5414)(387,5264)
\path(5187,4064)	(5186.505,4169.082)
	(5188.390,4259.815)
	(5200.181,4404.384)
	(5262.000,4589.000)

\path(5262,4589)	(5377.367,4710.573)
	(5538.855,4820.334)
	(5711.916,4908.178)
	(5862.000,4964.000)

\path(5862,4964)	(6013.256,4988.847)
	(6200.370,4995.489)
	(6387.049,4986.386)
	(6537.000,4964.000)

\path(6537,4964)	(6681.132,4921.820)
	(6852.596,4854.556)
	(7016.263,4767.014)
	(7137.000,4664.000)

\path(7137,4664)	(7190.028,4576.927)
	(7230.360,4457.278)
	(7261.512,4290.989)
	(7274.744,4185.961)
	(7287.000,4064.000)

\path(3237,4064)	(3218.898,4205.606)
	(3201.536,4309.449)
	(3162.000,4439.000)

\path(3162,4439)	(3080.314,4572.025)
	(2962.432,4722.927)
	(2831.835,4863.117)
	(2712.000,4964.000)

\path(2712,4964)	(2626.149,5012.720)
	(2519.340,5059.942)
	(2398.244,5104.811)
	(2269.534,5146.475)
	(2139.882,5184.079)
	(2015.961,5216.771)
	(1904.443,5243.695)
	(1812.000,5264.000)

\path(1812,5264)	(1694.303,5285.847)
	(1558.735,5305.930)
	(1400.240,5324.688)
	(1213.766,5342.562)
	(1108.458,5351.305)
	(994.259,5359.991)
	(870.538,5368.676)
	(736.665,5377.414)
	(592.005,5386.261)
	(435.929,5395.271)
	(267.805,5404.499)
	(87.000,5414.000)

\path(10137,4064)	(10153.313,4205.482)
	(10170.079,4309.284)
	(10212.000,4439.000)

\path(10212,4439)	(10295.701,4557.416)
	(10416.079,4686.474)
	(10546.917,4804.295)
	(10662.000,4889.000)

\path(10662,4889)	(10748.068,4934.738)
	(10855.033,4980.455)
	(10976.223,5024.983)
	(11104.965,5067.159)
	(11234.586,5105.814)
	(11358.414,5139.784)
	(11469.776,5167.901)
	(11562.000,5189.000)

\path(11562,5189)	(11671.979,5208.187)
	(11821.920,5225.911)
	(11917.369,5234.774)
	(12029.401,5243.930)
	(12160.212,5253.598)
	(12312.000,5264.000)

\path(3237,4064)	(3244.791,3952.269)
	(3247.387,3869.712)
	(3237.000,3764.000)

\path(3237,3764)	(3177.591,3600.235)
	(3133.528,3505.608)
	(3082.972,3408.616)
	(2971.618,3225.939)
	(2862.000,3089.000)

\path(2862,3089)	(2685.941,2973.673)
	(2571.078,2920.290)
	(2448.611,2867.855)
	(2326.256,2814.938)
	(2211.726,2760.109)
	(2037.000,2639.000)

\path(2037,2639)	(1961.752,2553.393)
	(1880.462,2449.572)
	(1798.866,2331.463)
	(1722.697,2202.991)
	(1657.692,2068.083)
	(1609.584,1930.665)
	(1584.109,1794.662)
	(1587.000,1664.000)

\path(1587,1664)	(1610.888,1551.230)
	(1647.928,1438.888)
	(1696.897,1327.559)
	(1756.571,1217.829)
	(1825.725,1110.284)
	(1903.135,1005.509)
	(1987.577,904.091)
	(2077.826,806.615)
	(2172.659,713.666)
	(2270.852,625.831)
	(2371.180,543.695)
	(2472.419,467.844)
	(2573.344,398.863)
	(2672.732,337.338)
	(2769.359,283.855)
	(2862.000,239.000)

\path(2862,239)	(3008.279,181.835)
	(3169.583,133.591)
	(3343.735,93.989)
	(3528.561,62.750)
	(3624.297,50.181)
	(3721.884,39.598)
	(3821.053,30.967)
	(3921.529,24.252)
	(4023.043,19.421)
	(4125.320,16.436)
	(4228.090,15.265)
	(4331.081,15.871)
	(4434.021,18.221)
	(4536.637,22.279)
	(4638.657,28.011)
	(4739.810,35.381)
	(4839.824,44.356)
	(4938.427,54.900)
	(5035.347,66.978)
	(5130.311,80.556)
	(5313.286,112.073)
	(5485.176,149.171)
	(5643.806,191.573)
	(5787.000,239.000)

\path(5787,239)	(5915.083,309.043)
	(6044.385,416.479)
	(6174.750,548.675)
	(6306.019,693.001)
	(6438.035,836.828)
	(6570.643,967.523)
	(6703.683,1072.458)
	(6837.000,1139.000)

\path(6837,1139)	(6933.912,1159.563)
	(7046.633,1166.576)
	(7169.483,1162.538)
	(7296.788,1149.946)
	(7422.868,1131.299)
	(7542.049,1109.093)
	(7648.652,1085.828)
	(7737.000,1064.000)

\path(7737,1064)	(7922.038,995.058)
	(8024.804,944.980)
	(8133.177,886.947)
	(8246.240,822.796)
	(8363.075,754.367)
	(8482.764,683.495)
	(8604.390,612.020)
	(8727.035,541.778)
	(8849.781,474.608)
	(8971.712,412.347)
	(9091.908,356.833)
	(9209.452,309.904)
	(9323.428,273.397)
	(9432.916,249.150)
	(9537.000,239.000)

\path(9537,239)	(9666.320,238.729)
	(9805.659,244.258)
	(9953.193,255.959)
	(10107.098,274.204)
	(10265.550,299.365)
	(10426.725,331.814)
	(10588.799,371.924)
	(10749.949,420.065)
	(10908.350,476.611)
	(11062.179,541.933)
	(11209.612,616.403)
	(11348.825,700.393)
	(11477.994,794.276)
	(11595.296,898.424)
	(11698.906,1013.208)
	(11787.000,1139.000)

\path(11787,1139)	(11868.707,1313.366)
	(11896.885,1407.433)
	(11917.005,1504.965)
	(11929.362,1605.120)
	(11934.249,1707.056)
	(11931.961,1809.929)
	(11922.791,1912.899)
	(11907.034,2015.121)
	(11884.984,2115.755)
	(11856.935,2213.956)
	(11823.181,2308.883)
	(11739.736,2485.545)
	(11637.000,2639.000)

\path(11637,2639)	(11531.394,2731.363)
	(11396.558,2794.292)
	(11242.105,2836.006)
	(11077.646,2864.724)
	(10912.794,2888.664)
	(10757.161,2916.044)
	(10620.359,2955.083)
	(10512.000,3014.000)

\path(10512,3014)	(10399.202,3131.983)
	(10288.466,3293.180)
	(10195.747,3464.787)
	(10137.000,3614.000)

\path(10137,3614)	(10119.615,3773.052)
	(10123.961,3896.766)
	(10137.000,4064.000)

\path(5187,4064)	(5179.001,3952.325)
	(5176.335,3869.788)
	(5187.000,3764.000)

\path(5187,3764)	(5242.613,3639.642)
	(5331.638,3496.554)
	(5423.343,3355.939)
	(5487.000,3239.000)

\path(5487,3239)	(5521.312,3107.631)
	(5550.709,2938.951)
	(5585.751,2770.295)
	(5637.000,2639.000)

\path(5637,2639)	(5724.036,2536.793)
	(5847.045,2430.073)
	(5977.532,2334.066)
	(6087.000,2264.000)

\path(6087,2264)	(6214.154,2199.651)
	(6378.217,2130.564)
	(6546.672,2071.944)
	(6687.000,2039.000)

\path(6687,2039)	(6855.437,2027.665)
	(6956.757,2027.794)
	(7062.900,2032.707)
	(7168.930,2043.095)
	(7269.913,2059.647)
	(7437.000,2114.000)

\path(7437,2114)	(7567.603,2208.211)
	(7699.879,2344.854)
	(7813.215,2497.320)
	(7887.000,2639.000)

\path(7887,2639)	(7910.785,2740.192)
	(7919.752,2863.419)
	(7912.343,2986.936)
	(7887.000,3089.000)

\path(7887,3089)	(7780.380,3210.901)
	(7662.000,3314.000)

\path(7662,3314)	(7503.052,3492.451)
	(7419.405,3597.215)
	(7362.000,3689.000)

\path(7362,3689)	(7323.476,3818.491)
	(7305.861,3922.349)
	(7287.000,4064.000)

\end{picture}
}
\end{center}

Then the basic idea is to reflect $(M_\Sigma^3,g)$ through $\Sigma$ to 
get a manifold $(\bar{M}_\Sigma^3,\bar{g})$ with two asymptotically flat ends.
Then define $\phi(x)$ on $(\bar{M}_\Sigma^3,\bar{g})$ using 
equation \ref{eqn:phi} and $\varphi(x)$ on $(M_\Sigma^3,g)$ using equation
\ref{eqn:varphi}.  It follows from symmetry that $\phi(x) = \frac12$ on $\Sigma$,
so that 
\begin{equation}
  \phi(x) = \frac12 (\varphi(x) + 1) 
\end{equation}
on $(M_\Sigma^3,g)$.  Then 
\begin{equation}\label{eqn:symenergy}
     {\cal E}(\bar{g}) = \frac12 {\cal E}(\Sigma,g) 
\end{equation}
so that
theorem \ref{thm:green2} follows from theorem \ref{thm:green1}.

The only technicality is that theorem \ref{thm:green1} applies to smooth
manifolds with nonnegative scalar curvature, and $(\bar{M}_\Sigma^3,\bar{g})$
is typically not smooth along $\Sigma$, which also makes it unclear 
how to define the 
scalar curvature there.  However, it happens that because $\Sigma$ has zero
mean curvature, these issues can be resolved.

This idea of reflecting a manifold through its horizon is used 
by Bunting and Masood-ul-Alam in 
\cite{BM}, and the issue of the smoothness of the reflected manifold appears
in their paper as well.  However, in their setting they have the simpler
case in which the horizon not only has zero mean curvature but also has 
zero second fundamental form.  Hence, the reflected manifold is 
$C^{1,1}$, which apparently is sufficient for their purposes.

However, in our setting we can not assume that the horizon $\Sigma$ has
zero second fundamental form, so that $(\bar{M}_\Sigma^3,\bar{g})$ is only
Lipschitz.  To solve this problem, given $\delta > 0$ 
we will define a smooth manifold 
$(\tilde{M}^3_{\Sigma,\delta},\tilde{g}_\delta)$ 
with nonnegative scalar curvature 
which, in the limit as $\delta$ approaches zero, 
approaches $(\bar{M}_\Sigma^3,\bar{g})$ (meaning that there 
exists a diffeomorphism under which the metrics 
are arbitrarily uniformly close to each other and the total masses are 
arbitrarily close).  
Then by definition \ref{def:energy1} it follows that 
${\cal E}(\tilde{g_\delta})$ 
is close to ${\cal E}(\bar{g})$, from which we will
be able to conclude
\begin{equation}\label{eqn:approx}
   m \approx \tilde{m}_\delta \ge {\cal E}(\tilde{g}_\delta) 
     \approx {\cal E}(\bar{g})
     = \frac12 {\cal E}(\Sigma,g),
\end{equation}
where $\tilde{m}_\delta$ is the mass of 
$(\tilde{M}^3_{\Sigma,\delta},\tilde{g}_\delta)$ and the approximations in 
the above inequality can be made to be arbitrarily accurate by choosing 
$\delta$ small,
thereby proving inequality \ref{eqn:horizon_energy}.

The first step is to construct the smooth manifolds
\begin{equation}
(\tilde{M}^3_{\Sigma,\delta},\bar{g}_\delta)
\approx (M^3_\Sigma, g) \cup (\Sigma \times [0,2\delta], G) 
\cup (M^3_\Sigma, g), 
\end{equation}
where identifications are made along the boundaries
of these three manifolds as drawn below.  (To be precise, the second 
$(M^3_\Sigma, g)$ in the above union is meant to be a copy of the first
$(M^3_\Sigma, g)$ and therefore distinct.)  We will define the metric
$G$ such that the metric $\bar{g}_\delta$ is smooth, although it will not
have nonnegative scalar curvature.  Then we will define
\begin{equation}\label{eqn:gtilde}
   \tilde{g}_\delta = u_\delta(x)^4 \bar{g}_\delta
\end{equation}
so that $\tilde{g}_\delta$ is not only smooth but also has nonnegative scalar 
curvature, and we will show that because of our choice of 
the metric $G$, $u_\delta(x)$
approaches one in the limit as $\delta$ approaches zero.
\vspace{.5in}
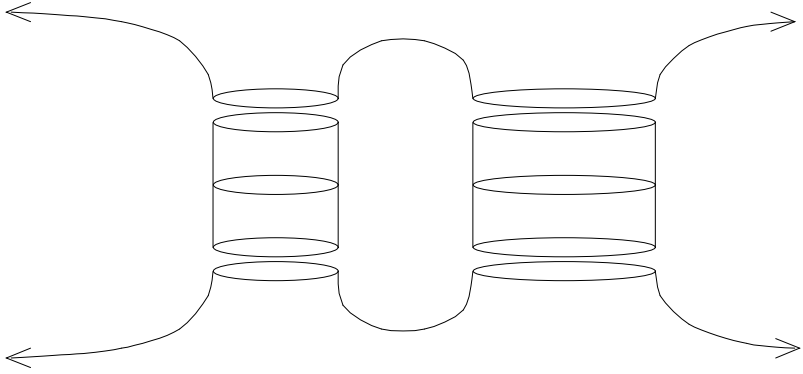
\begin{center}
\setlength{\unitlength}{0.00033333in}
\begingroup\makeatletter\ifx\SetFigFont\undefined%
\gdef\SetFigFont#1#2#3#4#5{%
  \reset@font\fontsize{#1}{#2pt}%
  \fontfamily{#3}\fontseries{#4}\fontshape{#5}%
  \selectfont}%
\fi\endgroup%
{\renewcommand{\dashlinestretch}{30}
\begin{picture}(12399,5739)(0,-10)
\put(4212,1512){\ellipse{1950}{300}}
\put(8712,1512){\ellipse{2850}{300}}
\put(4212,1887){\ellipse{1950}{300}}
\put(4212,2862){\ellipse{1950}{300}}
\put(4212,3837){\ellipse{1950}{300}}
\put(8712,1887){\ellipse{2850}{300}}
\put(8712,2862){\ellipse{2850}{300}}
\put(8712,3837){\ellipse{2850}{300}}
\put(4212,4212){\ellipse{1950}{300}}
\put(8712,4212){\ellipse{2850}{300}}
\path(387,312)(12,162)(387,12)
\path(12012,462)(12387,312)(12012,162)
\path(3237,3837)(3237,1887)
\path(5187,3837)(5187,1887)
\path(7287,3837)(7287,1887)
\path(10137,3837)(10137,1887)
\path(11937,5562)(12312,5412)(11937,5262)
\path(387,5712)(12,5562)(387,5412)
\path(3237,1512)	(3218.898,1370.394)
	(3201.536,1266.551)
	(3162.000,1137.000)

\path(3162,1137)	(3080.314,1003.975)
	(2962.432,853.073)
	(2831.835,712.883)
	(2712.000,612.000)

\path(2712,612)	(2626.149,563.280)
	(2519.340,516.058)
	(2398.244,471.189)
	(2269.534,429.525)
	(2139.882,391.921)
	(2015.961,359.229)
	(1904.443,332.305)
	(1812.000,312.000)

\path(1812,312)	(1694.303,290.153)
	(1558.735,270.070)
	(1400.240,251.312)
	(1213.766,233.438)
	(1108.458,224.695)
	(994.259,216.009)
	(870.538,207.324)
	(736.665,198.586)
	(592.005,189.739)
	(435.929,180.729)
	(267.805,171.501)
	(87.000,162.000)

\path(5187,1512)	(5186.505,1406.918)
	(5188.390,1316.185)
	(5200.181,1171.616)
	(5262.000,987.000)

\path(5262,987)	(5377.367,865.427)
	(5538.855,755.666)
	(5711.916,667.822)
	(5862.000,612.000)

\path(5862,612)	(6013.256,587.153)
	(6200.370,580.511)
	(6387.049,589.614)
	(6537.000,612.000)

\path(6537,612)	(6681.132,654.180)
	(6852.596,721.444)
	(7016.263,808.986)
	(7137.000,912.000)

\path(7137,912)	(7190.028,999.073)
	(7230.360,1118.723)
	(7261.512,1285.011)
	(7274.744,1390.039)
	(7287.000,1512.000)

\path(10137,1512)	(10153.313,1370.518)
	(10170.079,1266.716)
	(10212.000,1137.000)

\path(10212,1137)	(10295.701,1018.584)
	(10416.079,889.526)
	(10546.917,771.705)
	(10662.000,687.000)

\path(10662,687)	(10748.068,641.262)
	(10855.033,595.545)
	(10976.223,551.017)
	(11104.965,508.841)
	(11234.586,470.186)
	(11358.414,436.216)
	(11469.776,408.099)
	(11562.000,387.000)

\path(11562,387)	(11671.979,367.813)
	(11821.920,350.089)
	(11917.369,341.226)
	(12029.401,332.070)
	(12160.212,322.402)
	(12312.000,312.000)

\path(5187,4212)	(5186.505,4317.082)
	(5188.390,4407.815)
	(5200.181,4552.384)
	(5262.000,4737.000)

\path(5262,4737)	(5377.367,4858.573)
	(5538.855,4968.334)
	(5711.916,5056.178)
	(5862.000,5112.000)

\path(5862,5112)	(6013.256,5136.847)
	(6200.370,5143.489)
	(6387.049,5134.386)
	(6537.000,5112.000)

\path(6537,5112)	(6681.132,5069.820)
	(6852.596,5002.556)
	(7016.263,4915.014)
	(7137.000,4812.000)

\path(7137,4812)	(7190.028,4724.927)
	(7230.360,4605.278)
	(7261.512,4438.989)
	(7274.744,4333.961)
	(7287.000,4212.000)

\path(3237,4212)	(3218.898,4353.606)
	(3201.536,4457.449)
	(3162.000,4587.000)

\path(3162,4587)	(3080.314,4720.025)
	(2962.432,4870.927)
	(2831.835,5011.117)
	(2712.000,5112.000)

\path(2712,5112)	(2626.149,5160.720)
	(2519.340,5207.942)
	(2398.244,5252.811)
	(2269.534,5294.475)
	(2139.882,5332.079)
	(2015.961,5364.771)
	(1904.443,5391.695)
	(1812.000,5412.000)

\path(1812,5412)	(1694.303,5433.847)
	(1558.735,5453.930)
	(1400.240,5472.688)
	(1213.766,5490.562)
	(1108.458,5499.305)
	(994.259,5507.991)
	(870.538,5516.676)
	(736.665,5525.414)
	(592.005,5534.261)
	(435.929,5543.271)
	(267.805,5552.499)
	(87.000,5562.000)

\path(10137,4212)	(10153.313,4353.482)
	(10170.079,4457.284)
	(10212.000,4587.000)

\path(10212,4587)	(10295.701,4705.416)
	(10416.079,4834.474)
	(10546.917,4952.295)
	(10662.000,5037.000)

\path(10662,5037)	(10748.068,5082.738)
	(10855.033,5128.455)
	(10976.223,5172.983)
	(11104.965,5215.159)
	(11234.586,5253.814)
	(11358.414,5287.784)
	(11469.776,5315.901)
	(11562.000,5337.000)

\path(11562,5337)	(11671.979,5356.187)
	(11821.920,5373.911)
	(11917.369,5382.774)
	(12029.401,5391.930)
	(12160.212,5401.598)
	(12312.000,5412.000)

\end{picture}
}
\end{center}

We will use the local coordinates $(z,t)$ to describe points on 
$\Sigma \times [0,2\delta]$, where 
$z = (z_1, z_2) 
\in \mbox{a local coordinate chart for }\Sigma$ and $t \in [0,2\delta]$.
Then we define $G(\partial_t, \partial_t) = 1$, 
$G(\partial_t, \partial_{z_1}) = 0$, and 
$G(\partial_t, \partial_{z_2}) = 0$.  Then
it follows that $\Sigma \times t$ is obtained by flowing 
$\Sigma \times 0$ in the unit normal direction for a time $t$, and that 
$\partial_t$ is orthogonal to $\Sigma \times t$.  Hence, all that remains
to fully define the metric $G$ is to define it smoothly on the tangent planes 
of $\Sigma \times t$ for $0 \le t \le 2\delta$.

Let $\bar{G}(z,t)$ be the metric $G$ restricted to $\Sigma \times t$.
Then
\begin{equation}\label{eqn:GODE}
   \frac{d}{dt} \bar{G}_{ij}(z,t) = 2 \bar{G}_{ik}(z,t) h^k_j(z,t)  
\end{equation}
where $h^k_j(z,t)$ is the second fundamental form of $\Sigma \times t$ in
$(\Sigma \times [0,2\delta], G)$ with respect to the normal vector 
$\partial_t$.  
Furthermore, since $\Sigma \times 0$ is
identified with $\Sigma \in (M^3_\Sigma,g)$, we can extend the coordinates
$(z,t)$ for $t$ slightly less than zero into $(M^3_\Sigma,g)$, thereby
giving us smooth initial data for $\bar{G}_{ij}(z,t)$ and $h^k_j(z,t)$ for 
$-\epsilon < t \le 0$, for some positive $\epsilon$.

Now we extend $h^k_j(z,t)$ smoothly for $0 \le t \le 2\delta$ 
in such a way that $h^k_j(z,t)$
is an odd function about $t = \delta$, 
meaning that $h^k_j(z,t) = - h^k_j(z, 2\delta-t)$.  
Naturally there are many ways 
to accomplish this smooth extension.  

Then we define $\bar{G}_{ij}(z,t)$ to be the 
smooth solution to the o.d.e.~given in equation \ref{eqn:GODE} using the
initial data for $\bar{G}_{ij}(z,t)$ at $t=0$.  
By the oddness of $h^k_j(z,t)$ about $t=\delta$ 
it follows that $\bar{G}_{ij}(z,t)$
is symmetric about $t=\delta$, that is, 
$\bar{G}_{ij}(z,t) = \bar{G}_{ij}(z,2\delta-t)$.
Hence, the identification of $\Sigma \times (2\delta)$ with 
$\Sigma \in \mbox{the second copy of }(M^3_\Sigma,g)$ is smooth by 
symmetry.  This completes the smooth construction of the metric 
$(\tilde{M}^3_{\Sigma,\delta},\bar{g}_\delta)$.

Now define $H(z,t) = \Sigma_j h^j_j(z,t)$ to be the mean curvature of 
$\Sigma \times t$ in $(\Sigma \times [0,2\delta], G)$, and let 
$\dot H(z,t) = \frac{d}{dt} H(z,t)$.  We note that 
\begin{equation}\label{eqn:zeroH}
   H(z,0) = 0 = H(z,2\delta)
\end{equation}
since $\Sigma$ is a horizon (and hence has zero mean curvature) in 
$(M^3_\Sigma,g)$.   
Let $\alpha = \sup_z \dot H(z,0)$ and let
$\beta = \sup_z \sum_{jk} h^k_j(z,0) h^j_k(z,0)$, 
which we note are functions of the 
metric $g$ on $M^3_\Sigma$ and are independent of $\delta$.  
Then we require that the smooth extension
we choose for $h^k_j$ satisfies
\begin{equation}\label{eqn:2alpha}
   \dot H(z,t) \le 2|\alpha| + 1,
\end{equation}
(which is possible because of equation \ref{eqn:zeroH}) and
\begin{equation}\label{eqn:2beta}
   \sum_{jk} h^k_j(z,t) h^j_k(z,t) \le 2\beta + 1 .
\end{equation}
Then combining equations
\ref{eqn:zeroH} and \ref{eqn:2alpha} also yields
\begin{equation}\label{eqn:2alphadelta}
   |H(z,t)| \le (2|\alpha| + 1) \delta.
\end{equation}

These estimates allow us to bound the scalar curvature of 
$(\tilde{M}^3_{\Sigma,\delta},\bar{g}_\delta)$ from below since 
by the second variation formula and the Gauss equation
(see equations \ref{eqn:2var} and \ref{eqn:Gauss}) we have that
\begin{equation}
  R = -2\dot{H} + 2K - |h|^2 - H^2
\end{equation}
where $R(z,t)$ is scalar curvature and $K(z,t)$ is the Gauss curvature of
$\Sigma \times t$.  At this point we realize that we also need a lower bound
$K_0$ for $K(z,t)$ which is independent of $\delta$, which follows from 
imposing an upper bound on the $C^2$ norm (in the $z$ variable) of
our smooth choice of $h^k_j(z,t)$.  Then using this combined with 
inequalities \ref{eqn:2alpha}, \ref{eqn:2beta},
and \ref{eqn:2alphadelta} we get
\begin{equation}\label{eqn:r>r0} 
  R(z,t) \ge R_0
\end{equation}
where $R_0$ is independent of $\delta$ (for $\delta < 1$). 

Now we are ready to define $\tilde{g}_\delta$ using equation 
\ref{eqn:gtilde}.  We already know that 
$(\tilde{M}^3_{\Sigma,\delta},\bar{g}_\delta)$ is smooth and has nonnegative
scalar curvature everywhere except possibly in $\Sigma \times [0,2\delta]$
where it has $R \ge R_0$.  If $R_0 \ge 0$, then we just let $u_\delta(x) = 1$
so that $\tilde{g} = \bar{g}$.  Otherwise, we define $u_\delta(x)$ such that
\begin{equation}\label{eqn:udelta}
   ( - 8 \Delta_{\bar{g}} + {\cal R}_\delta(x) ) u_\delta(x) = 0   
\end{equation}
and $u_\delta(x)$ goes to one in both asymptotically flat ends, where
${\cal R}_\delta(x)$ equals $R_0$ in $\Sigma \times [0,2\delta]$,  
equals zero for $x$ more than a distance $\delta$ from    
$\Sigma \times [0,2\delta]$, is smooth, 
and takes values in $[R_0, 0]$ everywhere.  Then it follows that 
for sufficiently small $\delta$, $u_\delta(x)$
is a smooth superharmonic function.  Furthermore, since 
${\cal R}_\delta$ is zero 
everywhere except on an open set whose volume is going to zero as $\delta$
goes to zero, and since ${\cal R}_\delta$ is uniformly bounded from below
on this small set, it follows from bounding Green's functions from above that 
\begin{equation}\label{eqn:udelta_bound}
   1 \le u_\delta(x) \le 1 + \epsilon(\delta)
\end{equation}
where $\epsilon$ goes to zero as $\delta$ approaches zero.

Furthermore, by equations \ref{eqn:r>r0}, \ref{eqn:gtilde}, 
and \ref{eqn:scalar_curv}, $(\tilde{M}^3_{\Sigma,\delta},\tilde{g}_\delta)$
has nonnegative scalar curvature, and since $u_\delta(x)$ and 
and $(\tilde{M}^3_{\Sigma,\delta},\bar{g}_\delta)$ are smooth, so is
$(\tilde{M}^3_{\Sigma,\delta},\tilde{g}_\delta)$.  In addition, it follows 
from the construction of $(\tilde{M}^3_{\Sigma,\delta},\bar{g}_\delta)$
that there exists a diffeomorphism into $(\bar{M}^3,\bar{g})$ with respect
to which the metrics are arbitrarily uniformly close to each other in the
limit as $\delta$ goes to zero.  Hence, by equation \ref{eqn:udelta_bound},
we see that the same statement is true for 
$(\tilde{M}^3_{\Sigma,\delta},\tilde{g}_\delta)$.  
Finally, it follows from equation \ref{eqn:udelta} that 
$\tilde{m}_\delta$, the mass 
of $(\tilde{M}^3_{\Sigma,\delta},\tilde{g}_\delta)$,
converges to $m$, the mass of $(\bar{M}^3_\Sigma, \bar{g})$, in the limit
as $\delta$ goes to zero.  Hence, inequality \ref{eqn:approx} follows, 
proving inequality \ref{eqn:horizon_energy}.

To prove the case of equality, we note that we can view the above proof
in a different way.  Since the singular manifold 
$(\bar{M}^3_{\Sigma},\bar{g})$ is the limit of the smooth manifolds
$(\tilde{M}^3_{\Sigma,\delta},\tilde{g}_\delta)$
which have nonnegative scalar curvature, it follows that  
$(\bar{M}^3_{\Sigma},\bar{g})$ has generalized nonnegative scalar curvature
as defined in \cite{BF}.  Then if we reexamine the proof of theorem
\ref{thm:green1}, we see that the theorem, including the case of equality, is
also true for singular manifolds like $(\bar{M}^3_{\Sigma},\bar{g})$ which 
have generalized nonnegative scalar curvature (see the discussion at the 
end of the proof of theorem \ref{thm:green1}).  Hence, by equation 
\ref{eqn:symenergy} and theorem \ref{thm:green1}, 
we get equality in inequality \ref{eqn:horizon_energy} if and only if
$(\bar{M}^3_{\Sigma},\bar{g})$ has zero scalar curvature and is conformal to 
$(\real^3,\delta)$ minus a finite number of points.

Since $(\bar{M}^3_{\Sigma},\bar{g})$ has two ends, it must be conformal to 
$(\real^3 \backslash \{0\},\delta)$, and since it has zero scalar curvature,
it follows from equation \ref{eqn:scalar_curv} that 
it must be a Schwarzschild metric.  Hence, in the case of equality for 
inequality \ref{eqn:horizon_energy}, 
$(M^3,g)$ must be a Schwarzschild manifold outside
$\Sigma$.  \qed

\section{Proof That $m(t)$ Is Nonincreasing}
\label{sec:m(t)}

In this section we will finish the proof of theorem \ref{thm:monotone}
begun in section \ref{sec:A(t)} by proving that $m(t)$, the total mass of 
$(M^3,g_t)$, is non-increasing in $t$.  The fact that $m(t)$ is 
nonincreasing is of course central to the argument presented in this 
paper for proving the Riemannian Penrose conjecture and is perhaps the 
most important property of the conformal flow of metrics $\{g_t\}$.

We begin with a corollary to lemma \ref{lem:integral} 
and theorem \ref{thm:equal} in section 
\ref{sec:existence}.
\begin{corollary}
The left and right hand derivatives $\frac{d}{dt^\pm}$ 
of $u_t(x)$ 
exist for all $t>0$ and are equal except at a countable number of 
$t$-values.  Furthermore,
\begin{equation}
\frac{d}{dt^+}\, u_t(x) = v_t^+(x)
\end{equation}
and 
\begin{equation}
\frac{d}{dt^-}\, u_t(x) = v_t^-(x)
\end{equation}
where $v_t^{\pm}(x)$ equals zero inside $\Sigma^{\pm}(t)$ 
(see definition \ref{def:limits}) and outside 
$\Sigma^{\pm}(t)$ is the harmonic function which equals $0$ on 
$\Sigma^{\pm}(t)$ and goes to $-e^{-t}$ at infinity.
\end{corollary}

We will use this corollary to compute the left and right hand derivatives
of $m(t)$.  As proven at the end of appendix \ref{sec:harmonic}, the flow
of metrics $\{g_t\}$ we are considering has the property that the rate of 
change of the metric $g_t$ is just a function of $g_t$ and not of $t$ or 
$g_0$.  Hence, we will just prove that $m'(0) \le 0$, from which it will 
follow that $m'(t) \le 0$.  So without loss of generality, we will assume 
that the flow begins at some time $-t_0 < 0$, and then compute the left and 
right hand derivatives of $m(t)$ at $t = 0$.

Also, we remind the reader that we proved that $\Sigma^+(t)$ and $\Sigma^-(t)$
are horizons in $(M^3,g_t)$ in lemma \ref{lem:C1}.  Furthermore, 
$v_0^\pm(x)$ is harmonic in
$(M^3,g_0)$, equals $0$ on $\Sigma^\pm(0)$, and goes to $-1$ at infinity.
Hence, by 
equation \ref{eqn:phi_expansion} and theorem \ref{thm:green2} of the 
previous section,  
\begin{equation}\label{eqn:v_expansion}
   v_0^\pm(x) = -1 + \frac{{\cal E}(\Sigma^\pm(0),g_0)}{2|x|}+ 
                {\cal O}\left(\frac{1}{|x|^2}\right)
\end{equation}
where 
\begin{equation}\label{eqn:m(0)lower}
   m(0) \ge \frac12 {\cal E}(\Sigma^\pm(0),g_0).
\end{equation}

Now we are ready to compute $m'(t)$.  As in section \ref{sec:definitions}, 
let $g_0 = {\cal U}_0(x)^4 g_{flat}$, where
$(M^3,g_{flat})$ is isometric to $(\real^3 \backslash B_{r_0}(0), \delta)$
in the harmonically flat end, where we have chosen $r_0$ and scaled  
the harmonically flat coordinate chart such that ${\cal U}_0(x)$
goes to one at infinity.  Then by definition \ref{def:totalmass} for the total
mass, we have that
\begin{equation}
   {\cal U}_0(x) = 1 + \frac{m(0)}{2|x|} +{\cal O}\left(\frac{1}{|x|^2}\right).
\end{equation}
We will also let $g_t = {\cal U}_t(x)^4 g_{flat}$ in the harmonically flat end.
Then since $g_t = u_t(x)^4 g_0$, it follows that 
\begin{equation}\label{eqn:mult_u}
   {\cal U}_t(x) = u_t(x) \,{\cal U}_0(x).
\end{equation}
Now we define $\alpha(t)$ and $\beta(t)$ such that 
\begin{equation}
   u_t(x) = \alpha(t) + \frac{\beta(t)}{|x|} 
            +{\cal O}\left(\frac{1}{|x|^2}\right).
\end{equation}
Then since $u_0(x) \equiv 1$ and 
$\frac{d}{dt^\pm}\, u_t(x) |_{t=0} = v_0^\pm(x)$, it follows from equation
\ref{eqn:v_expansion} that 
\begin{equation}\label{eqn:alphabeta}\begin{array}{ll}
\alpha(0) = 1, \hspace{.5in} & \frac{d}{dt^\pm}\,\alpha(t) |_{t=0} = -1, \\
& \\
\beta(0) = 0,  \hspace{.5in} & \frac{d}{dt^\pm}\,\beta(t) |_{t=0} =
                \frac12 {\cal E}(\Sigma^\pm(0),g_0).
\end{array}\end{equation}
Thus, by equation \ref{eqn:mult_u},
\begin{equation}
   {\cal U}_t(x) = \alpha(t) + \frac{1}{|x|}\left(\beta(t) + 
                   \frac{m(0)}{2}\alpha(t)\right)
                   +{\cal O}\left(\frac{1}{|x|^2}\right)
\end{equation}
so that by definition \ref{def:totalmass}
\begin{equation}
   m(t) = 2 \alpha(t)\left(\beta(t) + \frac{m(0)}{2}\alpha(t)\right). 
\end{equation}
Hence, by equation \ref{eqn:alphabeta}
\begin{equation}\label{eqn:mprime}
   \frac{d}{dt^\pm}\, m(t)|_{t=0} \;\; = \;\; 
   {\cal E}(\Sigma^\pm(0),g_0) - 2m(0) \;\; \le \;\; 0 
\end{equation}
by equation \ref{eqn:m(0)lower}.  Then since we were able to choose 
$t=0$ without loss of generality as previously discussed, we have proven
the following theorem.

\begin{theorem}\label{thm:massnonincreasing}
The left and right hand derivatives $\frac{d}{dt^\pm}$ of $m(t)$  
exist for all $t>0$ and are equal except at a countable number of 
$t$-values.  Furthermore, 
\begin{equation}
\frac{d}{dt^\pm}\, m(t) \le 0
\end{equation}
for all $t > 0$ (and the right hand derivative of $m(t)$ at $t=0$ exists
and is nonpositive as well).
\end{theorem}

Hence, $m'(t)$ exists almost everywhere, and since the 
left and right hand derivatives of $m(t)$ are all nonpositive,
$m(t)$ is nonincreasing.  Since we proved that $A(t)$ 
is constant in section \ref{sec:A(t)}, this completes the proof
of theorem \ref{thm:monotone}.

\section{The Stability of $\Sigma(t)$}
\label{sec:stability}

A very interesting and important property of the horizons $\Sigma(t)$
follows from the fact that they are locally stable with respect to outward
variations.
Since $\Sigma(t)$ is strictly outer-minimizing in $(M^3,g_t)$, we know that
each component of $\Sigma(t)$ must be stable under outward variations (holding
$t$ fixed).
In particular, choose any component $\Sigma_i(t)$ of $\Sigma(t)$, and flow it 
outwards in the unit normal direction at constant speed one.  
Since $\Sigma(t)$ is smooth, we can do this for some positive 
amount of time without the surface forming 
singularities.
Let $A(s)$ be the area of the surface after being flowed at speed one for 
time $s$.  Then since
$\Sigma_i(t)$ has zero mean curvature, $A'(0) = 0$.  And since $\Sigma(t)$ is
strictly outer-minimizing, which we recall 
means that all surfaces in $(M^3,g_t)$ 
which enclose $\Sigma(t)$ 
have strictly larger area, we must have $A''(0) \ge 0$.   

On the other hand, 
\begin{equation}
A'(s) = \int H d\mu
\end{equation}
where the integral is being taken over the surface resulting from 
flowing $\Sigma_i(t)$ for time $s$, $H$ is the mean curvature of the
surface, and $d\mu$ is the area form of the surface.
Then since $H \equiv 0$ at $s=0$, 
\begin{equation}
A''(0) = \int_{\Sigma_i(t)} \frac{d}{ds} (H) d\mu.
\end{equation}
We will then use the second variation formula
\begin{equation}\label{eqn:2var}
\frac{d}{dt} H = - |h|^2 - Ric(\vec{\nu},\vec{\nu}),
\end{equation}
the Gauss equation,
\begin{equation}\label{eqn:Gauss}
Ric(\vec{\nu},\vec{\nu}) = \frac12 R - K + \frac12 H^2 - \frac12 |h|^2,
\end{equation}
and the fact that $|h|^2 = \frac12 (\lambda_1-\lambda_2)^2 + \frac12 H^2$, 
where $h$ is the second fundamental form of $\Sigma_i(t)$
(so that $H = \mbox{trace}(h)$), $Ric$ is the Ricci curvature tensor
of $(M^3,g_t)$,
$\vec{\nu}$ is the outward pointing normal vector to $\Sigma_i(t)$, $R$ is the 
scalar curvature of $(M^3,g_t)$, 
$K$ is the Gauss curvature of $\Sigma_i(t)$, and 
$\lambda_1$ and $\lambda_2$ are the principal curvatures of $\Sigma_i(t)$, 
to get
\begin{equation}\label{eqn:App}
A''(0) =  \int_{\Sigma_i(t)} -\frac12 R + K - \frac14 (\lambda_1-\lambda_2)^2
\end{equation}
since $H=0$.
Hence, since $A''(0) \ge 0$, $R \ge 0$, and $\int_{\Sigma_i(t)} K \le 4\pi$ 
by the Gauss-Bonnet formula (actually we have equality in the last inequality 
since $\Sigma_i(t)$ is a sphere since $\Sigma(t)$ is strictly outer-minimizing
in $(M^3,g_t)$ - see below), we conclude that
\begin{equation}\label{eqn:Willmore-like_bound}
\int_{\Sigma_i(t)} (\lambda_1-\lambda_2)^2 d\mu  \le 16\pi 
\end{equation}
with respect to the metric $g_t$.  However, it happens that the left hand
side of inequality \ref{eqn:Willmore-like_bound} is 
conformally invariant, and $g_t$ is conformal to $g_{flat}$
(defined in section \ref{sec:definitions}). 
Thus, inequality \ref{eqn:Willmore-like_bound} is also 
true with respect to the metric $g_{flat}$. 
\begin{theorem}\label{thm:stability}
Each component $\Sigma_i(t)$ of the surface $\Sigma(t)$ (which is a strictly 
outer-minimizing horizon in $(M^3,g_t)$) is a sphere and satisfies
\begin{equation}\label{eqn:stability}
\int_{\Sigma_i(t)} (\lambda_1-\lambda_2)^2 d\mu  \le 16\pi 
\end{equation}
with respect to the fixed metric $g_{flat}$.
\end{theorem}

We also comment that 
equation \ref{eqn:App} can be used to prove that each component
of a strictly outer-minimizing surface $\Sigma$ in a manifold
$(M^3,g)$ with nonnegative scalar curvature is a sphere.
First, we choose a superharmonic function $u(x)$ defined 
outside $\Sigma$ in $(M^3,g)$ 
approximately equal to one
which has negative outward derivative on $\Sigma$.  Then the metric
$(M^3,u(x)^4 g)$ has {\it strictly} positive scalar curvature by equation
\ref{eqn:scalar_curv}.  Also, since $\Sigma$ now has negative mean curvature 
in $(M^3,u(x)^4 g)$ by the Neumann condition on $u(x)$,  $\Sigma$
acts a barrier so that the outermost minimal area enclosure of $\Sigma$ does
not touch $\Sigma$ and hence is a stable horizon with zero mean curvature.  
Then by equation \ref{eqn:App} and the Gauss-Bonnet formula, 
each component of this outermost minimal
area enclosure of $\Sigma$ must be
a sphere.  In the limit as $u(x)$ approaches one, the area of this outermost
minimal area enclosure in $(M^3,g)$ must approach the area of $\Sigma$.  Then 
since $\Sigma$ is {\it strictly} outer-minimizing in $(M^3,g)$, 
these outermost minimal area enclosures must be approaching $\Sigma$ in the
limit as $u(x)$ goes to one, 
proving that each component of $\Sigma$ is a sphere.  

Unfortunately, this argument doesn't quite 
work for non-strictly outer-\\ minimizing
horizons, but equation \ref{eqn:App} does prove that each component of the 
horizon is either a sphere or a torus.  The torus possibility is then 
ruled out by \cite{GG}.  

We will
need theorem \ref{thm:stability} in sections 
\ref{sec:bounded} and \ref{sec:asymptotic}.

\section{The Exponential Growth Rate of $\mbox{Diam}(\Sigma(t))$}
\label{sec:exponential}

In this section we will show that the 
diameter of $\Sigma(t)$ is growing approximately exponentially
with respect to the fixed metric $g_0$.  This fact will be used 
in section \ref{sec:bounded} to
prove that $\Sigma(t)$ eventually encloses any bounded set.  It will
also be used in section \ref{sec:asymptotic} to help bound the behavior 
of $\Sigma(t)$ for large $t$.

We recall from section \ref{sec:definitions} that we may write
$g_0 = {\cal U}_0(x)^4 g_{flat}$, where the harmonically flat end of 
$(M^3,g_{flat})$ is isometric to $(\real^3 \backslash  B_{r_0}(0), \delta)$,
where now we have chosen $r_0$ such that we can require 
${\cal U}_0(x)$ to go to one at infinity.    
Let $S(r)$ denote the sphere of radius
$r$ in $\real^3$ centered at zero.  Then for $r \ge r_0$, we can 
define $S(r)$ to be a sphere in $(M^3, g_{flat})$.  

\begin{theorem}\label{thm:diam}
Given any $\tilde{t} \ge 0$, there exists a $t \ge \tilde{t}$ such that 
$\Sigma(t)$ is not entirely enclosed by $S(r(t))$ (and
$r(t) \ge r_0$), where 
\begin{equation}
r(t) = \left(\frac{A_0}{65\pi}\right)^{1/2} e^{2t} ,
\end{equation}
and $A_0 = A(0)$.
\end{theorem}
{\it Proof.}
The proof is a proof by contradiction.  We will assume that for some 
$\tilde{t} \ge 0$ (with $r(\tilde{t}) \ge r_0$)
that $\Sigma(t)$ is entirely enclosed by $S(r(t))$ for all $t \ge \tilde{t}$.
Then we will use this to prove that $A(t)$ is not constant, contradicting 
theorem \ref{thm:monotone}.

In fact, we will show that $|S(r(t))|_{g_t}$ (which is the area of $S(r(t))$
with respect to the metric $g_t$) is less than $A_0 \equiv A(0)$ for some 
$t \ge \tilde{t}$.  Then since $A(t)$ is defined to be the area of $\Sigma(t)$
which is defined to be the outermost minimal area enclosure of $\Sigma_0$, 
this will force $A(t) < A_0$, a contradiction.

As in section \ref{sec:m(t)}, we define 
$g_t = u_t(x)^4 g_0$ and $g_t = {\cal U}_t(x)^4 g_{flat}$ so that
${\cal U}_t(x) = u_t(x) {\cal U}_0(x)$.  In addition, we define
${\cal V}_t(x) = v_t(x){\cal U}_0(x)$ so that
\begin{equation}
   \frac{d}{dt} {\cal U}_t(x) = {\cal V}_t(x)
\end{equation}
since $\frac{d}{dt} u_t(x) = v_t(x)$.
Next we define ${\cal Q}_t(x)$ outside $S(r(t))$ in $(M^3,g_{flat})$ such that
\begin{equation}\label{eqn:Q}
\left\{ \begin{array}{rll}
Q_t = & 0 & \mbox{on } S(r(t)) \\
\Delta_{g_{flat}} Q_t \equiv & 0 & \mbox{outside } S(r(t)) \\
Q_t \rightarrow & -e^{-t} & \mbox{at infinity.}
\end{array}\right.
\end{equation}
so that 
\begin{equation}\label{eqn:Qexp}
{\cal Q}_t(x)= -e^{-t}\left( 1 - \frac{r(t)}{|x|} \right)
\end{equation}
outside $S(r(t))$.

Then since we are assuming that $\Sigma(t)$ is entirely enclosed by $S(r(t))$, 
it follows from equations \ref{eqn:ODE3} and \ref{eqn:identity} that
\begin{equation}\label{eqn:VV}
\left\{ \begin{array}{rll}
{\cal V}_t \le & 0 & \mbox{on } S(r(t)) \\
\Delta_{g_{flat}} {\cal V}_t \equiv & 0 & \mbox{outside } S(r(t)) \\
{\cal V}_t \rightarrow & -e^{-t} & \mbox{at infinity.}
\end{array}\right.
\end{equation}
Hence, by the maximum principle,
\begin{equation}
  {\cal V}_t(x) \le {\cal Q}_t(x)
\end{equation}
for all $x$ outside $S(r(t))$.  Consequently,
\begin{eqnarray} 
{\cal U}_t(x) &=& {\cal U}_{\tilde{t}}(x) + 
                  \int_{\tilde{t}}^t {\cal V}_s(x) \; ds \\
            &\le& {\cal U}_{\tilde{t}}(x) + 
                  \int_{\tilde{t}}^t {\cal Q}_s(x) \; ds 
\end{eqnarray}
for $x$ outside $S(r(t))$.
Now choose $k$ such that 
${\cal U}_{\tilde{t}}(x) \le e^{-\tilde{t}} + \frac{k}{|x|}$
outside $S(r(\tilde{t}))$.
Then using this and equation \ref{eqn:Qexp}, we get
\begin{equation}\label{eqn:calUbound}
   {\cal U}_t(x) \le e^{-t} + \frac{1}{|x|} \left[ k +
   \sqrt{\frac{A_0}{65\pi}} \left(e^t - e^{\tilde{t}}\right) \right]
\end{equation}
for $x$ outside $S(r(t))$.

Now we define ${\cal A}(t)$ to be the area of $S(r(t))$ in $(M^3,g_t)$.
Then 
\begin{eqnarray}
{\cal A}(t) &=& \int_{S(r(t))} {\cal U}_t(x)^4 \; dA_{g_{flat}} \\
  & \le & 4\pi r(t)^2 \left(\sup_{S(r(t))} {\cal U}_t(x)\right)^4 \\
  & \le & \frac{4}{65} A_0 \left[2 + e^{-t} \left(
 k\sqrt{\frac{65\pi}{A_0}} - e^{\tilde{t}}\right)\right]^4
\end{eqnarray} 
by inequality \ref{eqn:calUbound}, where $dA_{g_{flat}}$ is the 
area form of $S(r(t))$ in $(M^3,g_{flat})$.  Hence, 
\begin{equation}
   \lim_{t \rightarrow \infty} {\cal A}(t) = \frac{64}{65} A_0,
\end{equation}
which means that ${\cal A}(t) < A_0$ for some $t \ge \tilde{t}$.
But by theorem \ref{thm:monotone}, 
$|\Sigma(t)|_{g_t} = A_0$ for all $t \ge 0$,
and since $\Sigma(t)$ was defined to be a minimal area enclosure of the
original horizon $\Sigma_0$ in $(M^3,g_t)$, 
we have a contradiction.  \qed

\section{Proof That $\Sigma(t)$ Eventually Encloses Any \\ Bounded Set} 
\label{sec:bounded}

In this section we will prove that $\Sigma(t)$ eventually encloses any
bounded set in a finite amount of time.  In particular, this will allow 
us to conclude that $\Sigma(t)$ flows into the harmonically flat end 
of $(M^3,g_0)$ in a finite amount of time, which greatly simplifies
the discussion in the next two sections.

\begin{theorem}\label{thm:bounded}
Given any bounded set $B \subset M^3$, there exists
a $t \ge 0$ such that $\Sigma(t)$ encloses $B$.
\end{theorem}
{\it Proof.}
The proof of this theorem is a proof by contradiction.  We will show that 
if the theorem were false, then we would be able to find a $\bar{t}$ such 
that $A(\bar{t})$, the area of $\Sigma(\bar{t})$ in $(M^3,g_{\bar{t}})$, 
was greater than $A_0$, contradicting theorem \ref{thm:monotone}.

So suppose
there exists a bounded set $B \subset M^3$ such that $\Sigma(t)$ does not 
(entirely) 
enclose $B$ for any $t \ge 0$.  Then since $B$ is bounded, there exists an 
$R_1 > r_0$ such that the coordinate sphere $S(R_1)$ (defined in the previous
section) encloses $B$.  Hence, $\Sigma(t)$ does not enclose $S(R_1)$ for any
$t$.  

Now choose any $R_2 > R_1$.  By theorem \ref{thm:diam}, there
exists a $\bar{t}$ such that $\Sigma(\bar{t})$ is not enclosed by $S(R_2)$.
Then it follows that at least one component $\Sigma_i(\bar{t})$ of 
$\Sigma(\bar{t})$ must satisfy 
\begin{equation}\label{eqn:SigmaInt}
\Sigma_i(\bar{t}) \cap S(R) \ne \varnothing
\end{equation}
for all 
$R_1 \le R \le R_2$.  This follows from the fact that each component of the 
region inside $\Sigma(\bar{t})$ must intersect the region inside $S(r_0)$ 
since otherwise the area of $\Sigma(\bar{t})$ could be decreased by simply 
eliminating that component.

We will show that it is possible to choose $R_2$ large enough to make the area
of $\Sigma_i(\bar{t})$ in $(M^3,g_{\bar{t}})$ as large as we want, which is a 
contradiction.  Since 
\begin{equation}\label{eqn:Uarea}
   |\Sigma_i(\bar{t})|_{g_{\bar{t}}} = \int_{\Sigma_i(\bar{t})} 
   {\cal U}_{\bar{t}}(x)^4 \; dA_{g_{flat}}
\end{equation}
where $dA_{g_{flat}}$ is the area form of $\Sigma_i(\bar{t})$ in 
$(M^3,g_{flat})$, we can show that $|\Sigma_i(\bar{t})|_{g_{\bar{t}}}$ is large
if we have adequate lower bounds on ${\cal U}_{\bar{t}}(x)$ and on the area
of $|\Sigma_i(\bar{t})|$ in $(M^3,g_{flat})$.

To find a lower bound 
on ${\cal U}_{\bar{t}}(x)$ we will use theorem \ref{thm:shf}
in appendix \ref{sec:shf}.  Let
\begin{equation}
   {\cal U}_\infty(x) = \lim_{t \rightarrow \infty} {\cal U}_t(x)
\end{equation}
Then since $u_t(x)$ is decreasing in $t$ by equation \ref{eqn:ODE4}, so is
${\cal U}_t(x)$.  Hence, 
\begin{equation}\label{eqn:Uge}
{\cal U}_{t}(x) \ge {\cal U}_\infty(x).
\end{equation}
for all $t \ge 0$.
However, to use theorem \ref{thm:shf}, we need a superharmonic function defined
on $\real^3$, whereas ${\cal U}_\infty(x)$ is defined on $(M^3,g_{flat})$ 
which is only isometric to $(\real^3 \backslash B_{r_0}(0),\delta)$ 
in the harmonically flat end.  Let $c>0$ be the minimum value of 
${\cal U}_\infty(x)$ on $S(r_0)$. 
and let $\phi$ be the isometry mapping 
$(\real^3 \backslash B_{r_0}(0),\delta)$ to the harmonically flat end of 
$(M^3,g_{flat})$.  Then we can define ${\cal U}(x)$ on $\real^3$ as
\begin{equation}
   {\cal U}(x) = \left\{\begin{array} {cl}
   \min(c/2, \; {\cal U}_\infty(\phi(x))) , & \mbox{ for } |x| \ge r_0 \\
   c/2, &  \mbox{ for } |x| < r_0. \end{array}\right.
\end{equation}

Since the minimum value of two superharmonic functions is superharmonic, it
follows that ${\cal U}(x)$ is superharmonic on 
$(\real^3,\delta)$.
Furthermore, since 
${\cal U}_\infty(x)$ 
goes to zero at infinity, there exists an $\tilde{r} \ge r_0$ 
such that
\begin{equation}\label{eqn:Uequal}
   {\cal U}(x) = {\cal U}_\infty(\phi(x))
\end{equation}
for $|x| > \tilde{r}$.  
(Without loss of generality, we will assume that we chose $R_1$ from before so
that $R_1 > \tilde{r}$.)  

Next, we observe that 
\begin{equation}
   |S(R)|_{g_t} = \int_{S(R)} 
   {\cal U}_t(x)^4 \; dA_{g_{flat}} \ge A_0
\end{equation}
by theorem \ref{thm:monotone} since $\Sigma(t)$ has area $A_0$ and 
is a minimal area enclosure of $\Sigma_0$ in $(M^3,g_t)$. 
Hence, taking the limit as $t$ goes to infinity, it follows from equation
\ref{eqn:Uequal} that  
\begin{equation}
   \int_{S_R(0)} {\cal U}(x)^4 \; dA \ge A_0
\end{equation}
in $(\real^3,\delta)$ for $R > \tilde{r}$, so that by theorem \ref{thm:shf} 
in appendix \ref{sec:shf}
\begin{equation}
   {\cal U}(x) \ge c \, A_0^{1/4} |x|^{-1/2}
\end{equation}
for $|x| \ge \tilde{r}$ and for some $c > 0$.  Thus, by equations 
\ref{eqn:Uequal} and \ref{eqn:Uge}
\begin{equation}\label{eqn:Ulb}
   {\cal U}_{t}(x) \ge c \, A_0^{1/4} |x|^{-1/2}
\end{equation}
for all $t \ge 0$ and for all $x$ in $M^3$ 
outside $S(\tilde{r})$.  This lower bound on 
${\cal U}_{t}(x)$ is the first of two steps needed to prove 
that integral in equation \ref{eqn:Uarea} is large.

The second step of the proof is to find the right lower bound on the area 
of $\Sigma_i(\bar{t})$ in $(M^3,g_{flat})$.  So far all we have is 
equation \ref{eqn:SigmaInt} which tells us that the diameter of 
$\Sigma_i(\bar{t})$ in $(M^3,g_{flat})$ is at least $R_2 - R_1$.  Naturally
this is not enough to bound the area from below without some control
on the possible geometries of $\Sigma_i(\bar{t})$ in $(M^3,g_{flat})$.

Fortunately, we do have very good control on the geometry 
of each component $\Sigma_i(\bar{t})$ of
$\Sigma(\bar{t})$ by virtue of theorem \ref{thm:stability}.  
Using the same notation as in the Gauss equation given
in equation \ref{eqn:Gauss}, we can substitute
\begin{equation}
   (\lambda_1 - \lambda_2)^2 = (2R - 4Ric(\nu,\nu)) - 4K + H^2
\end{equation}
into equation \ref{eqn:stability} to get
\begin{equation}\label{eqn:H^2_Bound}
   \int_{\Sigma_i(\bar{t})} H^2 \; dA_{g_{flat}} \le  32\pi + 
   \int_{\Sigma_i(\bar{t})} (4Ric(\nu,\nu) - 2R) \; dA_{g_{flat}}
\end{equation}
where we have used the Gauss-Bonnet theorem and the fact that every component
of $\Sigma(\bar{t})$ is a sphere.  In addition, $4Ric(\nu,\nu) - 2R$ is 
zero in the harmonically flat end of $M^3$ since $g_{flat}$ is flat.

Now we want to show that the right hand side of equation \ref{eqn:H^2_Bound}
is bounded.  
Let $K$ be the compact set of points outside (or on) 
the original horizon $\Sigma_0$ and inside (or on) $S(r_0)$.  Let
\begin{equation}
   R_{max} = \sup_K |4Ric(\nu,\nu) - 2R|
\end{equation}
which is finite since $(M^3,g_{flat})$ is smooth and $K$ is compact.
Then since
\begin{equation}\label{eqn:Uarea2}
   A_0 = |\Sigma(\bar{t})|_{g_{\bar{t}}} 
   \ge |\Sigma_i(\bar{t}) \cap K|_{g_{\bar{t}}}
   = \int_{\Sigma_i(\bar{t}) \cap K} 
   {\cal U}_{\bar{t}}(x)^4 \; dA_{g_{flat}}
\end{equation}
it follows that
\begin{equation}
   |\Sigma_i(\bar{t}) \cap K|_{g_{flat}} \le \frac{A_0}
   {\inf_K {\cal U}_{\bar{t}}(x)^4} \le \frac{A_0}
   {\inf_K u(x)^4 {\cal U}_0(x)^4}
\end{equation}
where we recall that ${\cal U}_{\bar{t}}(x) = u_{\bar{t}}(x) {\cal U}_0(x)$
and 
where $u(x)$ is defined to be the positive 
harmonic function which equals $1$ (and
hence also $u_{\bar{t}}(x)$) on the original horizon $\Sigma_0$
and goes to zero at infinity
and hence is a barrier function for the superharmonic function 
$u_{\bar{t}}(x)$ in $(M^3,g_{\bar{t}})$.  Thus, we have that 
\begin{equation}\label{eqn:Wbound}
   \int_{\Sigma_i(\bar{t})} H^2 \; dA_{g_{flat}} \le  32\pi + 
   \frac{R_{max} A_0}{\inf_K u(x)^4 {\cal U}_0(x)^4} \equiv k.
\end{equation}

We will need the Willmore functional bound in inequality \ref{eqn:Wbound}
to use the identity (equation 16.31 in \cite{GT})
\begin{equation}
   |\Sigma \cap B_r(x)| \ge \pi r^2 \left(1 - \frac{1}{16\pi}
   \int_{\Sigma \cap B_r(x)} H^2 \; dA  \right)
\end{equation}
where $\Sigma$ is any smooth, compact 
surface which is the boundary of a region 
in $\real^3$ and everything is with 
respect to the standard flat metric $\delta$.  Since $(M^3, g_{flat})$ is 
flat in the harmonically flat end region of $M^3$, we will be able to use this 
identity on $\Sigma_i(\bar{t})$.

Since our choice of $R_2$ could be arbitrarily large, let
$R_2 = 3(2^{n^2} - 1) R_1$ where $n$ is a positive integer which may be 
chosen to be arbitrarily large.  Then by equation \ref{eqn:SigmaInt} 
we can choose 
\begin{equation}
x_k \in \Sigma_i(\bar{t}) \cap S(3(2^k - 1)R_1)
\end{equation} 
for $1 \le k \le n^2$.  Then if we define $r_k = 2^k R_1$
it follows that the balls $B_{r_k}(x_k)$ are all disjoint.
Hence,
\begin{equation}
   \int_{\Sigma_i(\bar{t}) \cap B_{r_k}(x_k)} H^2 \; dA_{g_{flat}} \le 
   \frac{k}{n}
\end{equation}
except at most $n$ different values of $k$.  Then for these values of $k$,
it follows from equation \ref{eqn:Ulb} that 
\begin{eqnarray}\label{eqn:lb_for_int}
   \int_{\Sigma_i(\bar{t}) \cap B_{r_k}(x_k)} {\cal U}_{\bar{t}}(x)^4 \; 
   dA_{g_{flat}}  & \ge & c^4 A_0 \int_{\Sigma_i(\bar{t}) \cap B_{r_k}(x_k)} 
   |x|^{-2} \; dA_{g_{flat}} \\
   & \ge & c^4 A_0 \; |{\Sigma_i(\bar{t}) \cap B_{r_k}(x_k)}| \; 
   (2^{k+2}R_1)^{-2} \\
   & \ge &  c^4 A_0 \; \pi r_k^2 \left(1 - \frac{k}{16\pi n}\right)  
   (2^{k+2}R_1)^{-2} \\
   & = & c^4 A_0 \frac{\pi}{16} \left(1 - \frac{k}{16\pi n}\right) 
\end{eqnarray}
where $|x|$ is defined to be $r$ on $S(r)$ in the harmonically flat end 
of $(M^3,g_{flat})$.  Hence, we have that
\begin{eqnarray}
   A_0 & = & |\Sigma(\bar{t})|_{g_{\bar{t}}} \;\;\ge\;\;  
   |\Sigma_i(\bar{t})|_{g_{\bar{t}}} \\
   & = & \int_{\Sigma_i(\bar{t})} {\cal U}_{\bar{t}}(x)^4 \; dA_{g_{flat}} \\
   & \ge & \sum_{k = 1}^{n^2} \int_{\Sigma_i(\bar{t}) \cap B_{r_k}(x_k)} 
   {\cal U}_{\bar{t}}(x)^4 \; dA_{g_{flat}} \\
   & \ge & c^4 A_0 \frac{\pi}{16} \left(1 - \frac{k}{16\pi n}\right) (n^2 - n)
\end{eqnarray}
which is a contradiction since $n$ can be chosen to be arbitrarily large.
Hence, given any bounded set $B \subset M^3$, there must exist a $t \ge 0$
such that $\Sigma(t)$ encloses $B$. \qed

We immediately deduce a very useful corollary.  Since $\Sigma(t)$ always
flows outwards and 
must eventually entirely enclose $S(r_0)$ by the above
theorem, it follows that
after a certain point in time $\Sigma(t)$ is entirely in the harmonically
flat end of $(M^3,g_{flat})$.
Furthermore, since
$\Sigma(t)$ is defined to be the outermost surface with 
minimum area in $(M^3,g_t)$ which encloses the original
horizon,
these $\Sigma(t)$ only have one component
since having any additional components would only increase
the area of $\Sigma(t)$.  
Then since it follows from a stability argument that each component 
is a sphere (\cite{SY1} or see
the end of section \ref{sec:stability}), $\Sigma(t)$ is a single sphere.
Thus, we get the following
corollary to theorem \ref{thm:bounded}.

\begin{corollary}\label{cor:bounded}
There exists a $t_0 \ge 0$ such that for all $t \ge t_0$,
topologically $\Sigma(t)$ is a single sphere
and is in the 
harmonically flat end of $(M^3,g_{flat})$ which 
is isometric to $(\real^3 \backslash  B_{r_0}(0), \delta)$.
\end{corollary}

We remind the reader that $(M^3,g_{flat})$ was defined in section 
\ref{sec:definitions} and is conformal to the harmonically flat 
manifold $(M^3,g_0)$.

\section{Bounds on the Behavior of $\Sigma(t)$}
\label{sec:asymptotic}

The main objective of this section is to achieve upper and lower bounds 
for the diameter of $\Sigma(t)$ in $(M^3,g_0)$ for large $t$.  However,
in the process of deriving these bounds we also prove other interesting
although not essential bounds on the behavior of $\Sigma(t)$.

As in the previous two sections,
$g_0 = {\cal U}_0(x)^4 g_{flat}$, where the harmonically flat end of 
$(M^3,g_{flat})$ is isometric to $(\real^3 \backslash  B_{r_0}(0), \delta)$,
and ${\cal U}_0(x)$ goes to one at infinity.    
From corollary \ref{cor:bounded} of the previous section, we see that
$\Sigma(t)$ is in the flat region of $(M^3,g_{flat})$ for $t \ge t_0$ so that
the behavior of $\Sigma(t)$ for large $t$ reduces to a problem in 
$\real^3$.  

Furthermore, from section \ref{sec:exponential} we know that
the diameter of $\Sigma(t)$ is going up roughly as
some constant times $e^{2t}$.  Then since
$\real^3$ is linear, it makes since and is convenient 
to rescale distances by the factor $e^{-2t}$ so that in this new rescaled
$\real^3$ the diameter of $\Sigma(t)$ is approximately bounded. 
In fact,
the goal of this section is to show that the diameter of  
$\Sigma(t)$ in this rescaled $\real^3$ is bounded above and below by constants.
(From this point on $\Sigma(t)$ will 
refer to the rescaled $\Sigma(t)$.)

We will use capital letters to denote rescaled quantities.
Recall that $g_t = {\cal U}_t(x)
^4 g_{flat}$ and $g_t = u_t(x)^4 g_0$
so that ${\cal U}_t(x) = u_t(x) \, {\cal U}_0(x)$.  Let
\begin{equation}\label{eqn:defU}
U_t(x) = e^t \; {\cal U}_t(x e^{2t}) .
\end{equation}
so that $U_t(x)$ goes to one at infinity for all $t$ since 
we arranged ${\cal U}_0(x)$ to go to one at infinity and $u_t(x)$ goes
to $e^{-t}$ at infinity.  Analogously, we define
${\cal V}_t(x) = v_t(x) {\cal U}_0(x)$ (so that 
$\frac{d}{dt} {\cal U}_t(x) =  {\cal V}_t(x))$ and we define
\begin{equation}
   V_t(x) = e^t \; {\cal V}_t(x e^{2t}). 
\end{equation}
Then we observe that $V_t(x)$ goes to $-1$ at infinity (since $v_t(x)$
goes to $-e^{-t}$ at infinity) and $V_t(x)$ equals zero on $\Sigma(t)$.
Furthermore, differentiating equation \ref{eqn:defU} gives us
\begin{equation}\label{eqn:U_t(x)}
   \frac{d}{dt} U_t(x) = V_t(x) + U_t(x) + 2r\frac{\partial}{\partial r} U_t(x)
\end{equation}
where $r$ is the radial coordinate in $\real^3$.  

\begin{lemma}\label{lem:m}  
For $t \ge t_0$ (as defined in corollary \ref{cor:bounded}),
the Riemannian manifold 
$(\real^3 \backslash B_{r_0 e^{-2t}}(0), U_t(x)^4 \delta)$
is isometric to the harmonically flat end of $(M^3,g_t)$ and has total
mass $m(t)$.
\end{lemma}

Consequently, $\Sigma(t)$ has zero mean curvature in  
$(\real^3 \backslash B_{r_0 e^{-2t}}(0), U_t(x)^4 \delta)$, 
from which it follows that 
$U_t H + 4 \frac{dU_t}{d\vec{\nu}} = 0$, where $H$ is the mean
curvature of $\Sigma(t)$  
and $\vec{\nu}$ is the 
outward pointing unit normal vector of $\Sigma(t)$ in 
$(\real^3,\delta)$.
Also, since $(M^3,g_t)$ has 
zero scalar curvature outside $\Sigma(t)$ (for $t \ge t_0$
in corollary \ref{cor:bounded}), 
it follows from equation \ref{eqn:scalar_curv} that $U_t(x)$ 
is harmonic in $(\real^3, \delta)$
outside $\Sigma(t)$.  Then from the discussion and definitions
in the above paragraphs this in turn implies that ${\cal U}_t(x)$ 
is harmonic, which implies that ${\cal V}_t(x)$ is harmonic, 
which implies that $V_t(x)$ is also harmonic outside $\Sigma(t)$.
To summarize:
\vspace{.2in}
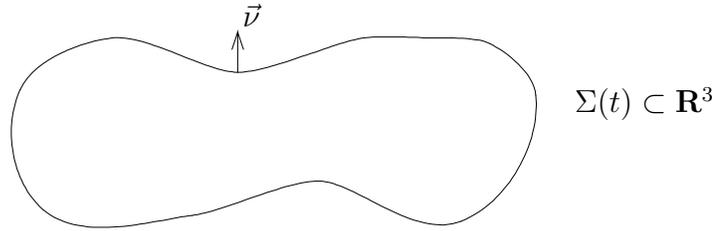
\begin{center}
\setlength{\unitlength}{0.0060in}
\begin{picture}(514,199)(0,-10)
\path(199,138)(199,173)(194,158)
\path(199,173)(204,158)
\path(309,168)	(306.451,167.604)
	(303.741,167.083)
	(300.881,166.445)
	(297.880,165.700)
	(294.750,164.857)
	(291.500,163.924)
	(288.141,162.911)
	(284.683,161.827)
	(281.137,160.681)
	(277.513,159.482)
	(273.821,158.239)
	(270.072,156.960)
	(266.275,155.656)
	(262.443,154.336)
	(258.584,153.007)
	(254.709,151.680)
	(250.828,150.363)
	(246.953,149.065)
	(243.092,147.796)
	(239.257,146.564)
	(235.458,145.379)
	(231.705,144.249)
	(228.009,143.184)
	(224.380,142.192)
	(220.828,141.284)
	(217.364,140.467)
	(213.998,139.751)
	(210.740,139.145)
	(207.601,138.657)
	(204.591,138.298)
	(201.721,138.076)
	(199.000,138.000)

\path(199,138)	(196.291,138.097)
	(193.454,138.378)
	(190.498,138.832)
	(187.430,139.443)
	(184.258,140.201)
	(180.991,141.091)
	(177.636,142.101)
	(174.202,143.217)
	(170.697,144.426)
	(167.129,145.716)
	(163.505,147.074)
	(159.835,148.485)
	(156.125,149.938)
	(152.385,151.418)
	(148.622,152.914)
	(144.844,154.412)
	(141.060,155.899)
	(137.277,157.361)
	(133.504,158.787)
	(129.748,160.162)
	(126.018,161.474)
	(122.322,162.710)
	(118.668,163.857)
	(115.064,164.901)
	(111.518,165.829)
	(108.038,166.630)
	(104.632,167.288)
	(101.309,167.792)
	(98.076,168.128)
	(94.941,168.284)
	(91.913,168.245)
	(89.000,168.000)

\path(89,168)	(84.205,167.298)
	(81.669,166.846)
	(79.054,166.327)
	(76.367,165.740)
	(73.617,165.084)
	(70.814,164.360)
	(67.966,163.567)
	(65.081,162.705)
	(62.170,161.773)
	(59.239,160.771)
	(56.299,159.699)
	(53.358,158.557)
	(50.425,157.343)
	(47.509,156.059)
	(44.618,154.703)
	(41.761,153.274)
	(38.947,151.774)
	(36.185,150.201)
	(33.484,148.555)
	(30.852,146.836)
	(28.298,145.044)
	(23.461,141.237)
	(19.042,137.131)
	(15.112,132.725)
	(11.741,128.016)
	(10.287,125.546)
	(9.000,123.000)

\path(9,123)	(7.665,119.943)
	(6.456,116.796)
	(5.375,113.567)
	(4.420,110.263)
	(3.593,106.893)
	(2.892,103.463)
	(2.318,99.981)
	(1.871,96.454)
	(1.550,92.891)
	(1.356,89.298)
	(1.289,85.683)
	(1.348,82.054)
	(1.534,78.417)
	(1.845,74.782)
	(2.284,71.154)
	(2.848,67.542)
	(3.539,63.953)
	(4.355,60.395)
	(5.298,56.875)
	(6.367,53.400)
	(7.562,49.979)
	(8.882,46.618)
	(10.329,43.325)
	(11.901,40.109)
	(13.599,36.975)
	(15.423,33.932)
	(17.372,30.987)
	(19.447,28.148)
	(21.647,25.423)
	(23.973,22.818)
	(26.424,20.341)
	(29.000,18.000)

\path(29,18)	(32.301,15.382)
	(35.773,13.048)
	(39.404,10.985)
	(43.181,9.180)
	(47.092,7.621)
	(51.126,6.296)
	(55.269,5.191)
	(59.510,4.296)
	(63.836,3.596)
	(68.236,3.081)
	(72.698,2.736)
	(77.208,2.551)
	(81.755,2.511)
	(86.327,2.606)
	(90.912,2.822)
	(95.497,3.147)
	(100.071,3.569)
	(104.620,4.074)
	(109.134,4.651)
	(113.599,5.288)
	(118.005,5.971)
	(122.337,6.688)
	(126.585,7.427)
	(130.736,8.175)
	(134.778,8.920)
	(138.699,9.649)
	(142.487,10.351)
	(146.129,11.012)
	(149.614,11.619)
	(152.928,12.162)
	(156.061,12.626)
	(159.000,13.000)

\path(159,13)	(161.569,13.353)
	(164.297,13.838)
	(167.172,14.445)
	(170.185,15.164)
	(173.325,15.987)
	(176.582,16.904)
	(179.947,17.905)
	(183.408,18.981)
	(186.956,20.123)
	(190.580,21.322)
	(194.270,22.568)
	(198.016,23.851)
	(201.809,25.163)
	(205.636,26.493)
	(209.490,27.833)
	(213.358,29.174)
	(217.231,30.505)
	(221.100,31.818)
	(224.953,33.102)
	(228.780,34.350)
	(232.571,35.551)
	(236.317,36.695)
	(240.006,37.775)
	(243.629,38.779)
	(247.176,39.700)
	(250.635,40.527)
	(253.998,41.251)
	(257.254,41.863)
	(260.392,42.353)
	(263.402,42.712)
	(266.275,42.931)
	(269.000,43.000)

\path(269,43)	(272.071,42.843)
	(275.250,42.407)
	(278.533,41.715)
	(281.914,40.787)
	(285.388,39.646)
	(288.948,38.313)
	(292.590,36.811)
	(296.308,35.161)
	(300.096,33.384)
	(303.950,31.504)
	(307.863,29.542)
	(311.829,27.519)
	(315.844,25.457)
	(319.902,23.379)
	(323.998,21.306)
	(328.125,19.259)
	(332.279,17.262)
	(336.453,15.335)
	(340.643,13.500)
	(344.843,11.780)
	(349.047,10.196)
	(353.249,8.769)
	(357.446,7.523)
	(361.629,6.478)
	(365.796,5.657)
	(369.939,5.081)
	(374.053,4.772)
	(378.133,4.752)
	(382.173,5.043)
	(386.168,5.667)
	(390.112,6.645)
	(394.000,8.000)

\path(394,8)	(397.762,9.627)
	(401.450,11.451)
	(405.061,13.465)
	(408.591,15.660)
	(412.035,18.028)
	(415.391,20.561)
	(418.655,23.251)
	(421.821,26.090)
	(424.888,29.070)
	(427.851,32.183)
	(430.705,35.420)
	(433.449,38.775)
	(436.076,42.238)
	(438.585,45.802)
	(440.970,49.459)
	(443.228,53.200)
	(445.356,57.018)
	(447.349,60.904)
	(449.204,64.851)
	(450.917,68.850)
	(452.484,72.894)
	(453.901,76.974)
	(455.165,81.082)
	(456.271,85.211)
	(457.216,89.351)
	(457.996,93.496)
	(458.608,97.637)
	(459.047,101.765)
	(459.309,105.874)
	(459.391,109.955)
	(459.290,114.000)
	(459.000,118.000)

\path(459,118)	(458.389,121.719)
	(457.313,125.375)
	(455.817,128.954)
	(453.948,132.445)
	(451.750,135.836)
	(449.269,139.115)
	(446.552,142.268)
	(443.644,145.284)
	(440.591,148.150)
	(437.438,150.855)
	(434.232,153.386)
	(431.018,155.730)
	(427.842,157.876)
	(424.750,159.811)
	(421.788,161.523)
	(419.000,163.000)

\path(419,163)	(416.356,164.142)
	(413.383,165.095)
	(410.127,165.874)
	(406.637,166.497)
	(402.961,166.982)
	(399.148,167.344)
	(395.245,167.602)
	(391.299,167.771)
	(387.361,167.870)
	(383.476,167.914)
	(379.694,167.922)
	(376.062,167.910)
	(372.629,167.894)
	(369.442,167.893)
	(366.550,167.922)
	(364.000,168.000)

\path(364,168)	(361.454,168.112)
	(358.554,168.233)
	(355.352,168.357)
	(351.896,168.479)
	(348.238,168.594)
	(344.428,168.698)
	(340.515,168.785)
	(336.551,168.851)
	(332.585,168.889)
	(328.668,168.897)
	(324.849,168.867)
	(321.180,168.797)
	(317.710,168.680)
	(314.490,168.511)
	(311.570,168.286)
	(309.000,168.000)

\put(194,178){\makebox(0,0)[lb]{\smash{{{\SetFigFont{12}{14.4}{rm}
$\vec{\nu}$ }}}}}
\put(484,103){\makebox(0,0)[lb]{\smash{{{\SetFigFont{12}{14.4}{rm}
$\Sigma(t) \subset \real^3$ }}}}}
\end{picture}
\end{center}
\begin{equation}\label{eqn:U}
\left\{ \begin{array}{rll}
U_t H + 4\frac{dU_t}{d\vec{\nu}} = & 0 & \mbox{on } \Sigma(t) \\
\Delta U_t \equiv & 0 & \mbox{outside } \Sigma(t) \\
U_t \rightarrow & 1 & \mbox{at infinity}
\end{array}\right.
\end{equation}
\begin{equation}\label{eqn:V}
\left\{ \begin{array}{rll}
V_t = & 0 & \mbox{on } \Sigma(t) \\
\Delta V_t \equiv & 0 & \mbox{outside } \Sigma(t) \\
V_t \rightarrow & -1 & \mbox{at infinity.}
\end{array}\right.
\end{equation}
\vspace{.2in}

\noindent
Equations \ref{eqn:U_t(x)}, \ref{eqn:U}, and \ref{eqn:V} characterize the
new rescaled first order o.d.e.~in $t$ for $U_t(x)$.   
On the one hand $U_t(x)$ determines
$\Sigma(t)$ since $\Sigma(t)$ is the outermost area minimizing horizon
of $(\real^3 \backslash B_{r_0 e^{-2t}}(0), U_t(x)^4 \delta)$, 
and on the other hand $\Sigma(t)$
determines $V_t(x)$ by equation \ref{eqn:V} which determines the first order
rate of change of $U_t(x)$ in equation \ref{eqn:U_t(x)}.  In the next 
section we will prove that $U_t(x)$ actually converges to $1 + \frac{M}{2r}$
for some positive $M$ in the limit as $t$ goes to infinity in this 
o.d.e.~and that $\Sigma(t)$ converges to a sphere of radius $m/2$.  
However, first
it is necessary to prove that the diameter of $\Sigma(t)$
is bounded, which is what we will do in this section.

In section \ref{sec:m(t)} we proved that $m(t)$ was nonincreasing.  In 
fact, by closely reexamining equations \ref{eqn:mtilde}, \ref{eqn:symenergy}, 
and \ref{eqn:mprime} we have that
\begin{equation}\label{eqn:mprime2}
  m'(t) = - 2 \tilde{m}(t) \le 0
\end{equation}
where $\tilde{m}(t)$ is the total mass of the manifold 
$(\bar{M}^3_{\Sigma(t)},\tilde{g}_t )$, where $\tilde{g}_t = 
\phi(x)^4 \bar{g_t}$, $M^3_{\Sigma(t)}$ is
the closed region of $M^3$ which is outside or on $\Sigma(t)$, 
$(\bar{M}^3_{\Sigma(t)}, \bar{g_t})$ is the manifold obtained by reflecting
$(M^3_{\Sigma(t)}, g_t)$ through $\Sigma(t)$, and $\phi(x)$ is the harmonic
function on $(\bar{M}^3_{\Sigma(t)}, \bar{g_t})$ which goes to one in the
original end and zero in the other end.  

Then since 
$(M^3_{\Sigma(t)}, g_t)$ is isometric to 
$(\real^3_{\Sigma(t)}, U_t(x)^4 \delta)$ (where $\real^3_{\Sigma(t)}$ is the
region in $\real^3$ outside $\Sigma(t)$),  
$(\bar{M}^3_{\Sigma(t)}, \tilde{g_t})$ is isometric to 
$(\real^3_{\Sigma(t)}, W_t(x)^4 \delta)$ where $W_t(x) = \phi(x)U_t(x)$.
Furthermore, by equation \ref{eqn:identity} 
it follows that $W_t(x)$ is harmonic
in $(\real^3,\delta)$,  
and since $\phi(x) = \frac12$ on $\Sigma(t)$ by symmetry and 
$U_t(x)$ and $\phi(x)$ both go to one at infinity, we have that
\vspace{.2in}
\begin{equation}\label{eqn:W}
\left\{ \begin{array}{rll}
W_t = & \frac12 U_t & \mbox{on } \Sigma(t) \\
\Delta W_t \equiv & 0 & \mbox{outside } \Sigma(t) \\
W_t \rightarrow & 1 & \mbox{at infinity.}
\end{array}\right.
\end{equation}
\vspace{.05in}

\begin{lemma}  \label{lem:mtilde}
For $t \ge t_0$,
the Riemannian manifold 
$(\real^3_{\Sigma(t)}, W_t(x)^4 \delta)$
is isometric to 
$(\bar{M}^3_{\Sigma(t)}, \tilde{g}_t)$ defined above
and has total
mass $\tilde{m}(t)$.
\end{lemma}

\begin{corollary}\label{cor:m,mtilde}
For $t \ge t_0$,
\begin{equation}  
m(t) = - \frac{1}{2\pi} \int_{\Sigma(t)} \frac{dU_t}{d\vec{\nu}} 
\end{equation}
and 
\begin{equation}
\tilde{m}(t) = - \frac{1}{2\pi} \int_{\Sigma(t)} \frac{dW_t}{d\vec{\nu}} .
\end{equation}
\end{corollary}
{\it Proof.}
In fact, by the divergence theorem, the above equations are true if we replace 
$\Sigma(t)$ with any homologous surface containing $\Sigma(t)$ since
$U_t(x)$ and $W_t(x)$ are harmonic in $\real^3$.  Then the corollary 
follows from the definition of total mass given in definition 
\ref{def:totalmass} where we consider the above statements with $\Sigma(t)$
replaced by a large sphere at infinity.  \qed

Furthermore, since $U_t(x)$, $V_t(x)$, and $W_t(x)$ are all
harmonic functions in $(\real^3,\delta)$ outside $\Sigma(t)$, 
it follows from their boundary values that
\begin{equation}\label{eqn:V=U-2W}
V_t(x) = U_t(x) - 2W_t(x)
\end{equation}
outside $\Sigma(t)$.
Then plugging this into equation \ref{eqn:U_t(x)} we get
\begin{equation}\label{eqn:U_t(x)2}
   \frac{d}{dt} U_t(x) = 2 ( U_t(x) - W_t(x) + 
                         r\frac{\partial}{\partial r} U_t(x) )
\end{equation}
outside $\Sigma(t)$.

The above equation reveals the key idea we will use to study the behavior
of the o.d.e.~for $U_t(x)$.  By equation \ref{eqn:mprime2}, it follows that
$\tilde{m}(t)$ must be going to zero for large $t$ since $m(t)$ cannot become
less than zero by the positive mass theorem.  Also, the positive mass 
theorem states that 
there is only one zero mass metric, namely $(\real^3,\delta)$, and
in section \ref{sec:limit} we will use this fact to prove that $W_t(x)$ is
approaching the constant function one.  Then it follows from studying 
equation \ref{eqn:U_t(x)2} that $U_t(x)$ approaches $1+\frac{M}{2r}$ for 
some positive $M$ in the limit as $t$ goes to infinity. With a few additional
observations this will prove 
that $(M^3,g_t)$ converges to a Schwarzschild metric 
outside the horizon $\Sigma(t)$ as claimed in theorem \ref{thm:limit}.

In the rest of this section, we will show that the rescaled horizon
$\Sigma(t)$ is very well behaved in $\real^3$ as $t \rightarrow \infty$.
In fact, we will show that both the areas and the diameters of the 
surfaces $\Sigma(t)$ have uniform upper and lower bounds.  
Later in section \ref{sec:limit} we will use this to
prove that the harmonic functions $U_t(x)$, $V_t(x)$, and 
$W_t(x)$ also have upper and lower bounds independent of $t$, which will
be needed when we take limits of these harmonic function.

We recall that by corollary \ref{cor:bounded}, $\Sigma(t)$ has 
only one component for $t \ge t_0$, and that this component is a sphere.  
Furthermore, by equation \ref{eqn:H^2_Bound},
it follows that the conformal-invariant (and hence scale-invariant) quantity
\begin{equation}\label{eqn:Willmore}
\int_{\Sigma(t)} H^2 d\mu \le 32\pi ,
\end{equation}
where $H$ is the mean curvature and $d\mu$ is the area form of $\Sigma(t)$
in $(\real^3,\delta)$.  

The Willmore functional of a surface in $\real^n$ is defined to
be one fourth of the integral of the mean curvature squared over the surface.
Surfaces with bounded Willmore functional have been
widely studied, and in particular, it was shown by L. Simon
in \cite{LS} that the ratio of the diameter squared to the 
area of a surface
is bounded both from above and from 
below by the Willmore functional of the surface.  
More precisely, for a surface $\Sigma$ in $\real^3$,    
\begin{equation}\label{eqn:adratio}
\frac{4}{\int_{\Sigma} H^2 d\mu} \le
\frac{\mbox{diam}(\Sigma)^2}{|\Sigma|} \le
\frac{C^2}{4} \int_{\Sigma} H^2 d\mu ,
\end{equation}
where $C$ is some positive constant.  
Hence, since theorem \ref{thm:diam} tells us that the
rescaled $\Sigma(t)$ satisfy
\begin{equation}\label{eqn:lower_diam}
\mbox{diam}(\Sigma(t)) \ge \left(\frac{A_0}{65\pi}\right)^{1/2} 
\end{equation}
for arbitrarily large values of $t$, we get the following corollary. 
\begin{corollary}\label{cor:lower_area}
Given any $\tilde{t} \ge 0$, there exists a $t \ge \tilde{t}$ such that
\begin{equation}  |\Sigma(t)| \ge \frac{A_0}{k} ,\end{equation} 
where $A_0 = A(0)$, 
$|\Sigma(t)|$ denotes the area of the rescaled 
$\Sigma(t)$ in $(\real^3,\delta)$, and $k = 8 \cdot 65 C^2 \pi^2$. 
\end{corollary}

Now going back to corollary \ref{cor:m,mtilde} and using equation 
\ref{eqn:V=U-2W}, we get that
\begin{equation}
m(t) - 2 \tilde{m}(t) = 
-\frac{1}{2\pi} \int_{\Sigma(t)} \frac{dV_t}{d\vec{\nu}} .
\end{equation}
Then since $\tilde{m}(t)$ is the total mass of 
$(\bar{M}_{\Sigma(t)}^3, \tilde{g}_t)$, by the positive mass theorem it
must be positive.  Hence,  
\begin{equation}
m(t) \ge 
- \frac{1}{2\pi} \int_{\Sigma(t)} \frac{dV_t}{d\vec{\nu}} .
\end{equation}
But from theorem \ref{thm:capacity_of_sigma} in appendix \ref{sec:lower} 
and equation \ref{eqn:Willmore}, we have that
\begin{equation}  - \frac{1}{2\pi} \int_{\Sigma(t)} \frac{dV_t}{d\vec{\nu}} 
\ge (24\pi)^{-1/2} |\Sigma(t)|^{1/2} \end{equation} 
which proves inequality \ref{eqn:lower_mass1} of the following theorem.
\begin{theorem} \label{thm:lower_mass}
For $t \ge t_0$ (as defined in corollary \ref{cor:bounded}),
\begin{equation}\label{eqn:lower_mass1}
m(t) \ge \left(\frac{|\Sigma(t)|}{24\pi}\right)^{1/2} 
\end{equation}
where $|\Sigma(t)|$ denotes the area of the rescaled 
$\Sigma(t)$ in $(\real^3,\delta)$.  
Also,
\begin{equation}\label{eqn:lower_mass2}
m(t) \ge \left( \frac{A_0}{24\pi k} \right)^{1/2} 
\end{equation}
for all $t \ge 0$, where again $A_0 = A(0)$
and $k = 8 \cdot 65 C^2 \pi^2$.  
\end{theorem}

Inequality \ref{eqn:lower_mass2} then follows from inequality 
\ref{eqn:lower_mass1}, corollary \ref{cor:lower_area}, and the fact that
$m(t)$ is nonincreasing.
We note that 
since $m(t)$ and $A_0$ are respectively the total mass of $(M^3,g_t)$ and 
the area of the horizon $\Sigma(t)$ in $(M^3,g_t)$, 
inequality \ref{eqn:lower_mass2} is a weak Penrose inequality
for $(M^3,g_t)$.  Furthermore, since the area of the horizon $A(t)$ is 
constant and $m(t)$ is non-increasing, we get this same weak Penrose inequality
for the original metric $(M^3,g_0)$.

Now going back to equation \ref{eqn:U}, we observe that 
\begin{equation}  m(t) = - \frac{1}{2\pi} \int_{\Sigma(t)} \frac{dU_t}{d\vec{\nu}} 
= \frac{1}{8\pi}\int_{\Sigma(t)} U_t(x) H \end{equation} 
Thus, by the Cauchy-Schwarz inequality,
\begin{equation}  m(t) \le \frac{1}{8\pi}\left(\int_{\Sigma(t)} U_t(x)^2\right)^{1/2}
\left(\int_{\Sigma(t)} H^2 \right)^{1/2} .\end{equation} 
Furthermore, since the area of $\Sigma(t)$ in $(M^3,g_t)$ 
equals $A_0$, by lemma \ref{lem:m}
\begin{equation}
 A_0 = \int_{\Sigma(t)} U_t(x)^4 .
\end{equation}
Thus, 
\begin{equation}  \int_{\Sigma(t)} U_t(x)^2 \le |\Sigma(t)|^{1/2} A_0^{1/2}, \end{equation} 
so that by equation \ref{eqn:Willmore} we have that 
\begin{equation}
m(t) \le (2\pi)^{-1/2} A_0^{1/4} |\Sigma(t)|^{1/4},
\end{equation}  
which, when combined with  
theorem \ref{thm:lower_mass}, gives inequality \ref{eqn:area_bounds} of the 
following theorem.
\begin{theorem}\label{thm:Areadiameter}
For $t \ge t_0$ (as defined in corollary \ref{cor:bounded}),
\begin{equation}\label{eqn:area_bounds}
\frac{1}{(12 k)^{2}} \le \frac{|\Sigma(t)|}{A_0} \le 12^2 
\end{equation}
and 
\begin{equation}\label{eqn:diam_bounds}
\frac{1}{8\pi (12 k)^2} \le \frac{\mbox{diam}(\Sigma(t))^2}{A_0}
\le 8\pi  (12 C)^2
\end{equation}
where $A_0 = A(0)$ and $k = 8 \cdot 65 C^2 \pi^2$. 
\end{theorem}
Inequality \ref{eqn:diam_bounds} then follows from inequalities 
\ref{eqn:area_bounds}, \ref{eqn:Willmore}, and \ref{eqn:adratio}, and
is important for section \ref{sec:limit}.

\section{The Limit Metric}
\label{sec:limit}

In this section we will prove that, outside the horizons $\Sigma(t)$, the 
metrics $(M^3,g_t)$ approach a Schwarzschild metric.  More
precisely, we will prove theorem \ref{thm:limit} by showing that
the rescaled $\Sigma(t)$, defined in the previous section as the original
$\Sigma(t)$ rescaled by the factor $e^{-2t}$, 
converge to a coordinate sphere of radius $M/2$ 
in $(\real^3,\delta)$ and that 
\begin{equation}
\lim_{t \rightarrow \infty} U_t(x) = 1 + \frac{M}{2|x|} 
\end{equation}
for $|x| \ge M/2$, where $M = \lim_{t \rightarrow \infty} m(t)$.

The first step is to bound $U_t(x)$ from above.  From inequality
\ref{eqn:diam_bounds} in section \ref{sec:asymptotic}, it follows that 
the rescaled $\Sigma(t)$ (defined for $t \ge t_0$) 
stay inside $S_{r_{max}}(0)$, where $r_{max} = 12C(8\pi A_0)^{1/2}$.
Hence, by the maximum principle and equation \ref{eqn:V},
\begin{equation}
   V_t(x) \le -1 + \frac{r_{max}}{|x|}.
\end{equation}
Now choose $b$ such that 
\begin{equation}
   U_{t_0}(x) \le 1 + \frac{b}{|x|}
\end{equation}
and let $c = \max(b, r_{max})$.  Then analyzing equation \ref{eqn:U_t(x)}
allows us to conclude that 
\begin{equation}\label{eqn:Uupperbound}
   U_{t}(x) \le 1 + \frac{c}{|x|}
\end{equation}
for all $t \ge t_0$ and all $x$ outside the ball of radius $r_0 e^{-2t}$.
Then since $W_t(x) \le (1 + U_t(x))/2$ by the
maximum principle and equation \ref{eqn:W}, it follows that 
\begin{equation}\label{eqn:Wupperbound}
   W_{t}(x) \le 1 + \frac{c}{2|x|}
\end{equation}
for all $t \ge t_0$ and all $x$ outside $\Sigma(t)$.

Since $W_t(x)$ is harmonic outside $S_{r_{max}}(0)$ in $(\real^3, \delta)$ 
and goes to one at infinity, it is completely determined in this region by 
its values on $S_{r_{max}+1}(0)$.  
\begin{definition}
For $\alpha \in (0,1/2)$, we  
define $H_\alpha$ to be the set of positive harmonic functions $h(x)$ 
defined outside 
$S_{r_{max}}(0)$ in $(\real^3, \delta)$ which go to one at infinity and which 
satisfy
\begin{equation}
   |h(x) - 1| \le \alpha
\end{equation}
for all $x \in S_{r_{max}+1}(0)$.
\end{definition}
Now let's put the supremum topology (for $x$ outside $S_{r_{max}+1}(0)$) 
on $H_\alpha$, so that $k(x)$ is in an $\epsilon$ neighborhood of 
$h(x)$ if and only if $|k(x) - h(x)| < \epsilon$ for all 
$x \in S_{r_{max}+1}(0)$, for
$k(x), h(x) \in H_\alpha$.  Then it follows 
that $H_\alpha$ is a compact
space with this topology \cite{GT}.

Next we define the following very useful continuous functional ${\cal F}$ on 
$H_\alpha$.  
\begin{definition}\label{def:F(h(x))}
Given $h(x) \in H_\alpha$, let $(P^3,k)$ be the Riemannian manifold 
isometric to $(\real^3 \backslash \bar{B}_{r_{max}+2}(0), h(x)^4 \delta)$.  
Then we define
\begin{equation}\label{eqn:functional}
   {\cal F}(h(x)) = \inf_{\psi(x)} \left\{ \frac{1}{4\pi}\int_{(P^3,k)}
   |\nabla \psi|^2 \; dV  \;\;\;|\;\;\; 
   \lim_{x \rightarrow \infty} \psi(x) = \psi_0 \right\} 
\end{equation}
where $\psi$ is a spinor, $\psi_0$ is a fixed constant spinor of norm one
defined at infinity,
$\nabla$ is the spin connection, and $dV$ is the volume form on $P^3$ 
with respect to the metric $k$.
\end{definition}
Furthermore, from standard theory there exists a minimizing spinor
for each $h(x)$ which satisfies 
\begin{equation}\label{eqn:elliptic_spinor}
   \left\{ \begin{array}{rl} \bar\nabla^j \nabla_j \psi(x) = 0, &  
   \mbox{ for } x \in P^3 \\
   \nu^j \nabla_j \psi = 0, &\mbox{ for } x \in \partial P^3 \\
   \lim_{x \rightarrow \infty} \psi(x) = \psi_0 &
   \end{array}\right.
\end{equation} 
where $\nabla_j$ is the spin connection  
in $(P^3,k)$, $\bar\nabla^j$ its formal 
adjoint, and $\vec{\nu}$ is the outward pointing unit normal vector to the
boundary of $P^3$.   
\begin{lemma}
The functional ${\cal F}$ is continuous on $H_\alpha$.
\end{lemma}
{\it Proof.}  
Working with respect to the standard flat metric on 
$\real^3 \backslash \bar{B}_{r_{max}+2}(0)$, we can write down explicit
formulas for the spin connection derived in \cite{F} for the case of 
two-component Weyl spinors.
Hence, we let  
$\psi(x) = (\psi^1(x), \psi^2(x))$ be a pair of complex-valued
functions and
\begin{equation}\label{eqn:EL}
   {\cal F}(h(x)) = \frac{1}{4\pi}
   \int_{\real^3 \backslash \bar{B}_{r_{max}+2}(0)}
   h(x)^2 \, |\vec{\nabla} \psi + i(\frac{\vec{\nabla} h}{h} 
   \times \vec{\sigma})\psi|^2
   dV
\end{equation}
where $\times$ is the cross product
in $\real^3$, $\sigma^i$ are the Pauli spin matrices  
\begin{equation}
   \sigma^1 = \left(\begin{array}{rr} 0 & 1 \\ 1 & 0 \end{array}\right),
   \sigma^2 = \left(\begin{array}{rr} 0 & -i \\ i & 0 \end{array}\right),
   \sigma^3 = \left(\begin{array}{rr} 1 & 0 \\ 0 & -1 \end{array}\right),
\end{equation}
$\vec{\nabla}$ is now the usual gradient in $\real^3$,
and $dV$ is the usual volume form in $\real^3$.  Using the fact that
\ref{eqn:elliptic_spinor} is the Euler-Lagrange equation for 
equation \ref{eqn:EL}, one can
show that $|\psi(x)|$ is uniformly bounded.
It also 
follows that $|\vec{\nabla}\psi(x)|
\le c/r^2$, and since $h(x)$ is harmonic, 
$|\vec{\nabla}h(x)| \le c/r^2$ too, for some uniform constant $c>0$.  

Using these facts, we can compute the 
derivative of ${\cal F}(h_t(x))$ with respect to $t$ at $t=0$.
Since the energy functional of the spinors 
in equation \ref{eqn:functional} is strictly convex, it follows that the
minimizing spinor varies smoothly for smooth variations of $h(x)$.  
Furthermore,
it follows from equation \ref{eqn:elliptic_spinor} that the contribution to
the first order rate of change of ${\cal F}(h_t(x))$ due to the variation of
the minimizing spinor is zero.  Hence, from the previous paragraph it
follows that the 
derivative of ${\cal F}(h_t(x))$ with respect to $t$ at $t=0$ 
is uniformly bounded (with respect to the supremum norm on $H_\alpha$) 
in all directions, from which it follows that ${\cal F}$ is continuous 
on $H_\alpha$.  \qed
 
\begin{definition}
Define ${\cal H}_\alpha$ to be the closure in the topological space 
$H_\alpha$ of the set of all $h(x)$ such that the corresponding manifold 
with boundary 
$(P^3,k)$ defined in the previous definition can be extended to be a 
complete, smooth, asymptotically flat manifold with
nonnegative scalar curvature (with possibly multiple ends 
but without boundary).
\end{definition}
Note that since ${\cal H}_\alpha$ is a closed subset of a compact topological
space, ${\cal H}_\alpha$ is also compact using this same supremum topology.

Using definition \ref{def:totalmass}, we can define the total mass functional
$m(h(x))$ for $h(x) \in H_\alpha$ to be the total mass of $(P^3,k)$ (defined
in definition \ref{def:F(h(x))}).  We note that $m$ is continuous on 
$H_\alpha$.
\begin{lemma}\label{lem:m>F}
Given $h(x) \in {\cal H}_\alpha$, 
\begin{equation}
   m(h(x)) \ge {\cal F}(h(x)). 
\end{equation}
\end{lemma}
{\it Proof.}  This lemma 
follows directly from Witten's proof of the positive mass
theorem \cite{Wi}, \cite{PT}.  In Witten's argument, he showed that the 
total mass is bounded below by such an integral over the whole manifold for
a spinor which satisfies the Dirac equation and goes to a constant spinor with
norm one at infinity.  Thus, the infimum of the same integral over a smaller
region and over more spinors must be smaller than the original integral.  This
proves the inequality for $h(x)$ which correspond to $(P^3,k)$ which can be 
extended smoothly.  Then the inequality follows for all $h(x) \in 
{\cal H}_\alpha$ since the total mass $m$ and the functional ${\cal F}$ are 
continuous functionals on $H_\alpha$.  \qed
\begin{lemma} 
For $h(x) \in H_\alpha $, 
${\cal F}(h(x )) = 0$ if and only if $h(x) \equiv 1 $.
\end{lemma}
{\it Proof.} 
By equation \ref{eqn:elliptic_spinor}, it follows that the norm of the 
minimizing spinor is not identically zero, from which it follows that 
${\cal F}(h(x )) = 0$ implies the existence of a parallel spinor.  
Hence, since $(P^3,k)$ is asymptotically flat, it must be flat everywhere, 
which implies that $h(x) \equiv 1$.
Conversely, if $h(x) \equiv 1$, choosing $\psi(x) \equiv \psi_0$ 
proves that ${\cal F}(h(x)) = 0$.  \qed
\begin{lemma}\label{lem:Fed}
For all $\delta > 0$, there exists an $\epsilon > 0$ such that 
\begin{equation}
{\cal F}(h(x)) < \epsilon  \;\;\;\Rightarrow\;\;\; 
\sup_{x \in S_{r_{max}+1}(0)} |h(x) - 1| < \delta
\end{equation}
for $h(x) \in {\cal H}_\alpha$.
\end{lemma}
{\it Proof.}  Since ${\cal F}$ is a continuous functional on the 
compact space ${\cal H}_\alpha$ and only equals zero if $h(x) \equiv 1$,
the lemma follows.  \qed

\begin{corollary}\label{cor:med}
For all $\delta > 0$, there exists an 
$\epsilon \in (0,2(r_{max}+1)\delta)$ such that 
\begin{equation}
m(h(x)) < \epsilon  \;\;\;\Rightarrow\;\;\; 
\sup_{x \in S_{r_{max}+1}(0)} |h(x) - 1| < \delta
\end{equation}
for $h(x) \in {\cal H}_\alpha$.
\end{corollary}
{\it Proof.}  The fact that there exists an $\epsilon > 0$ follows directly 
from lemmas \ref{lem:m>F} and \ref{lem:Fed}.  Then by considering   
$h(x) = 1 + a/2|x|$ for $a > 0$ which corresponds to a 
Schwarzschild metric with total mass $a$, it follows that 
$\epsilon < 2(r_{max}+1)\delta$.  \qed

Alternatively, the previous discussion with spinors beginning with 
definition \ref{def:F(h(x))} 
and ending with the above corollary could be replaced
by quoting theorem 1.1 of \cite{BF} (the proof of which also uses spinors) 
which proves corollary \ref{cor:med} as well.  

Now we are ready to apply these results to understand the asymptotic
behavior of $U_t(x)$ as $t$ goes to infinity.  
\begin{theorem}\label{thm:Ulimit}
For all $\delta > 0$, there exists a $\bar{t}$ such that for 
all $t \ge \bar{t}$
\begin{equation}
   |U_t(x) - \left(1 + \frac{M}{2|x|} \right)| \le \delta
\end{equation}
for $|x| \ge r_{max}+1$, where 
$M \equiv \lim_{t \rightarrow \infty} m(t) > 0$.  
\end{theorem}
{\it Proof.}
First, we let 
\begin{equation}
   \bar{U}_t(x) = U_t(x) - \left(1 + \frac{m(t)}{2|x|} \right)
\end{equation}
and 
\begin{equation}\label{eqn:Wbar}
   \bar{W}_t(x) = W_t(x) - \left(1 + \frac{\tilde{m}(t)}{2|x|} \right)
\end{equation}
so that by lemmas \ref{lem:m} and \ref{lem:mtilde} the harmonic functions 
$\bar{U}_t(x)$ and $\bar{W}_t(x)$ do not have a constant term or 
a $1/|x|$ term in their expansions using spherical harmonics.  Then 
substituting these two expressions into equation \ref{eqn:U_t(x)2} yields
\begin{equation}\label{eqn:Ubar_t(x)}
 \frac{d}{dt} \bar{U}_t(x) = 2 ( \bar{U}_t(x) - \bar{W}_t(x) + 
                         r\frac{\partial}{\partial r} \bar{U}_t(x) )
\end{equation}
by equation \ref{eqn:mprime2}.  

Also, by equation \ref{eqn:mprime2} and the positive mass theorem, 
$\int_0^\infty \tilde{m}(t) \; dt \le m(0)/2$,
so $T_{bad} = \{ t \;|\; \tilde{m}(t) \ge \epsilon \}$ has finite 
measure less than or equal to $m(0)/2\epsilon$.  Hence, by lemma 
\ref{lem:mtilde} and corollary \ref{cor:med}, for all $\delta > 0$, 
there exists an $\epsilon > 0$ such that 
\begin{equation}
  \sup_{x \in S_{r_{max}+1}(0)} |W_t(x) - 1| < \delta
\end{equation}
and hence  
\begin{equation}\label{eqn:Wbarbound}
  \sup_{x \in S_{r_{max}+1}(0)} |\bar{W}_t(x)| < 2 \delta
\end{equation}
since 
\begin{equation}\label{eqn:mtildebound}
\tilde{m}(t) < \epsilon < 2(r_{max}+1) \delta,
\end{equation}
for all $t \notin T_{bad}$.  Also, for all $t$ we have the uniform bound 
\begin{equation}
  \sup_{x \in S_{r_{max}+1}(0)} |\bar{W}_t(x)| < B
\end{equation}
(where $B = 1 + c/2(r_{max}+1)$) by equations \ref{eqn:Wbar}
and \ref{eqn:Wupperbound} and since $\tilde{m}(t) \le c$ by equation 
\ref{eqn:Wupperbound}.  Then since $\bar{W}_t(x)$ is a harmonic function
without a constant term or a $1/|x|$ term, it follows that
\begin{equation}
   |\bar{W}_t(x)| \le k \frac{(r_{max}+1)^2}{|x|^2}
    \left\{\begin{array}{rl}
   2\delta & \mbox{ for } t \notin T_{bad} \\ 
   B & \mbox{ for } t \in T_{bad} \end{array}\right.
\end{equation}
for $|x| \ge r_{max}+1$ and some positive constant $k$ (which we will not 
need to compute).   
Then from analyzing equation \ref{eqn:Ubar_t(x)} and using the fact
that $\bar{U}_t(x)$ does not have a constant term or a $1/|x|$ term and the
fact that $T_{bad}$ has finite measure, 
it follows that we can 
choose some $\bar{t}$ large enough such that 
\begin{equation}
   |\bar{U}_t(x)| \le k \frac{(r_{max}+1)^2}{|x|^2} (2.01) \delta
\end{equation}
for all $t \ge \bar{t}$ and $|x| \ge r_{max}+1$.  
Hence, since $\delta > 0$ was arbitrary and 
$M = \lim_{t \rightarrow \infty} m(t)$ (and is positive by
inequality \ref{eqn:lower_mass2}), the theorem follows.  \qed 

\begin{corollary}
For all $\delta > 0$, there exists a $\bar{t}$ such that 
\begin{equation}
   |V_t(x) - \left(\frac{M}{2|x|} - 1\right)| \le \delta
\end{equation}
for $|x| \ge r_{max}+1$ and all $t \ge \bar{t}$ except on a set with  
measure less than $\delta$.
\end{corollary}
{\it Proof.}
We recall from equation \ref{eqn:V=U-2W} that 
$V_t(x) = U_t(x) - 2 W_t(x)$.  Hence, since $U_t(x)$ is converging to 
$1 + \frac{M}{2|x|}$ 
by theorem \ref{thm:Ulimit} and $W_t(x)$ is converging to $1$
(for $t \notin T_{bad}$) by equations \ref{eqn:Wbar}, 
\ref{eqn:Wbarbound}, and \ref{eqn:mtildebound}, 
then the corollary follows from 
the fact that $T_{bad}$ has finite measure.  \qed

\begin{corollary}\label{cor:Vlimit}
For all $\delta > 0$, there exists a $\bar{t}$ such that for 
all $t \ge \bar{t}$,
\begin{equation} \label{eqn:VVlimit}
   |V_t(x) - \left(\frac{M}{2|x|} - 1\right)| \le \delta
\end{equation}
and
\begin{equation} \label{eqn:WWlimit}
   |W_t(x) - 1| \le \delta
\end{equation}
for $|x| \ge r_{max}+1$.
\end{corollary}
{\it Proof.}  
Equation \ref{eqn:VVlimit} follows from the previous corollary, 
equation \ref{eqn:V}, and theorem \ref{thm:equal}.  Equation \ref{eqn:WWlimit}
then follows from equation \ref{eqn:VVlimit}, theorem \ref{thm:Ulimit}, 
and equation \ref{eqn:V=U-2W}.  \qed 

The next two lemmas use corollary \ref{cor:Vlimit} to prove that $\Sigma(t)$
converges to the sphere of radius $M/2$.  The main idea is that since 
$V_t(x)$ converges to $\frac{M}{2|x|} - 1$ and equals zero on $\Sigma(t)$
by definition, then $\Sigma(t)$ must be converging to $S_{M/2}(0)$.

In the rest of this section we will want to take limits of certain 
sequences of surfaces, and naturally there are several ways to do this.
However, in
our case, all of the surfaces we are dealing with are boundaries of regions,
so it seems most natural to follow \cite{MM}.  In \cite{MM}, a very general
definition of the measure of the perimeter of a Lebesgue measurable set is
given on p. 64.  Then on p. 70, it is shown that the space of Lebesgue 
measurable sets with equally bounded perimeters in a compact region $K$ 
is compact with respect to the $L^1$ norm of the characteristic functions
of the regions, meaning that any sequence of such sets 
contains a subsequence such that the 
characteristic functions of the subsequence converge in $L^1(K)$.  
Thus, in what follows, we will 
say that $\Sigma_{\infty}$ is a limit of a sequence of 
surfaces $\{\Sigma(t_i)\}$ if it is the boundary 
of a Lebesgue measurable set with bounded perimeter 
whose characteristic function is the $L^1$ limit 
of the characteristic functions of a subsequence of the $\{\Sigma(t_i)\}$.
We note that \ref{eqn:area_bounds} implies that the 
$\{\Sigma(t_i)\}$ have equally bounded perimeters for $t_i \ge t_0$.
We also note that since this perimeter function is shown to be
lower semicontinuous in \cite{MM}, 
$\Sigma_{\infty}$ has bounded perimeter.
(We also note that we define the boundary of a region to be 
the set of boundary points such that every open ball around a boundary point
contains a positive measure of both the region and the complement of the 
region.)  

From the above considerations, it also turns out that the 
$\{\Sigma(t_i)\}$ will converge to $\Sigma_{\infty}$ in the Hausdorff 
distance sense as well.  This follows from the fact that each 
$\Sigma(t_i)$ minimizes area in 
$(\real^3 \backslash B_{r_0 e^{-2t}}(0), U_t(x)^4 \delta)$
and we have uniform upper and lower bounds on the $U_{t_i}(x)$.
(This is related to the proof of lemma \ref{lem:c1r2} in appendix 
\ref{sec:regularity}.)
\begin{lemma}\label{lem:limit1}
Let $\Sigma_\infty$ be any limit of $\Sigma(t)$ in $(\real^3,\delta)$.
Then no part of $\Sigma_\infty$ lies inside $S_{M/2}(0)$.
\end{lemma}
{\it Proof.}
Let $R$ be the open region outside $\Sigma_\infty$, and 
consider a sequence of $t_i$ going to infinity
such that $\Sigma(t_i)$ is converging to a limit $\Sigma_\infty$.
Since each $\Sigma(t_i)$ is outer-minimizing and we
have uniform upper and lower bounds on $U_{t_i}(x)$ (and hence on the 
corresponding metric), it follows that $R$ must be connected.
Then $V_{t_i}(x)$ must converge to a harmonic function in $R$, and 
by corollary \ref{cor:Vlimit} this limit harmonic function 
equals $\frac{M}{2|x|} - 1$.  But each $V_{t_i}(x) \le 0$
by definition, so $R \cap B_{M/2}(0) = \varnothing$, proving the lemma.  \qed
\begin{lemma}\label{lem:limit2}
Let $\Sigma_\infty$ be any limit of $\Sigma(t)$ in $(\real^3,\delta)$.	
Then no part of $\Sigma_\infty$ lies outside $S_{M/2}(0)$.
\end{lemma}
{\it Proof.}
Suppose otherwise.  Then we can consider a sequence of $t_i$ going to infinity
such that $\Sigma(t_i)$ is converging to a limit $\Sigma_\infty$
which lies at least 
partially outside $S_{M/2}(0)$. 
Again, defining the region $R$ as above, 
$V_{t_i}(x)$ must converge to a harmonic function in $R$, and 
by corollary \ref{cor:Vlimit} this limit harmonic function 
equals $\frac{M}{2|x|} - 1$.  

Since $V_{t_i}(x)$ is harmonic, it minimizes its energy among functions with
the same boundary data.  Thus, since $\Sigma(t_i)$ is contained inside 
$S_{r_{max}}(0)$, the energy of $V_{t_i}(x)$ is less than 
$4\pi r_{max}$.  Now choose $x_0 \in \Sigma_\infty$ which maximizes distance
from the origin, and let 
$|x_0| = M/2 + 2r$ for some $r>0$.  
Then by the co-area formula,   
\begin{eqnarray}
   4\pi r_{max} &>&
   \int_{R_i} |\nabla V_{t_i}(x)|^2 \; dV 
   > \int_{R_i \cap B_r(x_0)} |\nabla V_{t_i}(x)|^2 \; dV  \\
   &=& \int_{-1}^0 dz \int_{L_z \cap B_r(x_0)} |\nabla V_{t_i}(x)| \; dA_z \\
   &\ge&  \int_{-1}^0 dz \, |L_z \cap B_r(x_0)|^2 
   \left(\int_{L_z \cap B_r(x_0)} |\nabla V_{t_i}(x)|^{-1}  
   \; dA_z \right)^{-1} \label{eqn:inverse}
\end{eqnarray}
where $R_i$ is the region outside $\Sigma(t_i)$, $dV$ is the standard volume
form on $\real^3$, $L_z$ is the level set on which $V_{t_i}$ equals $z$, 
$dA_z$ is the area form of $L_z$, and $|\cdot|$ denotes area
in $\real^3$.  Next, we define
\begin{equation}
   {\cal V}(z) = \mbox{ the volume of the region } \{x \in B_r(x_0) 
   \;|\; V_{t_i}(x) > z \}
\end{equation}
for $-1 < z \le 0$.  Then from inequality \ref{eqn:inverse},
\begin{equation}\label{eqn:contra}
   4\pi r_{max} > \int_{-1}^0 |L_z \cap B_r(x_0)|^2 \, {\cal V}'(z)^{-1} \; dz.
\end{equation}

Since $V_{t_i}(x)$ equals zero on $\Sigma(t_i)$ and yet is converging 
to $\frac{M}{2|x|} - 1$ in $R$, we must have ${\cal V}'(z)$ approaching zero
and the surfaces $L_z \cap B_r(x_0)$ converging to 
$\Sigma(t_i) \cap B_r(x_0)$ for $z \in (\frac{M}{2|x_0|} - 1 ,0)$ 
as $t_i$ goes to infinity.  However, since $\Sigma(t_i)$ minimizes area in
$(\real^3 \backslash B_{r_0 e^{-2t_i}}(0), U_{t_i}(x)^4 \delta)$ and 
$U_{t_i}(x)$ is uniformly bounded, 
it follows that $|L_z \cap B_r(x_0)|$
is not going to zero for $z \in (\frac{M}{2|x_0|} - 1 ,0)$.  Hence, the 
right hand side of inequality \ref{eqn:contra} is going to infinity as 
$t_i$ goes to infinity, giving us a contradiction and proving the lemma.  
\qed

\begin{theorem}
   The surfaces $\Sigma(t)$ converge to $S_{M/2}(0)$ in the 
Hausdorff distance sense in the limit as 
$t$ goes to infinity,  
\begin{equation}\label{eqn:Ulimit1}
   \lim_{t \rightarrow \infty} U_t(x) = 1 + \frac{M}{2|x|}
\end{equation}
for $|x| \ge M/2$, and
\begin{equation}\label{eqn:Ulimit2}
   \lim_{t \rightarrow \infty} U_t(x) = \sqrt{\frac{2M}{|x|}}
\end{equation}
for $0 < |x| \le M/2$, where as usual 
$M = \lim_{t \rightarrow \infty} m(t) > 0$.
\end{theorem}
{\it Proof.}
First we observe that for all $\delta > 0$, there exists a $\bar{t}$ such that
\begin{equation}\label{eqn:Sigma_annulus}
   \Sigma(t) \subset B_{\frac{M}{2}+ \delta}(0) \backslash 
   B_{\frac{M}{2} - \delta}(0)
\end{equation}
for all $t \ge \bar{t}$.  Otherwise, we could choose a sequence of $t_i$
going to infinity such
that each $\Sigma(t_i)$ did not lie entirely in the annulus.  
By previous discussions, a limit $\Sigma_{\infty}$ would have to exist, and 
since this limit is valid in the Hausdorff distance sense as previously 
discussed, at least part of $\Sigma_{\infty}$ would have to lie off of 
$S_{M/2}(0)$, contradiction at least one of the two previous lemmas.
This proves 
equation \ref{eqn:Sigma_annulus}.  

As a corollary, we get that for all $\delta > 0$, there exists a $\bar{t}$
such that 
\begin{equation}\label{eqn:Vbounds}
   F_{\frac{M}{2}-\delta}(x) \le  V_t(x) \le F_{\frac{M}{2}+\delta}(x)     
\end{equation} 
where we define
\begin{equation}
   F_{a}(x) = \left\{\begin{array}{rl}\frac{a}{|x|} - 1 & \mbox{ for }
   |x| \ge a \\ 0  & \mbox{ for } |x| \le a \end{array}\right.      
\end{equation} 
for $a>0$.  Then equations \ref{eqn:Ulimit1} and \ref{eqn:Ulimit2} follow
(with uniform convergence on compact subsets of $\real^3 - \{0\}$) 
from inequality \ref{eqn:Vbounds} and analyzing the behavior of equation
\ref{eqn:U_t(x)}.  \qed

Theorem \ref{thm:limit} then follows from the above theorem and corollary 
\ref{cor:bounded}.

\section{Generalization to Asymptotically Flat \\ Manifolds and the 
Case of Equality}\label{sec:generalization}

Up to this point in the paper 
we have assumed that $(M^3,g_0)$ was harmonically flat
at infinity.  In particular, theorems \ref{thm:existence}, \ref{thm:monotone},
and \ref{thm:limit} only apply to harmonically flat manifolds as stated.
In this section, we will extend theorems 
\ref{thm:existence} and \ref{thm:monotone} and elements of theorem 
\ref{thm:limit} to asymptotically
flat manifolds.  This will prove the main theorem, theorem \ref{thm:Penrose},
except for the case of equality, which we will see follows 
from the case of equality of theorem \ref{thm:green2}.

It is worth noting that the main reason for initially considering only
harmonically flat manifolds was convenience.  Alternatively, we 
could have ignored harmonically flat manifolds and dealt only with 
asymptotically flat manifolds.  However, this would have complicated some of 
the arguments unnecessarily, so we chose to delay these considerations 
until now.

\begin{definition}\label{def:asymptotically_flat}
$(M^n,g)$ is said to be {\bf asymptotically flat} if 
there is a compact set $K\subset
M$ such that $M \backslash K$ is the disjoint union of ends 
$\{E_k\}$, such that
for each end there exists 
a diffeomorphism $\Phi_k:E_k \to \real^n \backslash \ B_1(0)$ such that, in
the coordinate chart defined by $\Phi_k$,
\[g=\sum_{i,j}g_{ij}(x)dx^i dx^j\]
where
\[g_{ij}(x)=\delta_{ij}+O(|x|^{-p})\]
\[|x||g_{ij,k}(x)| + |x|^2|g_{ij,kl(x)}| = O(|x|^{-p})\]
\[|R(g)|=O(|x|^{-q})\]
for some $p>\frac{n-2}{2}$ and some $q>n$, where we have used commas
to denote partial derivatives in the coordinate chart, and $R(g)$ is
the scalar curvature of $(M^n,g)$.
\end{definition}
These assumptions on the asymptotic behavior of $(M^n,g)$ at infinity
imply the existence of the limit
\begin{equation}\label{eqn:ADM_mass}
M_{ADM}(g)= \frac{1}{16\pi} \lim_{\sigma\to\infty}
\int_{S_\sigma}\sum_{i,j}(g_{ij,i}\nu_j-g_{ii,j}\nu_j)\,d\mu
\end{equation}
where $\omega_{n-1}=Vol(S^{n-1}(1))$, $S_\sigma$ is the coordinate sphere
of radius $\sigma$, $\nu$ is the unit normal to $S_\sigma$, and $d\mu$ 
is the area element of $S_\sigma$ in the coordinate chart.  
The quantity $M_{ADM}$ is called the {\bf total mass} of
$(M^n,g)$ (see \cite{ADM}, \cite{Ba2}, \cite{S}, and \cite{SY5}), and agrees
with the definition of total mass for harmonically flat 3-manifolds
given in definition \ref{def:totalmass}.

First we observe that
the arguments in the proof of the existence theorem, theorem
\ref{thm:existence}, 
did not use harmonic flatness anywhere, so we immediately get existence
of the conformal flow of metrics for asymptotically flat manifolds.
Similarly, the arguments used in section \ref{sec:A(t)} 
to prove that $A(t)$ is constant still hold.
Next we reexamine the proof of theorem \ref{thm:massnonincreasing}
which proved that $m(t)$ was nonincreasing.  The only 
modification we need to make is to use the more general definition for the
total mass of an asymptotically flat manifold given by equation 
\ref{eqn:ADM_mass}.
It is then straight forward 
to check that equation \ref{eqn:mprime} and hence theorem
\ref{thm:massnonincreasing} are still true.  Hence,
\begin{theorem}
Theorems \ref{thm:existence} and \ref{thm:monotone} 
are true for asymptotically flat manifolds as well as 
harmonically flat manifolds.
\end{theorem}

We choose not to extend theorem \ref{thm:limit} to asymptotically flat
manifolds, but conjecture that it is still true.  Instead, we observe
that we must still have 
\begin{equation}\label{eqn:pi}
m(t) \ge \sqrt{\frac{A(t)}{16\pi}}
\end{equation}
for asymptotically flat
manifolds.  Otherwise, given an asymptotically flat counterexample, we could
use lemma \ref{lem:sy} to perturb the manifold slightly making it harmonically
flat at infinity such that it still violated equation \ref{eqn:pi}.  Then 
applying the conformal flow of metrics to this harmonically flat manifold
would violate theorems \ref{thm:monotone} and \ref{thm:limit}, which is a 
contradiction.  Setting $t=0$ in inequality \ref{eqn:pi} then 
proves the Riemannian Penrose inequality for 
asymptotically flat manifolds.

The case of equality of theorem \ref{thm:Penrose} then follows from 
equation \ref{eqn:mprime} and theorem \ref{thm:green2}.  If we have equality
in the Riemannian Penrose inequality, then applying the conformal flow of 
metrics to this initial metric must also give equality in inequality 
\ref{eqn:pi} for all $t \ge 0$.  Hence, the right hand derivative of $m(t)$
at $t=0$ equals zero, so by equation \ref{eqn:mprime}, 
\begin{equation}
   {\cal E}(\Sigma^+(0),g_0) \;\; = \;\; 2m(0).
\end{equation}
By definition \ref{def:limits} and equation \ref{eqn:ODE2}, 
$\Sigma^+(0)$ is the outermost minimal area enclosure of $\Sigma_0$
in $(M^3,g)$.  Furthermore, by the case of equality of theorem 
\ref{thm:green2}, $(M^3,g)$ is a Schwarzschild manifold outside
$\Sigma^+(0)$.  Hence, $\Sigma^+(0)$ is the outermost horizon of $(M^3,g)$,
so $(M^3,g)$ is isometric to a Schwarzschild manifold outside their
respective outermost horizons.  This completes the proof of theorem 
\ref{thm:Penrose} and the Riemannian Penrose inequality. 

The reader might also have noticed that none of the arguments in this paper
have used anything about the original manifold inside the original horizon
$\Sigma_0$.  Hence, we can generalize theorem 
\ref{thm:Penrose} to the following.
\begin{theorem}\label{thm:Penrose_boundary}
Let $(M^3,g)$ be a complete, smooth, asymptotically flat   
$3$-manifold with boundary which has nonnegative scalar curvature and total 
mass $m$.  Then if the boundary is an outer-minimizing horizon (with one or
more components) of total area $A$, then
\begin{equation}\label{eqn:Penrose_boundary}
  m \ge \sqrt{\frac{A}{16\pi}} 
\end{equation}
with equality if and only if $(M^3,g)$ is isometric to a 
Schwarzschild manifold outside their respective outermost horizons.
\end{theorem}

\section{New Quasi-Local Mass Functions}\label{sec:qlm}

Given a region in a space-like slice of a space-time, it is natural to ask
how much energy and momentum is contained in that region.  As described in 
the introduction, there does exist a well-defined notion of total mass and 
also of energy-momentum density.  
However, when the region in question is some finite
region which is not the entire manifold or just a single point, it is not
very well understood how to define how much energy and momentum is in that
region.

Various definitions of ``quasi-local'' mass exist, such as the Hawking mass,
which was used by Huisken and Ilmanen in \cite{HI} to prove their Riemannian 
Penrose inequality, for example.  Good definitions of quasi-local mass 
should satisfy certain reasonable properties \cite{CY} such as positivity and
some kind of monotonicity, either under a flow or by inclusion.  In addition,
for a large region, the quasi-local mass function should approach the total
mass of the manifold, and for horizons with area $A$ it is thought that the 
mass should be $\sqrt{A/16\pi}$.

Let $(M^3,g)$ be a complete asymptotically flat manifold with nonnegative
scalar curvature and total mass $m$.  Let $\Sigma$ be any surface in $M^3$
which is in the class of surfaces ${\cal S}$ defined in section 
\ref{sec:definitions}.

\begin{definition}\label{def:qlm1}  
Suppose $u(x)$ is a 
positive harmonic function in $(M^3,g)$ outside $\Sigma$
going to a constant at infinity scaled such that $(M^3,u(x)^4 g)$ has the 
same total mass as $(M^3, g)$. 

Then if  $\Sigma$ is an outer-minimizing 
horizon with area $A$ in $(M^3, u(x)^4 g)$, we define
the quasi-local mass of $\Sigma$ in $(M^3, g)$ to be 
\begin{equation}
   m_g(\Sigma) = \sqrt{\frac{A}{16\pi}}.
\end{equation}
\end{definition}
\begin{definition}
We define $S$ to be the subset of ${\cal S}$ of surfaces $\Sigma$ for which 
such a conformal factor $u(x)$ exists, and we note that 
(by equation \ref{eqn:identity} mostly) this conformal factor
is unique for each $\Sigma$ when it exists.
\end{definition}
As usual $\Sigma$ could have multiple components.
It is also interesting that 
\begin{equation}
   m_g(\tilde{\Sigma}) = m_{u(x)^4 g}(\tilde{\Sigma}) 
\end{equation}
for all surfaces $\tilde{\Sigma} \in S$
where the conformal factor $u(x)$ is any harmonic function 
in $(M^3,g)$ defined outside $\tilde{\Sigma}$ which goes to a constant at 
infinity scaled such that $(M^3,u(x)^4 g)$ has the 
same total mass as $(M^3, g)$. 

\begin{lemma}
The quasi-local mass function $m_g(\Sigma)$ defined for $(M^3,g)$ 
is nondecreasing for the family of surfaces $\Sigma(t)$ defined by
equation \ref{eqn:ODE2}.  That is, $m_g(\Sigma(t))$ is nondecreasing in $t$.
Furthermore, 
\begin{equation}
   m_g(\Sigma(0)) = \sqrt{\frac{A}{16\pi}}
\end{equation}
where $A$ is the area of the original outer-minimizing horizon $\Sigma_0$
in $(M^3,g)$, and 
\begin{equation}
   \lim_{t \rightarrow \infty} m_g(\Sigma(t)) = m,
\end{equation}
the total mass of $(M^3,g)$.
\end{lemma}
{\it Proof.}  We consider the conformal flow of metrics $g_t$ beginning 
with $(M^3,g)$ as discussed throughout this paper and in the previous section
for asymptotically flat manifolds.
Then we note that by equation \ref{eqn:ODE2}, $\Sigma(t)$ is an 
outer-minimizing horizon in $(M^3,u_t(x)^4 g)$ with area $A(t)$.  
Hence, $u_t(x)$ satisfies the conditions in definition \ref{def:qlm1} except 
that it is not scaled to have the correct mass.  Hence, since mass has units 
of length, it follows that
\begin{equation}
   m_g(\Sigma(t)) = \frac{m}{m(t)}\sqrt{\frac{A(t)}{16\pi}}
\end{equation}
where again 
$m$ is the total mass of $(M^3,g)$.  Then the lemma follows from theorems
\ref{thm:monotone} and \ref{thm:limit} and the fact that $m(0) = m$.  \qed

There is a trick which allows us to extend the definition of this quasi-local
mass function to all surfaces in ${\cal S}$.
\begin{definition}
Define
\begin{equation}
   \tilde{m}_g(\Sigma) = \sup \left\{m_g(\tilde{\Sigma}) \;\;|\;\;
   \Sigma \mbox{ (entirely) encloses } \tilde{\Sigma} \in S \right\}  
\end{equation}
where $\Sigma$ is any surface in ${\cal S}$.
\end{definition}
It follows trivially that $\tilde{m}_g(\Sigma)$ is monotone with respect to 
enclosure, which is a desirable property for quasi-local mass functions 
to have since larger regions should contain more mass in the nonnegative 
energy density setting which we are in.
We also note that this same construction can be used with the Hawking mass to 
make it monotone with respect to enclosure too, where the original Hawking mass
should only be defined to exist for surfaces which equal their own maximal 
minimal area enclosures, as motivated by the results of Huisken and Ilmanen 
in \cite{HI}.

We also notice that the existence of the Penrose inequality allows us to 
define another new quasi-local mass function which is similar in nature to the
Bartnik mass \cite{Ba4}, \cite{HI}.  
In fact, whereas the Bartnik mass could also be called an 
outer quasi-local mass function, it makes sense to call 
the new quasi-local mass function defined below
the inner quasi-local mass function, which is clear from the definition.  

\begin{definition}\label{def:qlm3}
Given a surface $\Sigma \in {\cal S}$ in $(M^3,g)$, consider 
all other asymptotically flat, complete, Riemannian manifolds 
$(\tilde{M}^3,\tilde{g})$ with nonnegative scalar curvature  
which are isometric to $(M^3,g)$ outside $\Sigma$.   Then we define
\begin{equation}
   m_{inner}(\Sigma) = \sup \sqrt{\frac{\tilde{A}}{16\pi}}
\end{equation}
where $\tilde{A}$ is the infimum of
the areas of all of the surfaces in $(\tilde{M}^3,\tilde{g})$
in ${\tilde{\cal S}}$.
\end{definition}
We note here that $\tilde{\cal S}$ is defined the same way as ${\cal S}$
in definition \ref{def:calS} and that the surface in $\tilde{\cal S}$ with
minimum area may have multiple components.  We also note that for $\tilde{A}$
to be nonzero that $(\tilde{M}^3,\tilde{g})$ must have more than one 
asymptotically flat end.

\begin{lemma}
Let $(M^3,g)$ be an asymptotically flat, complete, Riemannian manifold with 
nonnegative scalar curvature, and let $\Sigma_1, \Sigma_2 \in {\cal S}$ 
such that $\Sigma_2$ (entirely) encloses $\Sigma_1$.  Then
\begin{equation}
   m \ge m_{inner}(\Sigma_2) \ge m_{inner}(\Sigma_1)
\end{equation}  
where $m$ is the total mass of $(M^3,g)$.
\end{lemma}
{\it Proof.}  Follows directly from the Penrose inequality and definition
\ref{def:qlm3}.  \qed

Also, if $\Sigma$ is outer-minimizing 
(see definition \ref{def:outerminimizing}), then 
\begin{equation}
   m_{outer}(\Sigma) \ge m_{inner}(\Sigma)
\end{equation}
where $m_{outer}(\Sigma)$ is basically the Bartnik mass \cite{Ba4}
except that we only consider extensions of the metric in which 
$\Sigma$ continues to be outer-minimizing.  The proof of this inequality 
and related discussions 
will be included in a paper on quasi-local mass which is
currently in progress.

\section{Open Problems and Acknowledgments}\label{sec:acknowledgments}

Even though the original Penrose conjecture concerned only three dimensional, 
space-like slices of a space-time, it is easy to generalize the conjecture 
to higher dimensions using the same motivation as the three-dimensional case.
In this paper we have restricted our attention to proving 
the three dimensional 
case of the Riemannian Penrose conjecture, which is perhaps the most physically
interesting dimension.
However, it appears to the author that 
the techniques presented here generalize to higher dimensions, which will
be the subject of a future paper.

Of course, an even more important problem is to 
prove the Penrose conjecture for arbitrary space-like slices 
(as opposed to totally-geodesic slices treated here) of 
a space-time as described in section 
\ref{sec:introduction}.  There seem to be several natural ideas to try in
this regard.  First, one could attempt to solve a variant of Jang's equation
which was used to extend the Riemannian positive mass theorem to the general
case \cite{SY5}.  Otherwise, one could try to modify the flow of 
metrics defined in this paper to define a flow of the Cauchy data $(M^3,g,h)$ 
which has good monotonicity properties.  Both of these approaches are 
speculative but could yield interesting results.

We also note that in this paper we did not ever 
prove the uniqueness of the conformal
flow of metrics defined by equations \ref{eqn:ODE1}, \ref{eqn:ODE2}, 
\ref{eqn:ODE3}, and \ref{eqn:ODE4}, although we conjecture that this is true.
It would also be interesting to understand the relationship between the 
Hawking mass used in the Huisken/Ilmanen paper and the quasi-local mass 
function $m_g$ defined in the previous section.  In the spherically symmetric
case, $m_g$ is bounded below by the Hawking mass on the spherically symmetric
spheres since the Hawking mass
equals the Bartnik mass in this case (outside the outermost horizon).  It is
unclear if such a relationship generally holds.  

I have many people to thank who have helped make 
this paper possible.  I am very 
fortunate to have been  
introduced to the Penrose conjecture by Richard Schoen, who 
has shared many of his thoughts on the problem with me and has given me much
helpful feedback about this paper.  I would also like to  
thank Shing-Tung Yau, who hosted me during my stay at Harvard University
during the academic year 1997-98, 
and the NSF (grant \#DMS-9706006)
for supporting my research 
during that time in which the idea of this paper materialized.
I am also indebted to David Jerison, Leon Simon, and  
Gang Tian who all gave me very helpful 
suggestions for working out many of the important details of this paper.  
I also thank
Robert Bartnik for 
explaining the Bunting/Masood-ul-Alam paper \cite{BM} to me, which turned out
to be very useful.  I must also thank Felix Finster and Kevin Iga who both
helped me think about certain aspects of this paper so much that in both cases
it resulted in other papers \cite{BF}, \cite{BI}.  Finally, I would like to 
thank Gerhard Huisken, Tom Ilmanen, Dan Stroock, and Brian White for
enthusiastic and helpful conversations.

\newpage
\appendix

\section{The Harmonic Conformal Class of a Metric}
\label{sec:harmonic}

In this appendix we define a new equivalence class
and partial ordering of conformal metrics.  This    
provides a natural motivation for studying 
conformal flows of metrics to try to prove the Riemannian Penrose inequality.

Let 
\begin{equation}\label{eqn:conf}
g_2 = u(x)^{\frac{4}{n-2}} g_1
\end{equation}
where $g_2$ and $g_1$ are metrics on 
an $n$-dimensional manifold $M^n$, $n \ge 3$.  Then we get the surprisingly 
simple identity that
\begin{equation}\label{eqn:identity}
\Delta_{g_1}(u \phi) = u^\frac{n+2}{n-2} \Delta_{g_2}(\phi) 
                       + \phi \Delta_{g_1}(u)
\end{equation}
for any smooth function $\phi$.  

This motivates us to define the following relation.
\begin{definition}
Define
\[ g_2 \sim g_1 \]
if and only if 
equation \ref{eqn:conf} is satisfied with $\Delta_{g_1}(u) = 0$ and $u(x) > 0$.
\end{definition}
Then from equation \ref{eqn:identity} we get the following lemma.
\begin{lemma}
The relation $\sim$ is reflexive, symmetric, and transitive, and hence 
is an equivalence relation.
\end{lemma}
Thus, we can define the following equivalence class of metrics.
\begin{definition}
Define 
\[ [g]_H = \{\bar{g} \;\;|\;\; \bar{g} \sim g \} \]
to be the {\bf harmonic conformal class} of the metric $g$.
\end{definition}
Of course, this definition is most interesting when $(M^n,g)$ has 
nonconstant positive harmonic functions, which happens for example
when $(M^n,g$) has a boundary.  

Also, we can modify the relation $\sim$ to get another relation $\succeq$.
\begin{definition}
Define
\[ g_2 \succeq g_1 \]
if and only if equation \ref{eqn:conf} is satisfied
with $- \Delta_{g_1}(u) \ge 0$ and $u(x) > 0$.
\end{definition}
Then from equation \ref{eqn:identity} we get the following lemma.
\begin{lemma}
The relation $\succeq$ is reflexive and transitive, and hence 
is a partial ordering.
\end{lemma}
Since $\succeq$ is defined in terms of superharmonic functions, we 
will call it the superharmonic partial ordering of metrics on $M^n$.
Then it is natural to define the following set of metrics.
\begin{definition}
Define 
\[ [g]_S = \{\bar{g} \;\;|\;\; \bar{g} \succeq g \} .\]
\end{definition}
This set of metrics has the property that if $\bar{g} \in [g]_S$, then
$[\bar{g}]_S \subset [g]_S$

Also, the scalar curvature transforms nicely under a conformal change of 
the metric.  In fact, assuming equation \ref{eqn:conf} again, 
\begin{equation}\label{eqn:scalar_curv}
  R(g_2) = u(x)^{-(\frac{n+2}{n-2})}\left(-c_n \Delta_{g_1} + R(g_1)\right)u(x)
\end{equation}
where $c_n = \frac{4(n-1)}{n-2}$ \cite{S}.  This gives us the following lemma.
\begin{lemma}
The sign of the scalar curvature is preserved pointwise by $\sim$.  That is,
if $g_2 \sim g_1$, then $sgn(R(g_2)(x)) = sgn(R(g_1)(x))$ for all $x \in M^n$.

Also, if $g_2 \succeq g_1$, and $g_1$ has non-negative scalar curvature, then
$g_2$ has non-negative scalar curvature.  
\end{lemma}

Hence, the harmonic conformal equivalence relation $\sim$ and the superharmonic
partial ordering $\succeq$ are useful for studying questions about 
scalar curvature.  In particular, these notions are useful for studying the
Riemannian Penrose inequality which concerns asymptotically flat 3-manifolds
$(M^3,g)$ with non-negative scalar curvature.  Given such a manifold, define
$m(g)$ to be the total mass of $(M^3,g)$ and $A(g)$ to be the area of the 
outermost horizon (which could have multiple components) of $(M^3,g)$.  Define
$P(g) = \frac{m(g)}{\sqrt{A(g)}}$ to be the Penrose quotient of $(M^3,g)$.
Then an interesting question is to ask which metric in $[g]_S$ minimizes
$P(g)$.      

This paper can be viewed as an answer to the above question.  
We showed that
there exists a conformal flow of metrics (starting with $g_0$) 
for which the Penrose quotient
was non-increasing, and in fact this conformal flow stays inside $[g_0]_S$.
Furthermore, $g_{t_2} \in [g_{t_1}]_S$ for all $t_2 \ge t_1 \ge 0$.  We showed
that no matter which metric we start with, 
the metric converges to a Schwarzschild metric 
outside its horizon.  Hence, the minimum value of $P(g)$ in $[g]_S$ is 
achieved in the limit 
by metrics converging to a Schwarzschild metric (outside their 
respective horizons).

In the case that the $g$ is harmonically flat at infinity, a Schwarzschild
metric (outside the horizon) 
is contained in $[g]_S$.  More generally, given any asymptotically
flat manifold $(M^3,g)$, we can use $\real^3 \backslash  B_r(0)$ as a 
coordinate chart for the asymptotically flat end of $(M^3,g)$ which we are
interested in, where the metric $g_{ij}$ approaches $\delta_{ij}$
at infinity in this 
coordinate chart.
Then we can consider the conformal metric
\begin{equation}
   g_C = \left(1 + \frac{C}{|x|}\right)^4 g
\end{equation}
in this end.  In the limit as $C$ goes to infinity, 
the horizon will approach the coordinate sphere of radius $C$.  Then
outside this horizon in the limit as $C$ goes to infinity, the function  
$(1 + \frac{C}{|x|})$ will be close to a superharmonic function
on $(M^3,g)$ and the
metric $g_C$ will approach a Schwarzschild metric
(since the metric $g$ is approaching the standard metric on 
$\real^3$).  Hence, 
the Penrose quotient of $g_C$ will approach $(16 \pi)^{-1/2}$, which is the
Penrose quotient of a Schwarzschild metric.

As a final note, we prove that the first order o.d.e.~for $\{g_t\}$ defined 
in equations \ref{eqn:ODE1},\ref{eqn:ODE2}, \ref{eqn:ODE3}, and \ref{eqn:ODE4} 
is naturally defined in the sense that the rate of change of $g_t$ is a 
function only of $g_t$ and not of $g_0$ or $t$.  
To see this, given any solution $g_t = u_t(x)^4 g_0$ to equations 
\ref{eqn:ODE1},\ref{eqn:ODE2}, \ref{eqn:ODE3}, and \ref{eqn:ODE4},
choose any $s > 0$ and define 
$\bar{u}_t(x) = u_t(x)/u_s(x)$ so that
\begin{equation}\label{eqn:ODE1_bar}
   g_t = \bar{u}_t(x)^4 g_s 
\end{equation}
and $\bar{u}_s(x) \equiv 1$.  
Then  define $\bar{v}_t(x)$ such that
\begin{equation}\label{eqn:ODE3_bar}
\left\{
\begin{array}{r l l l}
\Delta_{g_s} \bar{v}_t(x) & \equiv & 0 & \mbox{ outside } \Sigma(t) \\
\bar{v}_t(x) & = & 0 & \mbox{ on } \Sigma(t) \\
\lim_{x \rightarrow \infty} \bar{v}_t(x) & = & -e^{-(t-s)} & \\
\end{array}
\right.
\end{equation}
and $\bar{v}_t(x) \equiv 0$ inside $\Sigma(t)$.  
Then what we want to show is 
\begin{equation}\label{eqn:ODE4_bar}
\bar{u}_t(x) = 1 + \int_s^t \bar{v}_r(x) dr
\end{equation}
To prove the above equation, we observe that
from equations \ref{eqn:identity}, \ref{eqn:ODE3_bar}, and \ref{eqn:ODE3} 
it follows that 
\begin{equation}
   v_t(x) = \bar{v}_t(x) \, u_s(x)
\end{equation}
since $\lim_{x \rightarrow \infty} u_s(x) = e^{-s}$.
Hence, since 
\begin{equation}
u_t(x) = u_s(x) + \int_s^t v_r(x) dr 
\end{equation}
by equation \ref{eqn:ODE4}, dividing through by $u_s(x)$ yields
equation \ref{eqn:ODE4_bar} as desired.
Thus, we see that 
the rate of change of $g_t(x)$ at 
$t=s$ is a function of $\bar{v}_s(x)$ which in turn
is just a function of $g_s(x)$ and the horizon $\Sigma(s)$.  Hence, to 
understand properties of the flow we need only analyze the behavior of the 
flow for $t$ close to zero, since any metric in the flow may be chosen to be
the base metric.  This point is used many times throughout the paper.

\section{An Example Solution to the Conformal Flow of Metrics}
\label{sec:Schwarzschild}

In this section we give the simplest example of a solution to the first order
o.d.e.~conformal flow of metrics defined by equations \ref{eqn:ODE1},
\ref{eqn:ODE2}, \ref{eqn:ODE3}, and \ref{eqn:ODE4}.  The initial metric in this
example is the three dimensional, space-like Schwarzschild metric which 
represents a single, non-rotating black hole in vacuum.  The
Schwarzschild metrics are 
also very natural from a geometric standpoint 
as well since they are spherically symmetric and have zero scalar 
curvature.

Since the flow does
not change the metric inside the horizon, we will define this metric to have
its horizon as a boundary, which is always allowable.  Then $(M^3,g_0)$ will
be defined to be isometric to 
$(\real^3-B_{m/2}(0), {\cal U}_0(x)^4 \delta_{ij})$, where
\begin{equation}
   {\cal U}_0(x) = 1 + \frac{m}{2r}
\end{equation}
where $r$ is the distance from the origin in $(\real^3, \delta_{ij})$ and 
$m$ is a positive constant equal to the mass of the black hole.

Next we define $(M^3,g_t)$ to be isometric to 
$(\real^3-B_{m/2}(0), {\cal U}_t(x)^4 \delta_{ij})$, where
\begin{equation}
   {\cal U}_t(x) = \left\{ \begin{array} {ll} 
e^{-t} + \frac{m}{2r} e^t, & \mbox{ for } r \ge \frac{m}{2} e^{2t}  \\ 
\sqrt{\frac{2m}{r}}, &\mbox{ for } r < \frac{m}{2} e^{2t}.  \end{array} 
\right.
\end{equation}
We note that on this metric the outermost horizon (and also the 
outermost minimal area enclosure of the original horizon) is the 
coordinate sphere given by $r = \frac{m}{2} e^{2t}$, so by 
equation \ref{eqn:ODE2} we define 
this horizon to be $\Sigma(t)$.  

Next we recall from equation \ref{eqn:ODE1} that $g_t = u_t(x)^4 g_0$.
Hence, $u_t(x) = {\cal U}_t(x)/ {\cal U}_0(x)$.  Furthermore,  
by equation \ref{eqn:ODE4} we must have $v_t(x) = \frac{d}{dt} u_t(x)$,
so 
\begin{equation}
   v_t(x) = \frac{1}{{\cal U}_0(x)}
\left\{ \begin{array} {ll} 
- e^{-t} + \frac{m}{2r} e^t, & \mbox{ for } r \ge \frac{m}{2} e^{2t}  \\ 
0, &\mbox{ for } r < \frac{m}{2} e^{2t}.  \end{array} 
\right.
\end{equation}
By equation \ref{eqn:identity}, $v_t(x)$ is harmonic on $(M^3,g_0)$
outside $\Sigma(t)$ since $a + b/r$ is harmonic in 
$(\real^3,\delta_{ij})$.
Then since $v_t(x)$ goes to $-e^{-t}$ at infinity, is continuous, 
and equals zero 
inside $\Sigma(t)$, it follows that equation \ref{eqn:ODE3} is 
satisfied.  Hence, $(M^3,g_t)$ is a solution to the first order conformal
flow of metrics defined by equations \ref{eqn:ODE1},
\ref{eqn:ODE2}, \ref{eqn:ODE3}, and \ref{eqn:ODE4}. 

This example is a good example to keep in mind when considering 
the main theorems of this paper.  For example, we notice that by
definition \ref{def:totalmass} the total mass $m(t)$ of $(M^3,g_t)$
equals $m$ and hence is nonincreasing as claimed in section 
\ref{sec:m(t)}, and the area $A(t)$
of the horizon $\Sigma(t)$ in $(M^3,g_t)$ is constant as claimed in
section \ref{sec:A(t)}.  Also, we see that the diameter of $\Sigma(t)$
is growing exponentially as claimed in theorem \ref{thm:diam}
and contains any given bounded set in a finite amount of time as 
claimed by theorem \ref{thm:bounded}.

Finally, we note that for all $t \ge 0$ in this example,  
$(M^3,g_t)$ is isometric to a 
Schwarzschild metric of total mass $m$ outside their respective 
horizons.  Hence, even though the metric is shrinking pointwise,
it is not changing at all outside its horizon, 
{\it after a reparametrization of the metric}.  It is in this sense that 
theorem \ref{thm:limit} states that no matter what the initial metric is, 
it eventually converges to a Schwarzschild metric outside its horizon.

\section{A Nonlinear Property of Superharmonic \\ Functions in $\real^3$}
\label{sec:shf}

In this appendix we present a nonlinear property of superharmonic functions 
in $\real^3$ which we needed in section \ref{sec:bounded}.  
However, this result
is of independent interest and, as is proven in \cite{BI}, directly implies 
the Riemannian Penrose inequality 
with suboptimal constant for manifolds which are 
conformal to $\real^3$.
Furthermore, the following theorem  
can be generalized to higher dimensions (for certain powers which 
unfortunately are not applicable to the Riemannian Penrose inequality
in higher dimensions) and is 
fully discussed and proven in \cite{BI}.  

\begin{theorem}\label{thm:shf}
There exists a constant $c>0$ such that if $u(x)$ is any positive, continuous, 
superharmonic function in $R^3$ satisfying   
\begin{equation}\label{eqn:shf1}
   \int_{S_r(0)} u(x)^4 \; dA \ge a
\end{equation}
for all $r > r_0$, then 
\begin{equation}\label{eqn:shf2}
   u(x) \ge c \, a^{1/4} |x|^{-1/2}
\end{equation}
for $|x| \ge r_0$.
\end{theorem}
{\it Discussion of proof.}
We refer the reader to \cite{BI}
(which is a joint work with Kevin Iga) 
for the details of the proof.  In that paper
we show that without loss of generality we may assume that $u(x)$ goes to zero
at infinity.  Then the next important step 
is a symmetrization argument to argue
that without loss of generality we may also assume 
that the support of $\Delta u$ is on the $x$-axis.  Then it follows that
\begin{equation}
   u(\vec{x}) = \int_0^\infty \frac{d\mu(t)}{|\vec{x}-(t,0,0)|}  
\end{equation}
where $\mu(t)$ is a positive measure on $[0, \infty)$.  The remainder
of the proof then 
involves converting inequality \ref{eqn:shf1} to a lower bound
on an integral expression of $d\mu(t)$ 
which is then used to prove inequality \ref{eqn:shf2}.  \qed

\section{Lower Bound for the Capacity of a Surface in $\real^3$ 
         with Bounded Willmore Functional}\label{sec:lower}
  
In this appendix we consider surfaces $\Sigma$ 
which are smooth, compact boundaries of open sets in $\real^3$ and which
have bounded Willmore functional.  That is, we will assume that 
\begin{equation}
 \int_\Sigma H^2 d\mu \le w ,
\end{equation}
where $H$ is the mean curvature (equal to the trace of the second fundamental
form) and $d\mu$ is the area form of $\Sigma$.

\vspace{.5in}
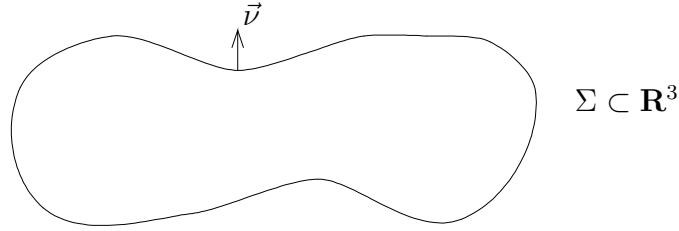
\begin{center}
\setlength{\unitlength}{0.0060in}
\begin{picture}(514,199)(0,-10)
\path(199,138)(199,173)(194,158)
\path(199,173)(204,158)
\path(309,168)	(306.451,167.604)
	(303.741,167.083)
	(300.881,166.445)
	(297.880,165.700)
	(294.750,164.857)
	(291.500,163.924)
	(288.141,162.911)
	(284.683,161.827)
	(281.137,160.681)
	(277.513,159.482)
	(273.821,158.239)
	(270.072,156.960)
	(266.275,155.656)
	(262.443,154.336)
	(258.584,153.007)
	(254.709,151.680)
	(250.828,150.363)
	(246.953,149.065)
	(243.092,147.796)
	(239.257,146.564)
	(235.458,145.379)
	(231.705,144.249)
	(228.009,143.184)
	(224.380,142.192)
	(220.828,141.284)
	(217.364,140.467)
	(213.998,139.751)
	(210.740,139.145)
	(207.601,138.657)
	(204.591,138.298)
	(201.721,138.076)
	(199.000,138.000)

\path(199,138)	(196.291,138.097)
	(193.454,138.378)
	(190.498,138.832)
	(187.430,139.443)
	(184.258,140.201)
	(180.991,141.091)
	(177.636,142.101)
	(174.202,143.217)
	(170.697,144.426)
	(167.129,145.716)
	(163.505,147.074)
	(159.835,148.485)
	(156.125,149.938)
	(152.385,151.418)
	(148.622,152.914)
	(144.844,154.412)
	(141.060,155.899)
	(137.277,157.361)
	(133.504,158.787)
	(129.748,160.162)
	(126.018,161.474)
	(122.322,162.710)
	(118.668,163.857)
	(115.064,164.901)
	(111.518,165.829)
	(108.038,166.630)
	(104.632,167.288)
	(101.309,167.792)
	(98.076,168.128)
	(94.941,168.284)
	(91.913,168.245)
	(89.000,168.000)

\path(89,168)	(84.205,167.298)
	(81.669,166.846)
	(79.054,166.327)
	(76.367,165.740)
	(73.617,165.084)
	(70.814,164.360)
	(67.966,163.567)
	(65.081,162.705)
	(62.170,161.773)
	(59.239,160.771)
	(56.299,159.699)
	(53.358,158.557)
	(50.425,157.343)
	(47.509,156.059)
	(44.618,154.703)
	(41.761,153.274)
	(38.947,151.774)
	(36.185,150.201)
	(33.484,148.555)
	(30.852,146.836)
	(28.298,145.044)
	(23.461,141.237)
	(19.042,137.131)
	(15.112,132.725)
	(11.741,128.016)
	(10.287,125.546)
	(9.000,123.000)

\path(9,123)	(7.665,119.943)
	(6.456,116.796)
	(5.375,113.567)
	(4.420,110.263)
	(3.593,106.893)
	(2.892,103.463)
	(2.318,99.981)
	(1.871,96.454)
	(1.550,92.891)
	(1.356,89.298)
	(1.289,85.683)
	(1.348,82.054)
	(1.534,78.417)
	(1.845,74.782)
	(2.284,71.154)
	(2.848,67.542)
	(3.539,63.953)
	(4.355,60.395)
	(5.298,56.875)
	(6.367,53.400)
	(7.562,49.979)
	(8.882,46.618)
	(10.329,43.325)
	(11.901,40.109)
	(13.599,36.975)
	(15.423,33.932)
	(17.372,30.987)
	(19.447,28.148)
	(21.647,25.423)
	(23.973,22.818)
	(26.424,20.341)
	(29.000,18.000)

\path(29,18)	(32.301,15.382)
	(35.773,13.048)
	(39.404,10.985)
	(43.181,9.180)
	(47.092,7.621)
	(51.126,6.296)
	(55.269,5.191)
	(59.510,4.296)
	(63.836,3.596)
	(68.236,3.081)
	(72.698,2.736)
	(77.208,2.551)
	(81.755,2.511)
	(86.327,2.606)
	(90.912,2.822)
	(95.497,3.147)
	(100.071,3.569)
	(104.620,4.074)
	(109.134,4.651)
	(113.599,5.288)
	(118.005,5.971)
	(122.337,6.688)
	(126.585,7.427)
	(130.736,8.175)
	(134.778,8.920)
	(138.699,9.649)
	(142.487,10.351)
	(146.129,11.012)
	(149.614,11.619)
	(152.928,12.162)
	(156.061,12.626)
	(159.000,13.000)

\path(159,13)	(161.569,13.353)
	(164.297,13.838)
	(167.172,14.445)
	(170.185,15.164)
	(173.325,15.987)
	(176.582,16.904)
	(179.947,17.905)
	(183.408,18.981)
	(186.956,20.123)
	(190.580,21.322)
	(194.270,22.568)
	(198.016,23.851)
	(201.809,25.163)
	(205.636,26.493)
	(209.490,27.833)
	(213.358,29.174)
	(217.231,30.505)
	(221.100,31.818)
	(224.953,33.102)
	(228.780,34.350)
	(232.571,35.551)
	(236.317,36.695)
	(240.006,37.775)
	(243.629,38.779)
	(247.176,39.700)
	(250.635,40.527)
	(253.998,41.251)
	(257.254,41.863)
	(260.392,42.353)
	(263.402,42.712)
	(266.275,42.931)
	(269.000,43.000)

\path(269,43)	(272.071,42.843)
	(275.250,42.407)
	(278.533,41.715)
	(281.914,40.787)
	(285.388,39.646)
	(288.948,38.313)
	(292.590,36.811)
	(296.308,35.161)
	(300.096,33.384)
	(303.950,31.504)
	(307.863,29.542)
	(311.829,27.519)
	(315.844,25.457)
	(319.902,23.379)
	(323.998,21.306)
	(328.125,19.259)
	(332.279,17.262)
	(336.453,15.335)
	(340.643,13.500)
	(344.843,11.780)
	(349.047,10.196)
	(353.249,8.769)
	(357.446,7.523)
	(361.629,6.478)
	(365.796,5.657)
	(369.939,5.081)
	(374.053,4.772)
	(378.133,4.752)
	(382.173,5.043)
	(386.168,5.667)
	(390.112,6.645)
	(394.000,8.000)

\path(394,8)	(397.762,9.627)
	(401.450,11.451)
	(405.061,13.465)
	(408.591,15.660)
	(412.035,18.028)
	(415.391,20.561)
	(418.655,23.251)
	(421.821,26.090)
	(424.888,29.070)
	(427.851,32.183)
	(430.705,35.420)
	(433.449,38.775)
	(436.076,42.238)
	(438.585,45.802)
	(440.970,49.459)
	(443.228,53.200)
	(445.356,57.018)
	(447.349,60.904)
	(449.204,64.851)
	(450.917,68.850)
	(452.484,72.894)
	(453.901,76.974)
	(455.165,81.082)
	(456.271,85.211)
	(457.216,89.351)
	(457.996,93.496)
	(458.608,97.637)
	(459.047,101.765)
	(459.309,105.874)
	(459.391,109.955)
	(459.290,114.000)
	(459.000,118.000)

\path(459,118)	(458.389,121.719)
	(457.313,125.375)
	(455.817,128.954)
	(453.948,132.445)
	(451.750,135.836)
	(449.269,139.115)
	(446.552,142.268)
	(443.644,145.284)
	(440.591,148.150)
	(437.438,150.855)
	(434.232,153.386)
	(431.018,155.730)
	(427.842,157.876)
	(424.750,159.811)
	(421.788,161.523)
	(419.000,163.000)

\path(419,163)	(416.356,164.142)
	(413.383,165.095)
	(410.127,165.874)
	(406.637,166.497)
	(402.961,166.982)
	(399.148,167.344)
	(395.245,167.602)
	(391.299,167.771)
	(387.361,167.870)
	(383.476,167.914)
	(379.694,167.922)
	(376.062,167.910)
	(372.629,167.894)
	(369.442,167.893)
	(366.550,167.922)
	(364.000,168.000)

\path(364,168)	(361.454,168.112)
	(358.554,168.233)
	(355.352,168.357)
	(351.896,168.479)
	(348.238,168.594)
	(344.428,168.698)
	(340.515,168.785)
	(336.551,168.851)
	(332.585,168.889)
	(328.668,168.897)
	(324.849,168.867)
	(321.180,168.797)
	(317.710,168.680)
	(314.490,168.511)
	(311.570,168.286)
	(309.000,168.000)

\put(194,178){\makebox(0,0)[lb]{\smash{{{\SetFigFont{12}{14.4}{rm}
$\vec{\nu}$ }}}}}
\put(484,103){\makebox(0,0)[lb]{\smash{{{\SetFigFont{12}{14.4}{rm}
$\Sigma \subset \real^3$ }}}}}
\end{picture}
\end{center}
Define the potential function $f$ of 
$\Sigma$ to be the constant function $1$ inside $\Sigma$ and 
the function satisfying
\begin{equation}\label{eqn:f}
\left\{ \begin{array}{rll}
f(x) = & 1 & \mbox{on } \Sigma \\
\Delta f \equiv & 0 & \mbox{outside } \Sigma \\
f \rightarrow & 0 & \mbox{at infinity}
\end{array}\right.
\end{equation}
outside $\Sigma$.
Expanding in terms of spherical harmonics, we see that
\begin{equation}
 f(x) = \frac{a}{r} + O \left(\frac{1}{r^2}\right) .
\end{equation}
We define $a$ to be the capacity 
of both $f(x)$ and $\Sigma$.  Also, 
since $f$ is nonnegative, $a$ is always nonnegative.  Furthermore, we 
notice that 
\begin{equation}
 a = \lim_{r \rightarrow \infty} - \frac{1}{4\pi} \int_{S_r(0)}
\frac{df}{dr}d\mu \,\,\, ,
\end{equation}
where $S_r(0)$ is the sphere of radius $r$ centered around
zero.  But since $f(x)$ is harmonic, it follows from the divergence theorem
that if we perform the above integral over any surface enclosing $\Sigma$, 
we will get the same result.  Thus, 
\begin{equation}
a = - \frac{1}{4\pi} \int_{\Sigma}
\frac{df}{d\vec{\nu}}d\mu \,\,\, ,
\end{equation}
where $\vec{\nu}$ is the outward pointing unit normal vector to $\Sigma$.
The goal of this appendix is to find a lower bound for $a$ 
in terms of the area of $\Sigma$ and the Willmore
bound $w$.

Since $f$ is harmonic off of $\Sigma$, 
the support of the distribution $\Delta f$ is on $\Sigma$.  In fact
\begin{equation}
 \Delta f (\phi) = \int_{\real^3} f \Delta \phi 
                   = \int_{\Sigma} \phi \frac{df}{d\vec{\nu}} d\mu 
\end{equation}
where $\frac{df}{d\vec{\nu}}$ is defined to be the outward directional
derivative of $f$ (which does not equal the inward directional derivative
which is zero).  Then since $f$ equals $\Delta f$ convolved with the Green's
function $- \frac{1}{4\pi |x|}$, we get that
\begin{equation}\label{eqn:f(y)}
f(y) = - \frac{1}{4\pi} \int_{\Sigma} 
          \frac{df}{d\vec{\nu}} \cdot \frac{1}{|y-x|} dx ,
\end{equation}
where $dx$ is the area form of $\Sigma$ with respect to the variable $x$.  
In a moment $dy$ will be the area form of $\Sigma$ with respect to the variable
$y$.
Then using the fact that $f(y) = 1$ on $\Sigma$, we find that
\begin{equation}\label{eqn:avg}
|\Sigma| = \int_{\Sigma} f(y)dy 
         = - \frac{1}{4\pi} \int_{\Sigma}  
             \frac{df}{d\vec{\nu}} \cdot 
             \left( \int_{\Sigma} \frac{1}{|y-x|} dy \right) dx .
\end{equation}

The next step is to find an upper bound for $\int_{\Sigma} \frac{1}{|y-x|} dy$.
To do this, we need the following lemma which follows from  
equation 16.34 in \cite{GT} (when $R$ is chosen to go to infinity).  We note
that the definition of mean curvature in \cite{GT} is half that of ours.
\begin{lemma}
Given a surface $\Sigma$ which is a smooth, compact boundary of
an open set in $\real^3$, 
\begin{equation}
 |\Sigma \cap B_r(x)| \le \left(\frac34 \int_{\Sigma} H^2 d\mu \right) r^2 
\end{equation}
for all $r > 0$ and all $x \in \Sigma$.
\end{lemma}
Hence, in our case, 
\begin{equation}\label{eqn:area_growth}
|\Sigma \cap B_r(x)| \le (\frac34 w)r^2 ,
\end{equation}
for any $x \in \Sigma$.  Then since $\frac{1}{|y-x|}$ is maximized when
$x$ and $y$ are closest together, we get that
$\int_{\Sigma} \frac{1}{|y-x|} dy$ is maximized subject to the constraint
in equation \ref{eqn:area_growth} when $|\Sigma \cap B_r(x)| = (\frac34 w)r^2$,
which gives us  
\begin{equation}
\int_{\Sigma} \frac{1}{|y-x|} dy \le 
\int_0^{\sqrt{\frac{4|\Sigma|}{3w}}} \frac{1}{r} \cdot 2(\frac34 w)r dr
= (3w)^{1/2} |\Sigma|^{1/2} .
\end{equation}
Thus, from equation \ref{eqn:avg},
\begin{equation}
|\Sigma| \le - \frac{1}{4\pi} \int_{\Sigma}  
             \frac{df}{d\vec{\nu}} dx \cdot (3w)^{1/2} |\Sigma|^{1/2} ,
\end{equation}
so that we have 
\begin{equation}
 a = - \frac{1}{4\pi} \int_{\Sigma}
\frac{df}{d\vec{\nu}}d\mu \ge (3w)^{-1/2} |\Sigma|^{1/2} .
\end{equation}
Thus, in summary, we have the following theorem.
\begin{theorem}\label{thm:capacity_of_sigma}
Let $\Sigma$ be a smooth, compact boundary of an open set in $\real^3$
with 
$\int_\Sigma H^2 d\mu \le w$.  Let $f(x)$ be the harmonic function
equal to one on $\Sigma$ and going to zero at infinity.  Then the capacity
\begin{equation}
a \equiv - \frac{1}{4\pi} \int_{\Sigma}
\frac{df}{d\vec{\nu}}d\mu \ge (3w)^{-1/2} |\Sigma|^{1/2} .
\end{equation}
\end{theorem}
We need this lower bound on the capacity of $\Sigma(t)$ in section
\ref{sec:asymptotic}.

\section{Regularity of the horizons $\Sigma^\epsilon(t)$}
\label{sec:regularity}

In this appendix we compute upper bounds on the $C^{k,\alpha}$ ``norms''
(see definition \ref{def:k,alpha,S}) of the surfaces $\Sigma^\epsilon(t)$,
for some $\alpha \in (0,1)$.  
We will also compute upper bounds on the $C^{k,\alpha}$ norms
of the metrics $g^\epsilon_t$ outside $\Sigma^\epsilon(t)$.
These bounds are  
independent of $\epsilon$, allowing us to conclude in  
section \ref{sec:existence} that 
that the horizons $\Sigma(t)$ are smooth and the metrics $g_t$ are smooth
outside $\Sigma(t)$. 

The critical step will be to achieve a uniform (independent of $\epsilon$)
bound on the $C^{1,\alpha}$ norms
of the surfaces $\Sigma^\epsilon(t)$ in the coordinate charts, 
for $t \in [0,T]$.   Then by 
equations \ref{eqn:ODE3_ep} and 
\ref{eqn:ODE4_ep}, we will also have a uniform bound on the $C^{1,\alpha}$ 
norms of the metrics $g^\epsilon_t = u^\epsilon_t(x)^4 g_0$ outside
$\Sigma^\epsilon(t)$.  Schauder 
estimates applied
to the minimal surface equation will then give us a uniform bound on the 
$C^{2,\alpha}$ norms of the surfaces $\Sigma^\epsilon(t)$.  Repeating this 
bootstrapping argument yields the desired uniform bounds on the 
$C^{k,\alpha}$ norms
of the surfaces $\Sigma^\epsilon(t)$ and the metric 
$g^\epsilon_t$ outside
$\Sigma^\epsilon(t)$. 

We note that it follows from the definition of $\Sigma^\epsilon(t)$ 
given in equation \ref{eqn:ODE2_ep} (and from equation \ref{eqn:ODE3_ep}) 
that $\Sigma^\epsilon(t)$ globally 
minimizes area among surfaces in ${\cal S}$ (defined in section 
\ref{sec:definitions}) which contain the original horizon $\Sigma_0$.
Then as in section \ref{sec:green}, 
let us define $M^3_{\Sigma_0}$ to be the closed region of $M^3$ which
is outside (or on) $\Sigma_0$, and define $(\bar{M}^3_{\Sigma_0},\bar{g}_0)$
to be two distinct copies of $(M^3_{\Sigma_0},g_0)$ 
identified along $\Sigma_0$.
Hence, $(\bar{M}^3_{\Sigma_0},\bar{g}_0)$ 
has a reflection
symmetry which keeps $\Sigma_0$ fixed.

With out loss of generality, let us now replace $(M^3,g_0)$ with 
$(\bar{M}^3_{\Sigma_0},\bar{g}_0)$ for the remainder of this section.
We can do this since the portion of $(M^3,g_0)$ which is inside 
$\Sigma_0$ does not affect the conformal flow of metrics 
defined in equations \ref{eqn:ODE1_ep},
\ref{eqn:ODE2_ep}, \ref{eqn:ODE3_ep}, and \ref{eqn:ODE4_ep}.
It is true that this operation means that the new $(M^3,g_0)$ will not be
smooth along $\Sigma_0$, but this turns out not to be important.  Furthermore,
the big advantage (which we leave to the reader to check) 
is that now $\Sigma^\epsilon(t)$ globally minimizes area among 
{\it all} surfaces in ${\cal S}$, which is a fact we will use in the
next lemma.  However, because of this construction, all of the constants we
will be defining will depend on $\Sigma_0$, although this is not a problem.

We will continue with the same notation as defined in section 
\ref{sec:existence}.  Since the initial manifold $(M^3,g_0)$ is smooth 
(except possibly on $\Sigma_0$ because of the above construction, which we
will ignore) and 
harmonically flat at infinity, 
it may be covered by a finite number of smooth coordinate
charts $\{{\cal C}_i\}$, 
where each harmonically flat end (of which there can only be a finite
number) has $\real^3 \backslash B_1(0)$ as an asymptotically flat 
coordinate chart and all of the 
other coordinate charts are different copies of $B_1(0) \subset \real^3$.  
For the rest of this appendix 
we will work inside the coordinate charts $\{{\cal C}_i\}$ which have the 
standard $\real^3$ metric.

\begin{definition}
Let $\mbox{dis}_i(x,S)$ be the infimum (which could equal infinity) of the 
lengths of all of the paths in ${\cal C}_i$ between the point $x$ and 
the set $S$ with respect to the coordinate chart metric in ${\cal C}_i$.
Also, let $\mbox{dis}_{g_0}(x,S)$ be the infimum
of the 
lengths of all of the paths in $M^3$ between the point $x$ and 
the set $S$ with respect to the metric $g_0$.
\end{definition}

\begin{definition}
Let $X^\epsilon(t)$ be the three dimensional open 
region inside $\Sigma^\epsilon(t)$ so that
\begin{equation}\label{eqn:X^ep}
   \Sigma^\epsilon(t) = \partial X^\epsilon(t).
\end{equation}
\end{definition}
The region $X^\epsilon(t)$ always exists since $\Sigma^\epsilon(t) \in 
{\cal S}$.

\begin{lemma}\label{lem:c1r2}
There exists a constant $c_1 > 0$ depending only on $T$, $\Sigma_0$, $g_0$, 
and the choice of 
coordinate charts $\{{\cal C}_i\}$ such that for all $i$ and all
${x_0} \in \Sigma^\epsilon(t) \cap {\cal C}_i$,
\begin{equation}\label{eqn:area_r^2}
   |X^\epsilon(t) \cap S_r({x_0})|_{{\cal C}_i} \ge c_1 r^2
\end{equation}
for all $r \in (0, \mbox{dis}_i({x_0},\partial {\cal C}_i))$ and $t \in [0,T]$.
\end{lemma}
{\it Proof.} 
Let $S_r({x_0})$ and $B_r({x_0})$ respectively 
be the coordinate 2-sphere and closed 
coordinate 3-ball of radii $r$ centered at 
${x_0} \in \Sigma^\epsilon(t) \cap {\cal C}_i$.
Define
\begin{equation}
   A(r) = |\Sigma^\epsilon(t) \cap B_r({x_0})|_{{\cal C}_i},
\end{equation}
\begin{equation}
   L(r) = |\Sigma^\epsilon(t) \cap S_r({x_0})|_{{\cal C}_i},
\end{equation}
and 
\begin{equation}
   \bar{A}(r) = |X^\epsilon(t) \cap S_r({x_0})|_{{\cal C}_i}.
\end{equation}

Let us also define $f(a) = l$ for $a \in [0, 4\pi]$, where 
$l$ is the minimum length required to enclose a region of area $a$ in
the unit sphere $S^2$.
Hence, for small $a$, $f(a) \approx \sqrt{4\pi a}$.

Finally, we let $\gamma_1$ be the infimum of the smallest eigenvalue
and let $\gamma_2$ be the supremum of the largest eigenvalue of the 
metric $g_0$ over every point in all of the coordinate charts $\{{\cal C}_i\}$,
and then define $\gamma = \gamma_1/\gamma_2$.
Since the coordinate charts are smooth and $g_0$ approaches the standard metric
$\delta_{ij}$ in the ends of the noncompact coordinate charts, $\gamma > 0$.
Then since 
\begin{equation}\label{eqn:uepbound}
1 \ge u^\epsilon_t(x) \ge (1-\epsilon)^{\lb \frac{t}{\epsilon} \rb} 
\approx e^{-t}
\end{equation}
by the maximum principle and equations \ref{eqn:vvv} and \ref{eqn:v=-u}, it 
then follows that the corresponding ratio of eigenvalues for the metric 
$g^\epsilon_t$ in the coordinate charts
$\{{\cal C}_i\}$ is at least $\gamma_T$, where 
\begin{equation}
   \gamma_T = (1-\epsilon)^{\lb \frac{4T}{\epsilon} \rb} \gamma > 0,
\end{equation}
$t \in [0,T]$, and as usual we are requiring $t$ to be an integral multiple
of $\epsilon$.  In the limit as $\epsilon$ goes to zero, we note that
$\gamma_T$ approaches $e^{-4T} \gamma$, which is positive.

Now we are ready to begin the actual proof of the lemma.  
Since $\Sigma^\epsilon(t)$ is globally area 
minimizing among surfaces in ${\cal S}$ with respect to the metric 
$g^\epsilon_t$,
we know that it has area in $(M^3,g^\epsilon_t)$ less than or equal to 
the areas of the two comparison surfaces
$\partial (X^\epsilon(t) \backslash B_r({x_0}))$ and 
$\partial (X^\epsilon(t) \cup B_r({x_0}))$.  Putting this into terms of the 
above definitions, this yields
\begin{equation}\label{eqn:gamma_area}
   \gamma_T A(r) \le \bar{A}(r) \le 4\pi r^2 - \gamma_T A(r).
\end{equation}
Hence, it follows from the above equation that 
\begin{equation}
   f\left(\frac{\bar{A}(r)}{r^2}\right) \ge f\left(\frac{\gamma_T A(r)}{r^2}
   \right)
\end{equation}
since $f(a) = f(4\pi - a)$ and is monotone increasing from $0$ to $2\pi$. 
Furthermore, it follows from the definition of $f$ that
\begin{equation}
   L(r) \ge f\left(\frac{\bar{A}(r)}{r^2}\right) r.
\end{equation}
Hence, since from multivariable calculus we have that $A'(r) \ge L(r)$, we
deduce that
\begin{equation}\label{eqn:aprimebound}
   A'(r) \ge f\left(\frac{\gamma_T A(r)}{r^2}\right) r.
\end{equation}
It is then straightforward to show that inequalities \ref{eqn:gamma_area}
and \ref{eqn:aprimebound}
imply inequality \ref{eqn:area_r^2} with
$c_1 = 4 \gamma_T^2 / \pi$, proving the lemma.  \qed

\begin{corollary}\label{cor:w_Holder}
Let $w(x)$ be any nonnegative harmonic function in $(M^3,g_0)$ defined 
outside $\Sigma^\epsilon(t)$ which equals zero on $\Sigma^\epsilon(t)$.
Define $w(x)$ to be identically zero 
inside $\Sigma^\epsilon(t)$.  Then there exists a constant $c_2 \in (0,1)$
depending only on $T$, $\Sigma_0$, $g_0$, and $\{{\cal C}_i\}$ such that for 
all $i$ and all ${x_0} \in \Sigma^\epsilon(t) \cap {\cal C}_i$, 
\begin{equation}
   \sup_{S_{r/2}({x_0})} w(y) \le c_2 \sup_{S_{r}({x_0})} w(y) 
\end{equation}
for all $r \in (0, \mbox{dis}_i({x_0},\partial {\cal C}_i))$ and $t \in [0,T]$.
\end{corollary}
{\it Proof.}
The function $w(x)$ is subharmonic in $(M^3,g_0)$, 
and hence is bounded above by the harmonic function $h(x)$ 
defined in the coordinate ball $B_r({x_0})$ with Dirichlet boundary data
$h(x) = \sup_{S_{r}({x_0})} w(y)$ on $S_r({x_0}) \backslash X^\epsilon(t)$ and 
$h(x) = 0$ on $S_r({x_0}) \cap X^\epsilon(t)$.  The corollary then follows by 
estimating $h(x)$ on $S_{r/2}({x_0})$ using the Poisson kernel and 
lemma \ref{lem:c1r2}.  \qed

\begin{corollary}
There exist constants $c_3$, $c_4$, and $\beta \in (0,1)$ depending only on 
$T$, $\Sigma_0$, $g_0$, and $\{{\cal C}_i\}$ such that 
\begin{equation}\label{eqn:ineq1}
|v^\epsilon_t(x)| \le c_3 \mbox{dis}_{g_0}(x, X^\epsilon(t))^\beta
\end{equation}
and for all $y \in {\cal C}_i$ with 
$\mbox{dis}_{g_0}(x,y) \le \mbox{dis}_{g_0}(x, X^\epsilon(t))$
\begin{equation}\label{eqn:ineq2}
|v^\epsilon_t(x) - v^\epsilon_t(y)| \le c_4 \mbox{dis}_{g_0}(x,y)^\beta
\end{equation}
for all $t \in [0,T]$
\end{corollary}
{\it Proof.}
First we note that since distances with respect to the coordinate chart metrics
and $g_0$ are within a bounded factor of each other, we only need to 
prove each of the above inequalities in each coordinate chart.

Inequality \ref{eqn:ineq1} follows from the definition of $v^\epsilon_t(x)$
given in equation \ref{eqn:ODE3_ep} and from recursively applying 
corollary \ref{cor:w_Holder}.  Inequality \ref{eqn:ineq1} then implies the
interior gradient estimate
\begin{equation}
   |\nabla v^\epsilon_t(x)| \le c_5 \, \mbox{dis}_{g_0}
   (x, X^\epsilon(t))^{\beta-1},
\end{equation}
which, when integrated along the straight path 
connecting $x$ and $y$, implies inequality \ref{eqn:ineq2}.  \qed

\begin{corollary}\label{cor:c_6}
Let $c_6 = 2 \max(c_3,c_4)$.  Then 
\begin{equation}\label{eqn:v_t_Holder}
   |v^\epsilon_t(x) - v^\epsilon_t(y)| \le c_6 \mbox{dis}_{g_0}(x,y)^\beta. 
\end{equation}
for all $x$ and $y$ and for $t \in [0,T]$.
\end{corollary}
{\it Proof.} 
Let $r_x = \mbox{dis}_{g_0}(x,X^\epsilon(t))$ and 
$r_y = \mbox{dis}_{g_0}(y,X^\epsilon(t))$, and let $B$ now denote geodesic
balls in $(M^3,g_0)$.  Then we consider two cases.

{\it Case 1:}  Suppose $B_{r_x}(x) \cap B_{r_y}(y) = \varnothing$.  Then
inequality \ref{eqn:v_t_Holder} follows from inequality \ref{eqn:ineq1} and
the triangle inequality.

{\it Case 2:}  Suppose $B_{r_x}(x) \cap B_{r_y}(y) \neq \varnothing$.  Then
choose $z \in B_{r_x}(x) \cap B_{r_y}(y)$ which is on the length minimizing
geodesic 
connecting the points $x$ and $y$.  Then inequality \ref{eqn:v_t_Holder} 
follows from the triangle inequality and 
inequality \ref{eqn:ineq2} applied to the points $x$ and $z$ and
to the points $y$ and $z$.  \qed

\begin{definition}\label{def:k,alpha}
Let $w(x)$ be any $C^k$ function defined on $(M^3,g_0)$.  
Then we define the following norm and seminorms (denoted by brackets), all
of which depend on our choice of coordinate charts $\{{\cal C}_i\}$,
$\alpha \in (0,1)$, and $k = 0, 1, 2, ...$
\begin{equation} 
[w]_{k;\Omega} = 
\sup_i \;\sup_{x \in {\cal C}_i}\;\sup_{|\gamma| = k} |D^\gamma w(x)|, 
\end{equation}
\begin{equation} 
[w]_{k,\alpha;\Omega} = \sup_{i} \; \sup_{x \ne y} \; \sup_{|\gamma| = k} 
             \left\{ \frac{|D^{\gamma}w(x)-D^{\gamma}w(y)|}{|x-y|^\alpha} 
             \;\;|\;\;  x,y \in {\cal C}_i \right\} ,
\end{equation}
and
\begin{equation}
||w||_{C^{k,\alpha}(\Omega)} = [w]_{k,\alpha;\Omega} + 
                               \sum_{j=1}^k [w]_{j;\Omega} ,
\end{equation}
where we also require that $x,y \in \Omega \subset M^3$ in the above equations.
\end{definition}
Hence, from corollary \ref{cor:c_6} 
and equation \ref{eqn:ODE4_ep}, it follows that 
\begin{equation}\label{eqn:uHolder}
   ||u^\epsilon_t||_{C^{0,\beta}(M^3)} \le c_7 
\end{equation}
for $t \in [0,T]$, where $c_7$ depends only on 
$T$, $\Sigma_0$, $g_0$, and $\{{\cal C}_i\}$.

\begin{definition}\label{def:k,alpha,S}
Let $S$ be any smooth surface in $(M^3,g_0)$ which is the boundary of 
a region.  Let $\eta$ be any vector field defined on all of $(M^3,g_0)$
such that on $S$ it equals the outward pointing unit normal vector of $S$.
Let $\eta = (\eta_1, \eta_2, \eta_3)$ be the pull back of $\eta$ on each 
coordinate chart.  We note that definition \ref{def:k,alpha} can be used 
for vector valued functions as well as real valued functions.  
Then abusing notation slightly (since the following is
not a norm), we define
\begin{equation}
||S||_{C^{k,\alpha}} = \inf_{\eta} ||\eta||_{C^{k-1,\alpha}(S)}
\end{equation}
for $\alpha \in (0,1)$ and $k = 1, 2, ...$
\end{definition}

The next part of the proof is to use inequality \ref{eqn:uHolder} and
the fact that $\Sigma^\epsilon(t)$ minimizes area in 
$(M^3,g^\epsilon_t = u^\epsilon_t(x)^4 g_0)$ to conclude that 
the $\Sigma^\epsilon(t)$ are uniformly
$C^{1,\beta/4}$ surfaces.  By this we mean that we will find an upper bound
on $||\Sigma^\epsilon(t)||_{C^{1,\beta/4}}$ which is independent of $\epsilon$.
Conveniently, the main theorem we need,
theorem \ref{thm:DeGiorgi}, was essentially 
proved already by De Giorgi \cite{DG} to understand the regularity of 
codimension one 
minimal surfaces in $\real^n$, and is summarized in \cite{MM}.  
However, the application of this theorem to this setting is quite 
interesting, so we summarize the arguments below.

\begin{definition}
Let $X$ and $Y$ be regions in $(M^3,g_0)$ (and hence also in the coordinate
charts $\{{\cal C}_i\}$) with smooth boundaries of finite area. 
Then we define
\[
   \psi_{g_0}(X,x,r) = |\partial X \cap B_r(x)|_{g_0} -
   \inf\{|\partial Y \cap B_r(x)|_{g_0} \;|\; Y=X 
   \mbox{ outside } B_r(x) \}
\]
for $x \in {\cal C}_i$ and 
$r \in (0, \mbox{dis}_i(x,\partial {\cal C}_i))$, where $B_r(x)$ is the
closed ball of radius $r$ in the coordinate chart ${\cal C}_i$.
\end{definition}

Hence, the function $\psi_{g_0}(X,x,r)$ can be thought of 
as the excess area of $\partial X$ in $(M^3,g_0)$ in the 
coordinate ball $B_r(x)$.  We note that the above definition for 
$\psi_{g_0}(X,x,r)$ in the coordinate chart ${\cal C}_i$ 
equals $\psi(X,x,r)$ defined in \cite{MM} in the 
special case that $g_0$ equals the coordinate chart metric.

\begin{definition}
Let $X$ be a region in $(M^3,g_0)$ which has a smooth boundary of finite area.
Then following \cite{MM}, we define $X$ to be
{\bf $(K,\lambda)_{g_0}$-minimal} 
in $\{{\cal C}_i\}$ if and only if 
\begin{equation}
   \psi_{g_0}(X,x,r) \le K r^{2+\lambda}
\end{equation}
for all $i$, $x \in {\cal C}_i$, and 
$r \in (0, \mbox{dis}_i(x,\partial {\cal C}_i))$.
\end{definition}

Again, we note that a region which is 
$(K,\lambda)_{g_0}$-minimal in the coordinate chart ${\cal C}_i$
as defined above is $(K,\lambda)$-minimal as 
defined in \cite{MM} in the special case that 
$g_0$ equals the coordinate chart metric.
The usefulness of the above definition can be seen in the next lemma
and theorem. 

\begin{lemma}\label{lem:uniform}
Suppose $X$ is a region in $(M^3,g_0)$ which has a smooth boundary of
finite area and which is $(K,\lambda)_{g_0}$-minimal in $\{{\cal C}_i\}$.
We note that the metric $g_0$ and the coordinate charts $\{{\cal C}_i\}$ 
are assumed to be smooth.
Then  for all $\epsilon > 0$, there exists a $\delta > 0$ 
(depending only on $K$, $\lambda$, $g_0$, and $\{{\cal C}_i\}$)
such that for all $i$ and $x_0 \in \partial X \cap {\cal C}_i$,
\begin{equation}
   \omega_i(X,x_0,r) \equiv |D\chi_X|(B_r(x_0)) - |D\chi_X(B_r(x_0))| < 
   \epsilon \, r^2
\end{equation}
for all $r \in (0,\min(\delta, \mbox{dis}_i(x_0,\partial {\cal C}_i))]$.
\end{lemma}
{\it Proof.}  
We note that we are adopting the notation of \cite{MM}, so that $\chi_X$
is the characteristic function of the region $X$ and $D\chi_X$ is the 
distributional derivative (with respect to the coordinate chart) 
of that characteristic function.  
Hence, $\omega_i$ equals zero for regions with flat
boundaries in the coordinate chart ${\cal C}_i$, 
and in general can be thought of as a way of measuring how 
far a boundary is from being flat inside the coordinate ball $B_r(x_0)$.

For convenience, we will also assume that $g_0$ equals the standard $\real^3$
metric in each of the finite number of coordinate charts.  
Then since $g_0$ is smooth and the proof uses a blow up argument, it is 
easy to adapt the proof we give here to the general case.

We proceed with a proof by contradiction.  Suppose that for some $\epsilon>0$,
there existed a counterexample region $X_\delta$ for all $\delta>0$ such that
\begin{equation}
   \omega_i(X_\delta, x_0, r_\delta) \ge \epsilon \, r_\delta^2
\end{equation}
for some $i$ and $x_0 \in \partial X_\delta \cap {\cal C}_i$ and for some 
$r_\delta \in (0,\min(\delta, \mbox{dis}_i(x_0,\partial {\cal C}_i)]$.
We may as well think of each region $X_\delta$ being in the same $\real^3$, 
and for convenience we translate each region by $-x_0$ (which depends on 
$\delta$) so that $0 \in \partial X_\delta$.  
Let $\bar{X}_\delta$ be $X_\delta$ rescaled by a factor of $1/r_\delta$.
Then we have that 
\begin{equation}\label{eqn:omegalb}
   \omega(\bar{X}_\delta, 0, 1) \ge \epsilon,
\end{equation}
where $\omega$ is defined for regions in $\real^3$.  Furthermore,
$\bar{X}_\delta$ is $(K r_\delta^\lambda, \lambda)$-minimal in 
$B_{1/r_\delta}(0) \subset \real^3$.

Let $\bar{X}$ be a limit region of $\{\bar{X}_\delta\}$, in the sense that
the characteristic function of $\bar{X}$ is the $L^1$ limit of the 
characteristic functions of $\{\bar{X}_{\delta_i}\}$, for some sequence
of $\{\delta_i\}$ converging to zero.  
The fact that such a limit
region exists is proven on p. 70 of \cite{MM}, and relies on the fact that
the areas of the boundaries of $\{\bar{X}_\delta\}$ are uniformly bounded
above, which follows from theorem 1 of section 2.5.3 of \cite{MM}.
Furthermore, by 
the lower semicontinuity of the area functional as defined in \cite{MM},
it follows that $\bar{X}$ has minimal boundary. 
Since the only minimizing boundaries in $R^3$ are planes, 
$\partial \bar{X}$ must be a plane
going through the origin, which we may as well assume is the $x$-$y$ plane
in $\real^3$ after a suitable rotation.

Let $\pi$ be the projection of $\real^3$ to the $x$-$y$ axis, and 
let $\vec\nu$ be the outward pointing normal vector of 
$\partial \bar{X}_{\delta_i}$.  Then
\begin{equation}
\lim_{\delta_i \rightarrow 0} 
|\int_{\partial \bar{X}_{\delta_i} \cap B_1(0)} \nu_z|
= \lim_{\delta_i \rightarrow 0} |\pi(\partial \bar{X}_{\delta_i} \cap B_1(0))|
= \lim_{\delta_i \rightarrow 0} |\partial \bar{X}_{\delta_i} \cap B_1(0)|
\end{equation}
since $\{\bar{X}_\delta\} \rightarrow \bar{X}$ and the areas converge as well
since each $\bar{X}_\delta$ is $(K r_\delta^\lambda, \lambda)$-minimal.
The above equation then implies that $\lim_{\delta_i \rightarrow 0} 
\omega(\bar{X}_{\delta_i}, 0, 1) = 0$, contradicting inequality 
\ref{eqn:omegalb} and proving the lemma.  \qed

\begin{theorem}\label{thm:DeGiorgi}
Suppose $X$ is a region in $M^3$ which has a smooth boundary of 
finite area and which is $(K,\lambda)_{g_0}$-minimal in $\{{\cal C}_i\}$. 
Then $\partial X$ is a 
$C^{1,\lambda/4}$ surface and 
\begin{equation}
|| \partial X ||_{C^{1,\lambda/4}} \le \bar{k}, 
\end{equation}
where $\bar{k}$ depends only on
$K$, $\lambda$, $g_0$, and $\{{\cal C}_i\}$.
\end{theorem}
{\it Proof.}  We restrict our attention to each coordinate chart
${\cal C}_i$ one at a time, and first consider the case that $g_0$ equals
the coordinate chart metric.  In this case, the theorem follows directly
from the proof of theorem 1 of section 2.5.4 of \cite{MM} and lemma 
\ref{lem:uniform}.  Then it is then a somewhat long but straightforward 
task to adapt the relevant theorems of \cite{MM} 
to verify that theorem 1 of section 2.5.4 of \cite{MM} is still true 
if we replace the standard metric on $\real^3$ with 
any smooth fixed metric $g_0$.  \qed

\begin{lemma}\label{lem:minimal}
The region $X^\epsilon(t)$, which we recall is defined to be the region
inside $\Sigma^\epsilon(t)$, is a 
$(K,\beta)_{g_0}$-minimal set in $\{{\cal C}_i\}$, for $t \in [0,T]$,
where $K$ depends only on  
$T$, $\Sigma_0$, $g_0$, 
and $\{{\cal C}_i\}$).
\end{lemma}
{\it Proof.}  We need to estimate   
\begin{equation}\begin{array}{ll} 
   \psi_{g_0}(X^\epsilon(t),x,r) & = 
   |\partial X^\epsilon(t) \cap B_r(x)|_{g_0} \\ &
   - \inf\{|\partial Y \cap B_r(x)|_{g_0} \;|\; Y=X^\epsilon(t) 
   \mbox{ outside } B_r(x) \} \end{array}
\end{equation}
from above.  Well,
since $\partial X^\epsilon(t)$ has minimal area in $(M^3,g^\epsilon_t)$, 
we have that
\begin{equation}
   |\partial Y \cap B_r(x)|_{g^\epsilon_t} \ge 
   |\partial X^\epsilon(t) \cap B_r(x)|_{g^\epsilon_t},
\end{equation}
which gives us
\begin{equation}
   |\partial Y \cap B_r(x)|_{g_0} \ge  
   |\partial X^\epsilon(t) \cap B_r(x)|_{g_0}
   \left(\frac{u_{min}}{u_{max}}\right)^4
\end{equation}
where $u_{min}$ and $u_{max}$ are respectively the minimum and maximum values
of $u^\epsilon_t(x)$ on the closed coordinate ball $B_r(x)$.  Hence,
\begin{equation}
\psi_{g_0}(X^\epsilon(t),x,r) \le |\partial X^\epsilon(t) \cap B_r(x)|_{g_0}
\left[1 - \left(\frac{u_{min}}{u_{max}}\right)^4 \right].
\end{equation}
The lemma then follows from equations \ref{eqn:uepbound}, 
\ref{eqn:gamma_area}, and \ref{eqn:uHolder}.  \qed

\begin{corollary}\label{cor:DeGiorgi}
The surface $\Sigma^\epsilon(t)$ is a 
$C^{1,\beta/4}$ surface, and  
\begin{equation}
   ||\Sigma^\epsilon(t)||_{C^{1,\beta/4}} \le \bar{c}_1,
\end{equation} 
for $t \in [0,T]$, 
where $\bar{c}_1$ depends only on
$T$, $\Sigma_0$, $g_0$, and $\{{\cal C}_i\}$.
\end{corollary}
{\it Proof.}  Follows directly from theorem \ref{thm:DeGiorgi}
and lemma \ref{lem:minimal}.  \qed

The remainder of the arguments presented in this appendix are mostly standard
applications of \cite{GT}.
For convenience, we now let $\alpha = \beta/4$.
\begin{lemma}\label{lem:bootstrap1}
For $k \ge 1$, 
\begin{equation}
||\Sigma^\epsilon(s)||_{C^{k,\alpha}} \le \bar{c}_k 
\mbox{ for all } s\in[0,T]
\;\;\Rightarrow\;\;
||u^\epsilon_t(x)||_{C^{k,\alpha}(M^3 \backslash X^\epsilon(t))}
\le \bar{\bar{c}}_k
\end{equation}
for $t \in [0,T]$,
where $\bar{\bar{c}}_k$ depends only on $\bar{c}_k$, 
$T$, $\Sigma_0$, $g_0$, and $\{{\cal C}_i\}$.
\end{lemma} 
{\it Proof.}  By the hypothesis of the lemma, 
the definition of $v^\epsilon_t(x)$ given in equation \ref{eqn:ODE3_ep},
and standard theorems found in \cite{GT} and \cite{Widman} (for $k=1$), 
we get an upper bound on 
$||v^\epsilon_s(x)||_{C^{k,\alpha}(M^3 \backslash X^\epsilon(t))}$.  
The lemma then follows from equation \ref{eqn:ODE4_ep}.
\qed 

\begin{lemma}\label{lem:bootstrap2}
For $k \ge 1$ and $t \in [0,T]$ an integer multiple of $\epsilon$,
\begin{equation}
||u^\epsilon_t(x)||_{C^{k,\alpha}(M^3 \backslash X^\epsilon(t))}
\le \bar{\bar{c}}_k,
\;\;\Rightarrow\;\;
||\Sigma^\epsilon(t)||_{C^{k+1,\alpha}} \le \bar{c}_{k+1}
\end{equation}
where $\bar{c}_{k+1}$ depends only on $\bar{\bar{c}}_k$, 
$T$, $\Sigma_0$, $g_0$, and $\{{\cal C}_i\}$.
\end{lemma}
{\it Proof.}  
Since $\Sigma^\epsilon(t)$ minimizes area in $(M^3,g_t)$ and since it 
can be viewed as a graph of a function of two variables over a uniformly large
domain by corollary \ref{cor:DeGiorgi}, we can apply 
Schauder estimates to the minimal surface equation to prove the lemma.

Let the graph of $x_3 = f(x_1,x_2)$ represent the minimal surface 
$\Sigma^\epsilon(t)$ in the coordinate chart ${\cal C}_i \subset \real^3$
with metric $g_{ij}(x_1,x_2,x_3)$ given by the metric $g^\epsilon_t
= u^\epsilon_t(x)^4 g_0$.  
Then we note that the minimal surface equation in this setting is
\begin{equation}\label{eqn:MSE}
a^{ij}(x) D_{ij}f(x) 
= p(x)
\end{equation}
where
\begin{equation}
a^{ij} = G^{ij} - \frac{G^{i\alpha}N_{\alpha}G^{j\beta}N_{\beta}}
{N^tGN},
\end{equation}
\begin{equation}
p = \frac12 N^t G_3 N
+ \left(\frac{D_i(G^{\alpha\beta})N_{\alpha}N_{\beta}G^{i\lambda}}{2N^tGN}
- D_i(G^{i\lambda})\right) N_\lambda ,
\end{equation}
\begin{equation}
G^{\alpha\beta} = g_{\alpha+1,\beta+1}g_{\alpha+2,\beta+2} - g_{\alpha+1,\beta+2}g_{\alpha+2,\beta+1}
\end{equation}
is the determinant of the $\alpha\beta$ minor matrix 
and the subscript addition in the above equation is modulo 3,  
\begin{equation}
   G_3^{\alpha\beta} = \frac{\partial}{\partial x_3} G^{\alpha\beta}, 
\end{equation}
and 
\begin{equation}
N = (D_1f, D_2f, -1),
\end{equation}
where $x = (x_1,x_2)$ and naturally $x_3 = f(x_1,x_2)$ in the above equations,
and we are using the convention that Latin indices are summed from 1 to 2 while
Greek indices are summed from 1 to 3.

We also note that since
\begin{equation}
   \Lambda_1 |v|^2 \le g_{\alpha\beta}v^\alpha v^\beta \le \Lambda_2 |v|^2 
\end{equation}
for some positive $\Lambda_1$ and $\Lambda_2$ which depend on $t$ 
by equation \ref{eqn:uepbound}, and since 
\begin{equation}
   |Df| \le B
\end{equation}
over a uniformly sized domain by corollary \ref{cor:DeGiorgi}, then it can
be shown that 
\begin{equation}\label{eqn:ellipticity_condition}
   a^{ij}v_iv_j \ge \lambda |v|^2,
\end{equation}
where
\begin{equation}
   \lambda = \frac{\Lambda_1^8}{4\Lambda_2^6 (1+B^2)^2}.
\end{equation}

Hence, since the $a_{ij}$ are uniformly positive definite and 
$a^{ij}(x)$ and $p(x)$ involve only first order derivatives of $f(x)$
and the metric $g_{ij}$, the lemma follows from corollary \ref{cor:DeGiorgi}
and bootstrapping with Schauder estimates applied
to equation \ref{eqn:MSE}.  \qed

Then combining corollary \ref{cor:DeGiorgi} with lemmas \ref{lem:bootstrap1}
and \ref{lem:bootstrap2}, we get the main result 
in this appendix.

\begin{corollary}\label{cor:regularity}
For $t \in [0,T]$, for some $\alpha \in (0,1/4)$, and for $k \ge 1$, 
we have uniform bounds on the $C^{k,\alpha}$ norms of $u^\epsilon_t(x)$
outside $\Sigma^\epsilon(t)$ and on the surfaces $\Sigma^\epsilon(t)$, 
as defined in definitions \ref{def:k,alpha} and \ref{def:k,alpha,S}.
Furthermore, these bounds depend 
only on $k$, $T$, $\Sigma_0$, 
$g_0$, and $\{{\cal C}_i\}$, and are independent of $\epsilon$.
\end{corollary}


\end{document}